\documentclass{siamonline171218} 

\usepackage[utf8]{inputenc}

\let\equation\relax
\let\equation*\relax
\usepackage{amsmath,amssymb,ifthen,fancyhdr,lastpage,mathrsfs,amsfonts,color,listings,dsfont}
\allowdisplaybreaks

\usepackage{graphicx}
\usepackage{xcolor}
\usepackage{subfig}
\usepackage{algorithm}
\usepackage{algorithmic}
\usepackage[colorinlistoftodos]{todonotes}

\usepackage{hyperref}
\hypersetup{
    colorlinks,
    linkcolor={blue!100!black},
    citecolor={blue!100!black},
    urlcolor={blue!100!black}
}

\usepackage{etoolbox}

\makeatletter
\DeclareRobustCommand{\barredF}{%
  \mathord{\vphantom{\mathfrak F}\mathpalette\@barredF\relax}%
}

\newcommand{\@barredF}[2]{%
  \ooalign{%
    \hidewidth\@barredFbar#1\hidewidth\cr
    $\m@th#1\mathfrak f$\cr
  }%
}

\newcommand{\@barredFbar}[1]{%
  \check@mathfonts
  \ifx#1\displaystyle
    \fontsize{\f@size}{\z@}%
    \def\@barredFbarkern{0.3}%
  \else
    \ifx#1\textstyle
      \fontsize{\f@size}{\z@}
      \def\@barredFbarkern{0.3}%
    \else
      \ifx#1\scriptstyle
        \fontsize{\sf@size}{\z@}
        \def\@barredFbarkern{0.4}%
      \else
        \fontsize{\ssf@size}{\z@}
        \def\@barredFbarkern{0.47}%
      \fi
    \fi
  \fi
  \usefont{OT1}{cmr}{m}{n}%
  \kern-\@barredFbarkern em 
  \raisebox{-.8ex}{\symbol{'26}}%
}
\makeatother

\newcommand{\sC}{\mathscr{C}}

\newcommand{\cD}{\mathcal{D}}

\newcommand{\cF}{\mathcal{F}}

\newcommand{\cM}{\mathcal{M}}
\newcommand{\cN}{\mathcal{N}}
\newcommand{\cO}{\mathcal{O}}

\newcommand{\cR}{\mathcal{R}}
\newcommand{\cS}{\mathcal{S}}

\newcommand{\mN}{\mathbb{N}}

\newcommand{\mR}{\mathbb{R}}

\newcommand{\mD}{\mathbb{D}}

\newcommand{\mV}{\mathbb{V}}

\newcommand{\mS}{\mathbb{S}}

 \newcommand{\be}{\ensuremath{\boldsymbol{e}}}

 \newcommand{\bx}{\ensuremath{\boldsymbol{x}}}

 \newcommand{\Bf}{\ensuremath{\boldsymbol{f}}}
 \newcommand{\Bg}{\ensuremath{\boldsymbol{g}}}
  \newcommand{\Bp}{\ensuremath{\boldsymbol{p}}}
  \newcommand{\BP}{\ensuremath{\boldsymbol{P}}}

  \newcommand{\BT}{\ensuremath{\boldsymbol{T}}}
  \newcommand{\bU}{\ensuremath{\boldsymbol{U}}}
  \newcommand{\bxi}{\ensuremath{\boldsymbol{\xi}}}
 
 \newcommand{\bn}{\ensuremath{\boldsymbol{n}}}
 \newcommand{\bw}{\ensuremath{\boldsymbol{w}}}
  
  \newcommand{\bu}{\ensuremath{\boldsymbol{u}}}
  \newcommand{\btu}{\ensuremath{\boldsymbol{\tilde u}}}
  \newcommand{\btheta}{\ensuremath{\boldsymbol{\theta}}}
  \newcommand{\bvarphi}{\ensuremath{\boldsymbol{\varphi}}}
  \newcommand{\bepsilon}{\ensuremath{\boldsymbol{\epsilon}}}
  \newcommand{\bnabla}{\ensuremath{\boldsymbol{\nabla}}}
  \newcommand{\src}{\ensuremath{\boldsymbol{s}}}

\newcommand{\Text}[1]{\text{\textnormal{#1}}}

\newcommand{\I}{\Text{i}}

\newcommand{\FA}{\;\; \forall \;}
\newcommand{\MTEXT}[1]{\;\;\;\;\;\text{#1}\;\;\;\;\;}

\newcommand{\D}{\Text{d}}

\newcommand{\closure}[2][3]{%
  {}\mkern#1mu\overline{\mkern-#1mu#2}}

\newcommand{\adj}[1]{#1^\ast}
\newcommand{\cc}[1]{\overline{#1}}

\newcommand{\transp}[1]{#1^\textup{T}}
\newcommand{\Nproj}{{N_{\textup{proj}}}}
\newcommand{\CInf}[1]{{\mathscr{C}^\infty ( #1 ) }}
\newcommand{\CInfO}{\CInf{\Omega}}
\newcommand{\CcInf}[1]{{\mathscr{C}^\infty_{\textup c} ( #1 ) }}
\newcommand{\CcInfO}{\CcInf{\Omega}}

\newcommand{\mRinf}{\mR \cup \{ \infty \}}

\newcommand{\FDGrad}{\nabla_{\textup{FD}}}
\newcommand{\F}{\barredF}
\newcommand{\FAlt}{\lower0.9ex\hbox{$\mathchar'26$}\mkern-8mu \mathfrak{f}}
\newcommand{\mF}{m_{\mathfrak{f}}}
\newcommand{\Mpx}{{m_{\textup{proj}}}}
\newcommand{\Mdata}{{m_{\textup{data}}}}
\newcommand{\Nvx}{n}
\newcommand{\nablaProj}{\nabla\!_\Proj}
\newcommand{\tot}{{\textup{tot}}}

\newcommand{\new}{\textup{new}}
\newcommand{\REF}{\textup{ref}}
\newcommand{\STOP}{\textup{stop}}

\newcommand{\obs}{{\textup{obs}}}

\newcommand{\grad}{\textup{grad}}

\newcommand{\IE}{i.e.\ }

\newcommand{\cSAbstr}{\tilde{\cS}}
\newcommand{\ProjAbstr}{{\tilde P}}
\newcommand{\isoBProjAbstr}{B_{\textup{iso}}}
\newcommand{\isoProj}{\Proj_{\textup{iso}}}
\newcommand{\backProjAbstr}{\ProjAbstr^{\textup{B}}}
\newcommand{\backProj}{\Proj^{\textup{B}}}
\newcommand{\backParProj}{\ParProj^{\textup{B}}}
\newcommand{\backCone}{\Cone^{\textup{B}}}
\newcommand{\Proj}{P}
\newcommand{\ParProj}{\mathscr{P}}
\newcommand{\isoParProj}{\ParProj_{\textup{iso}}}
\newcommand{\Cone}{\mathscr{D}}
\newcommand{\isoCone}{\Cone_{\textup{iso}}}

\newcommand{\kl}{\textup{KL}}

\newcommand{\uproj}{u_\Proj}
\newcommand{\uparproj}{u_\ParProj}
\newcommand{\tuproj}{{\tilde u_\Proj}}
\newcommand{\wproj}{{w_\Proj}}
\newcommand{\uprojj}[1]{u_{#1}}
\newcommand{\tuprojj}[1]{{\tilde u_{#1}}}
\newcommand{\wprojj}[1]{{w_{#1}}}

\newcommand{\DP}{{\mD\!_{\Proj}}}
\newcommand{\DPj}[1]{{\mD\!_{{#1}}}}
\newcommand{\DR}{{\mD\!_{\ParProj}}}
\newcommand{\DC}{{\mD\!_{\Cone}}}
\newcommand{\DetDom}{{L^2(\DP)}}
\newcommand{\DetDomR}{{L^2(\DR)}}

\newcommand{\DetDomj}[1]{{L^2(\DPj{#1})}}

\newcommand{\OmegaX}{{L^\ParProj_{\bx_\perp}}}
\newcommand{\OmegaPhi}{{L^\Cone_{\bvarphi}}}

\newcommand{\secref}[1]{$\S$\ref{#1}}
\newcommand{\ipref}[1]{Inverse Problem~\ref{#1}}
\newcommand{\nnl}{\nonumber \\}

\newcommand{\norm}[1]{\| #1 \|} 
\newcommand{\abs}[1]{| #1 |}
\newcommand{\parens}[1]{(#1 )}

\newcommand{\ip}[2]{\langle #1, #2\rangle} 

\newcommand{\Norm}[1]{\left\| #1 \right\|} 

\newcommand{\Parens}[1]{\left(#1 \right)}

\newcommand{\bnorm}[1]{\big\| #1 \big\|} 

\newcommand{\bparens}[1]{\big(#1 \big)}

\newcommand{\bip}[2]{\big\langle #1, #2\big\rangle} 

\newcommand{\Bparens}[1]{\Big(#1 \Big)}

\newcommand{\bbabs}[1]{\bigg| #1 \bigg|}
\newcommand{\bbparens}[1]{\bigg(#1 \bigg)}

\newcommand{\optspace}{\vspace{.5em}}

\DeclareMathOperator*{\argmin}{argmin}

\DeclareMathOperator*{\supp}{supp}
\DeclareMathOperator*{\Kern}{kern}
\DeclareMathOperator*{\Range}{range}

\DeclareMathOperator*{\identity}{id}

\DeclareMathOperator*{\imag}{Im}

\DeclareMathOperator*{\prox}{prox}

\newcommand{\enum}[1]{$\boldsymbol{(#1)}$}

\newcommand{\qedhere}{}

\numberwithin{theorem}{section}
\numberwithin{equation}{section}

\newsiamthm{MyIP}{Inverse Problem}
\newsiamthm{assump}{Assumption}
\newsiamthm{problem}{Problem}
\newsiamthm{convention}{Convention}
\newsiamremark{remark}{Remark}
\newsiamthm{scheme}{Scheme}

\begin{document}

\title{Generalized SART Methods for Tomographic Imaging\thanks{Submitted to the editors of Inverse Problems on January 7, 2019.\funding{This work was funded Deutsche Forschungsgemeinschaft DFG through Project C02
of SFB 755 - Nanoscale Photonic Imaging.}}}

\author{Simon Maretzke\thanks{Insitute for Numerical and Applied Mathematics, University of Goettingen, Lotzestrasse 16-18, 37083 G\"ottingen, Germany (\email{simon.maretzke@googlemail.com}).}} 

\headers{Generalized SART Methods for Tomographic Imaging}{Simon Maretzke}



\maketitle

\begin{abstract}
 Nowadays, the field computed tomography (CT) encompasses a large variety of settings, ranging from nanoscale to meter-sized objects imaged by different kinds of radiation in various acquisition modes. This experimental diversity challenges the flexibility of tomographic reconstruction methods.
 Kaczmarz-type methods, which exploit the natural block-structure of tomographic inverse problems, are a promising candidate to provide the required versatility in a computationally efficient manner.  
 In the present work, it is shown that indeed a surprisingly general class of tomographic Kaczmarz-iterations may be efficiently evaluated via  computational schemes of a similar structure as updates of the so-called simultaneous algebraic reconstruction technique (SART). 
 This enables regularized reconstructions with non-trivial image-formation models as well as non-quadratic or even non-convex data-fidelity terms at low computational costs. 
 Moreover, the proposed generalized SART schemes are equally applicable in parallel- and cone-beam settings and regardless of the choice of tomographic incident directions.
 Their potential is illustrated by outlining applications in several non-standard tomographic settings, including polychromatic CT and X-ray phase contrast tomography.
\end{abstract}

\begin{keywords}
tomographic imaging, Kaczmarz methods, 
variational reconstruction methods, cone-beam tomography, robust reconstruction, phase contrast tomography
\end{keywords}

\begin{AMS}
65R10, 65R32, 92C55, 94A08
\end{AMS}



%
%
%

\sloppy
\section{Introduction} \label{S:Intro}

Since the pioneering works of Cormack \cite{Cormack1963CT,Cormack1963CTII} and Hounsfield \cite{Hounsfield1973CT}, the field of computed tomography (CT) has broadened considerably. While classical CT based on the partial attenuation of X-rays by matter continues to be a principal workhorse of medical diagnosis, several other applications have emerged over the past decades: for instance, state-of-the-art transmission-electron-microscopes (TEM) may resolve unknown structures in three dimensions down to sub-nanometer resolutions by acquiring a series TEM-images under different incident directions of the electron-beam \cite{MidgleyWeyland2002_ETomoContrastMechanisms,Oektem2015MathETomo}. Moreover, the advent of {coherent} X-ray sources 
has enabled \emph{phase contrast} techniques and thereby extended the scope of X-ray tomography to quasi nonabsorbing micro- and nanoscale objects such as single biological cells or viruses \cite{Cloetens1999holotomography,Krenkel2014BCAandCTF,Guizar2015_QuantitativeROITomo,EkebergEtAl2015_Virus3DCDI,Bartels2015}.
Other CT applications rely on gamma-rays or even cosmic myons to image strongly attenuating objects such as oil pipelines \cite{Johansen2015GammaTomo} or ancient pyramids \cite{MorishimaEtAl2017PyramidMyonCT}.

All of the above settings come with their own peculiarities: In electron tomography for instance, tomographic views may be typically only acquired for a few incident directions in a limited range due to radiation damage and geometrical restrictions. In phase contrast tomography, the acquired data is given by  {diffraction patterns}, that relate to the actual projections of the object's refractivity in a highly non-trivial manner. Moreover, even the seemingly standard application of medical CT involves a generally nonlinear inverse problem. Disregarding this nonlinearity may lead to severe artifacts associated with so-called \emph{beam-hardening} or \emph{photon-starvation} effects \cite{BarrettEtAl2004ArtifactsInCT,Mori2013PhotonStarvation}. 

To address this variety of different settings, there exist for once a large number of different preprocessing methods that aim to make the data ``ready'' for reconstruction by standard algorithms such as filtered backprojection (FBP) or its approximate analogue for cone-beam tomography, the FDK-algorithm \cite{FeldkampDavisKress1984_FDKAlgorithm}. While this enables fast reconstruction, the applied pre-corrections are often heuristic and limited in their effectivity. For this reason, variational reconstruction methods have received increased attention recently, see e.g.\ \cite{SidkyEtAl2008CBTomoWithTV,Beck2009Teboulle_FISTA,YangEtAl2010_ROIByGenTV,Setzer2011ADMMForImageProcessing,SidkyEtAl2012TomoWithChambollePock,Kostenko2013_AllAtOncePCTWithTV,Burger2016VariationalImageProcessing}. These start from an image-formation model for the considered setting, supplement it with additional \emph{a priori} knowledge on the unknown object and expected data-errors and reconstruct by minimizing a cost functional that incorporates all of this information. While this approach is much more flexible in accounting for the peculiarities of a specific imaging modality, it suffers from computational complexity: typically, the optimization has to be carried out iteratively where each of the many iterations requires $\cO(N^4)$ arithmetic operations -- as much as a full FBP- or FDK-reconstruction. In current high-resolution CT applications with $N\gtrsim 10^3$ sampling-points along each dimension, this constitutes a major bottleneck.

Kaczmarz-type- or block-iterative methods like the (simultaneous) algebraic reconstruction technique (S)ART \cite{Gordon1970ART,AndersenKak1984SART} may provide a compromise between the flexibility of variational methods and the favorable complexity of FBP. 
By fitting only small data-\emph{blocks} in each iteration, 
Kaczmarz-updates can be computed at lower computational costs than bulk-iterations on the complete data. SART, for example, updates the reconstructed object by matching its projection to the measured data under one tomographic incident direction per step. Moreover, Kaczmarz-methods are often observed to exhibit fast \emph{semi-convergence} \cite{Natterer,ElfvingHansen2014SemiconvKaczmarz}, arriving at accurate reconstructions already after $\cO(1)$ fitting-cycles over the tomographic data set. 
If the iterations are sufficiently cheap to compute, this enables image-recovery at overall computational costs comparable to FBP or FDK.
Recent variants of Kaczmarz-type schemes may also incorporate advanced priors such as total-variation-penalties via interlacing gradient-descent- or proximal iterations \cite{Bertsekas2011_IncrProxMethods,AndersenHansen2014_genARTProx,DefriseEtAl2011_TVSurrogateSART} or block-structured   primal-dual-methods \cite{ChambolleEhrhardtEtAl2017StochasticPDHG}.

In the present work, it is shown that a very general class of regularized Kaczmarz-iterations (also called ``Tikhonov-Kaczmarz'' \cite{CezaroEtAl2011_TikhonovKaczmarz,KindermannLeitao2014_genARTConvergence} or ``incremental proximal'' iterations \cite{Bertsekas2011_IncrProxMethods2,AndersenHansen2014_genARTProx} by other authors) for tomographic problems may be evaluated via efficient computational schemes of a similar structure as classical SART-steps. The approach is therefore named \emph{generalized~SART} (GenSART). The underlying idea is that -- as long as the object is fitted to exactly \emph{one} tomographic view per iteration -- the computed updates will typically be uniform along the direction of the tomographic rays, i.e.\ a \emph{back-projection} of some increment in the lower-dimensional projection-space. 
This reasoning holds true for both parallel- and cone-beam settings, a large variety of data-fidelity functionals, complicated image-formation models such as diffraction operators and different regularization terms. The potential of GenSART is demonstrated by outlining applications in Poisson-noise-adapted and outlier-robust tomographic reconstruction, X-ray phase contrast tomography as well as for a polychromatic CT model. 

This manuscript is organized as follows: \secref{S:Background} introduces the basic tomographic imaging model and gives some background on Kaczmarz-type reconstruction methods. In \secref{S:SARTPrinciple}, the generalized SART principle is presented in an abstract setting, which is shown to be applicable to the considered Kaczmarz-iterations in \secref{S:AdmPenalties}. \secref{S:Applications} outlines several applications of the proposed GenSART-schemes, for some of which numerical examples are presented in \secref{S:NumExamples}. 


\optspace
\section{Background} \label{S:Background}

\subsection{Tomographic imaging model} \label{SS:TomoModel}

\begin{figure}[hbt!]
 \centering
 \includegraphics[width=.75\textwidth]{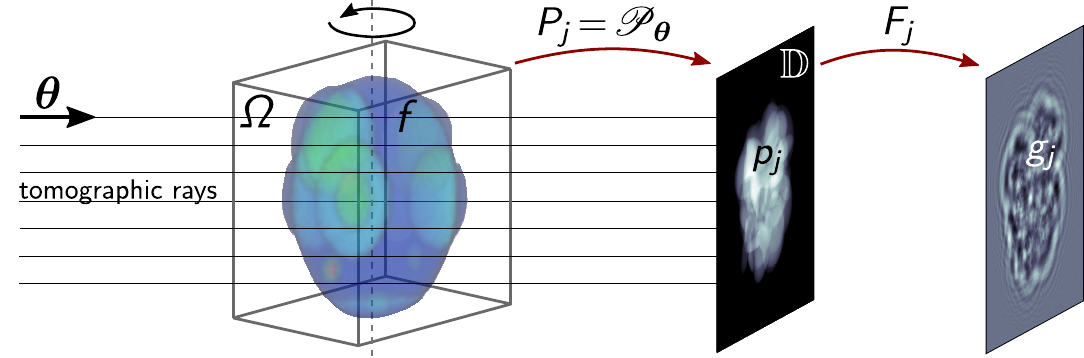} \vspace{.5em} \\
 \includegraphics[width=.75\textwidth]{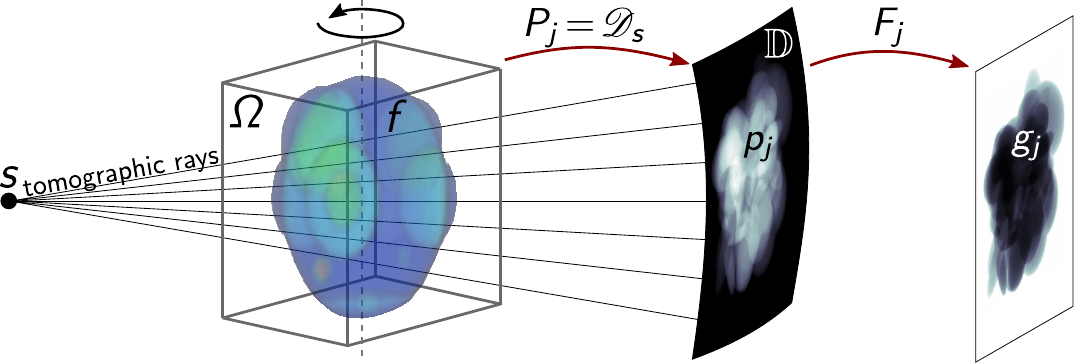}
 \caption{Sketch of the considered tomographic model: the measured data $g_j$ is given by a finite number of tomographic projections $p_j = \Proj_j(f)$ of an unknown 3D-object $f$, mapped under additional image-formation operators $F_j$. Top: parallel-beam example with diffractive image-formation (phase contrast). Bottom: cone-beam setup with absorption-contrast as in conventional CT.}
 \label{fig:sketch}
\end{figure}

We consider general tomographic inverse problems, for which the dependence of the data $g_{\tot}$ from the sought object $f$ can be modeled as
\begin{align}
 g_{\tot} = \begin{pmatrix}
               g_{1} \\
		\vdots 		\\
	       g_{\Nproj} 
             \end{pmatrix}
	       = \begin{pmatrix}
               F_1(\Proj_1 (f)) \\
		\vdots 		\\
	       F_{\Nproj}(\Proj_{\Nproj} (f))
             \end{pmatrix} = F_{\tot} ( \Proj_{\tot}(f) ) \label{eq:TomoModel}
\end{align}
$\Proj_j$ are parallel-beam- or cone-beam projectors acting on a 3D-object density $f \in L^2(\Omega)$ ($L^2(\Omega) := \{f: \mR^3 \to \mR : \norm{f}_{L^2} < \infty, \; \supp( f ) \subset \Omega\}$, $\norm{h}_{L^2}^2 := \int |h(x)| ^2 \, \D x$) with support $\supp(f)$ in a bounded open domain $\Omega \subset \mR^3$. The setting is sketched in \cref{fig:sketch}. 
Note that \eqref{eq:TomoModel} constitutes a \emph{semi-continuous} tomography model in that $f$ and the $g_{1}, \ldots, g_{\Nproj}$ are modeled as continuous images, but the number of tomographic views $\Nproj \in \mN$ is finite.

In a \emph{parallel-beam} tomography setup, each \smash{$\Proj_j = \ParProj_{\btheta_{(j)}}$} maps $f$ onto its line integrals along a certain incident direction $\btheta \in \mS^2 := \{\bx \in \mR^3: |\bx| = 1\}$: 
\begin{align}
 \ParProj_{\btheta} (f)( \bx_\perp ) := \int_{\mR}  f(x \bn_x + y \bn_y  + z \btheta) \, \D z \FA \bx_\perp = (x,y) \in \mR^2, 
 \label{eq:DefRadon}
\end{align}
 where $\bn_x, \bn_y \in \mS^2$ are chosen such that $\{\btheta, \bn_x, \bn_y \}$ forms an orthonormal basis of $\mR^3$. 
 
 In a \emph{cone-beam} tomography setup modeled by \smash{$\Proj_j = \Cone_{\src_{(j)}}$}, a single projection is given by line integrals along rays emanating from a point-source  $\src \in \mR^3  \setminus \closure \Omega$ along all directions:
\begin{align}
 \Cone_{\src} (f) ( \bvarphi )  := \int_0^\infty f(\src + t \bvarphi) \; \D t \FA \bvarphi \in \mS^2. \label{eq:DefCone}
\end{align}

The maps $F_j$ in \eqref{eq:TomoModel} denote (possibly nonlinear) \emph{image-formation operators}, modeling the relation between the raw projections $p_j = \Proj_j(f)$ and the data $g_j$ that is actually detected. For instance, the choice $F_j(p_j) = I_0 \cdot \exp( - p_j ) $ models absorption-based X-ray tomography with monochromatic illumination of known intensity distribution $I_0$. More complicated operators can be derived to model polychromatic CT or phase contrast tomography, for example. 

\subsection{Inverse problem and a priori constraints} \label{SS:InvProbConstraints}

This work is concerned with 
reconstructing the an unknown object $f$ from tomographic data $g_{\tot}$ modeled by \eqref{eq:TomoModel}, i.e.\ in solving the following inverse problem:
\vspace{.5em}
\begin{MyIP}[Tomographic reconstruction]\label{ip1}
 For a given domain $\Omega \subset \mR^3$, parallel- or cone-beam projectors $\Proj_{1}, \ldots, \Proj_{\Nproj}$ and image-formation operators $F_{1}, \ldots, F_{\Nproj}$, reconstruct an object  $f \in L^2(\Omega)$ from noisy data
 \begin{align} 
g_{j}^\obs =  F_{j}(\Proj_{j} (f)) + \bepsilon_{j} \MTEXT{for} j = 1, \ldots, \Nproj. 
 \end{align}
\end{MyIP}
\vspace{.5em}

As \ipref{ip1} is typically ill-posed even in the trivial case where the $F_j$ are just given by the identity, one typically aims to exploit available \emph{a priori} knowledge on the unknown density $f$ to be reconstructed. Common a priori constraints include:
\optspace
\begin{itemize}
  \item \emph{Support constraints:} $f$ is known to be identically zero outside $\Omega \subset \mR^3$. Note that such constraints are already incorporated in the tomographic imaging model from \secref{SS:TomoModel}.
\optspace
 \item \emph{Box constraints:} the admissible values of $f$ may be bounded from below and/or above. For example, it is often physically justified to assume \emph{non-negativity}, $f \geq 0$. \optspace
 \item \emph{Regularity constraints:} A certain smoothness of $f$ is assumed in the reconstruction, for example by imposing that the total variation seminorm $\int_{\Omega} |\nabla f| \, \D x$ is small.
\end{itemize}
\optspace
In addition to incorporating such constraints, one often seeks to account for an (approximately) known statistics of the data errors $\bepsilon$.

\subsection{Reconstruction methods} \label{SS:RecMethComplexity}

\paragraph{Variational methods:}  The widely used filtered back-projection- (FBP) and Feldkamp-Davis-Kress (FDK) reconstruction algorithms often lack flexibility to accurately account for specific tomographic settings and available a priori knowledge.  
As a remedy, {variational methods} have been proposed.
The idea is to minimize a \emph{generalized Tikhonov functional}:
\begin{align}
 f^{\textup{rec}} \in \argmin_ {f \in L^2(\Omega)} \cS_{\tot} \left( g_{\tot}^\obs; \, F_{\tot} \left( \Proj_{\tot} f \right) \right)  + \cR_{\tot}( f). \label{eq:genTikhonov}
\end{align}
The data-fidelity functional $\cS_{\tot} \left( g_{\tot};\, \cdot \right) $ may account for the expected data-error-statistics, while 
the penalty functional $\cR_{\tot}$ allows to impose a priori knowledge on the object.

The main drawback of variational methods lies in their computational costs. If the minimizer in \eqref{eq:genTikhonov} is computed by generic (convex) optimization methods, such as primal-dual methods \cite{ChambollePock2011,SidkyEtAl2012TomoWithChambollePock}, fast iterative soft shrinkage \cite{Beck2009Teboulle_FISTA} or the alternating direction method of multipliers \cite{Gabay1983ADMM,Setzer2011ADMMForImageProcessing}, many iterations are typically needed for convergence. Moreover, each iteration usually involves an evaluation of the full tomographic projector $\Proj_{\tot}$ and its adjoint $\adj \Proj_{\tot}$ and thus amounts to similar computational costs as a \emph{complete reconstruction} by FBP or FDK.

\paragraph{Kaczmarz-type methods:} One approach to decrease the number of expensive evaluations of $\Proj_{\tot}$ and $\adj \Proj_{\tot} $ compared to (bulk) variational methods is to exploit the block-structure of \ipref{ip1} by performing cyclic iterations on the sub-problems  $g_{j}^\obs = F_j(\Proj_j(f)) + \bepsilon_j$. We consider such \emph{Kaczmarz}- or \emph{block-iterative-}methods in a general variational form
\begin{align}
 f_{k+1} \in \argmin_{f \in L^2(\Omega)} \cS_k \left( g_{j_k}^\obs; \, F_{j_k} \left(\Proj_{j_k} ( f) \right) \right) +  \cR_k(f), \quad j_k \in \{1, \ldots , \Nproj \}.  \label{eq:genKaczmarz}
\end{align}
 The processing order $\{j_k\}$ is typically chosen as an integer number of \emph{cycles} over the whole data, given by permutations of the indices $\{1,\ldots, \Nproj\}$.
 Provided that the functionals $\cS_k, \cR_k$ in \eqref{eq:genKaczmarz} relate to $\cS_{\tot}, \cR_{\tot}$ in \eqref{eq:genTikhonov} in a suitable manner, Kaczmarz-schemes may yield similar reconstructions as the associated bulk variational methods, as will be seen in \secref{SS:KaczmarzVSTikhonov}.


The formulation \eqref{eq:genKaczmarz} generalizes classical Kaczmarz-iterations, see e.g.\ \cite[$\S$V.3]{Natterer}:
\begin{align}
   f_{k+1} &= \argmin_{f \in B_k} \norm{f - f_k}_{L^2}^2 \MTEXT{with}  B_k := \argmin_{f \in L^2(\Omega)} \norm{ \Proj_{j_k} (f ) - g_{j_k} }_{L^2}^2 \nnl
   &= \lim_{\alpha \to 0} \bbparens{ \argmin_{f \in L^2(\Omega)}  \norm{ \Proj_{j_k} (f ) - g_{j_k} }_{L^2}^2 +  \alpha \norm{f - f_k}_{L^2}^2 } \label{eq:classicalKaczmarz}
  \end{align}
  The well-known \emph{algebraic reconstruction technique} (ART) \cite{Kaczmarz1937ART,Gordon1970ART} is an analogue of this approach applied to a fully  forward model, where the projectors $\Proj_{j_k}$ are replaced by single rows of a projector matrix. On the other hand, \emph{simultaneous ART} (SART), proposed in \cite{AndersenKak1984SART} as a heuristic approach to stabilize ART by simultaneously updating all matrix-rows corresponding to a single tomographic projection, can be interpreted as the discretization of an analytical formula for the Kaczmarz-iterate in \eqref{eq:classicalKaczmarz}. This will be seen in \secref{SS:L2SART-example}.


\subsection{Convergence of Kaczmarz-iterations and relation to Tikhonov regularization} \label{SS:KaczmarzVSTikhonov}

Classical Kaczmarz-iterations of the form \eqref{eq:classicalKaczmarz} are known to exhibit fast \emph{semi-convergence}, typically yielding a regularized reconstruction within $\cO(1)$ cycles while increasingly amplifying data-noise if more iterations are preformed, see \cite{ElfvingHansen2014SemiconvKaczmarz} and references therein. Recently, Kindermann and Leit\~ao \cite{KindermannLeitao2014_genARTConvergence} obtained a much more concise characterization of convergence for \emph{regularized} (Tikhonov-)Kaczmarz-iterations, based on previous work by Elfving and Nikazad \cite{ElfvingNikazad2009KaczmarzAsPrecondLandweber}. One of their principal results (cf.\ \cite[Theorem 3]{KindermannLeitao2014_genARTConvergence}) can be cast to the following form:

\vspace{.5em}
   \begin{theorem}[Kindermann \& Leit\~ao \textup{\cite{KindermannLeitao2014_genARTConvergence}}] \label{thm:KindermannLeitao} 
 Let $X,Y_1,\ldots, Y_N$ be Hilbert spaces, $f_0 \in X$ and $g^{\textup{obs}} =  ( g^{\textup{obs}}_1, \ldots, g_N^{\textup{obs}}) \in Y_\tot := \prod_j Y_j$. Moreover, let $A_j: X \to Y_j$ linear-bounded and $A_\tot :=  \transp{ (A_1, \ldots, A_N) }: X\to Y_\tot$. Then, after a \emph{symmetric Kaczmarz-cycle},
\begin{align}
 f_{k+1} = \begin{cases}
           \argmin_{f \in X} \norm{A_{k+1} f - g^{\textup{obs}}_{k+1}}_{Y_j}^2 +  \alpha \norm{f - f_{k}}_X^2 &\text{for } k < N   \\
	   \argmin_{f \in X} \norm{A_{2N-k} f - g^{\textup{obs}}_{2N-k}}_{Y_j}^2 +  \alpha \norm{f - f_{k}}_X^2 &\text{for } k\geq N   \\
         \end{cases}
\end{align}
for $k = 0,1,\ldots, 2N-1$, 
the final iterate minimizes a \emph{bulk Tikhonov-functional}:
\begin{align}
 &f_{2N} = \argmin_{f \in X} \bnorm{W \cdot (A_\tot f - g^{\textup{obs}}) }_{Y}^2 + \frac \alpha 2  \bnorm{ f - f_0 }_X^2  \label{eq:KindermannLeitaoBulkTikh}  \\
 W &= \begin{pmatrix} \textup{id} + \frac 1 {2\alpha} A_1 A_1^\ast & & 0 \\  & \ddots  &\\ 0 & &  \textup{id} + \frac 1 {2\alpha} A_N A_N^\ast \end{pmatrix}^{1/2} 
 \begin{pmatrix} \textup{id}  & -\frac 1 {\alpha} A_1 A_2^\ast  & \ldots  & -\frac 1 {\alpha} A_1 A_N^\ast \\ &  \ddots & \ddots & \vdots \\ 0 & & \ddots & -\frac 1 {\alpha} A_{N-1} A_N^\ast \\ & & & \textup{id}  \end{pmatrix}^{-1} \nonumber
\end{align}
 \end{theorem}
 \vspace{.5em}
 
 \Cref{thm:KindermannLeitao}  is applicable to the considered Kaczmarz-iterations in \eqref{eq:genKaczmarz} if the data-fidelities $\cS_k$ are quadratic, the maps $F_j$ are linear and $\cR_k(f) = \alpha \norm{f - f_k}_X^2 $ for some Hilbert-space-norm $\norm{\cdot}_X$. For this case, it states that already \emph{one} symmetric Kaczmarz-cycle essentially yields a minimizer of the associated bulk Tikhonov functional, up to a slight modification of the data-fidelity term by the  operator $W$ in \eqref{eq:KindermannLeitaoBulkTikh}. Notably, $W$ deviates from the identity only by contributions of order $ \alpha^{-1}$. This implies that the result is not applicable to the limit $\alpha \to 0$ of the classical Kaczmarz-method \eqref{eq:classicalKaczmarz}. On the contrary,
 for sufficiently large $\alpha$, 
 \cref{thm:KindermannLeitao} shows that regularized Kaczmarz-iterations may be used to approximate Tikhonov-minimizers, \IE to \emph{emulate} bulk variational methods -- at least in the considered quadratic setting.
 Other convergence results pointing in a similar direction  are given in \cite{Bertsekas2011_IncrProxMethods,Bertsekas2011_IncrProxMethods2,AndersenHansen2014_genARTProx}.

 \subsection{Contribution} In order for Kaczmarz-methods to provide a truly efficient alternative to bulk variational approaches, essentially two conditions have to be satisfied:
  \optspace
 \begin{enumerate}
  \item The iterates $f_k$ arrive at a reasonable reconstruction after few cycles (ideally one).  \optspace
  \item The individual Kaczmarz-iterations may be evaluated at low computational costs.
 \end{enumerate}
  \optspace
 The convergence theory reviewed in \secref{SS:KaczmarzVSTikhonov} indicates that condition~1 is often fulfilled, in agreement with empirical observations for Kaczmarz-reconstructions in practice.
 The present article, on the other hand, does \emph{not} focus on convergence but contributes to meeting condition~\emph{2} by proposing computationally efficient solution-formulas, termed generalized SART-schemes, for Kaczmarz-iterations of the general variational form \eqref{eq:genKaczmarz}.

\optspace
\section{The Generalized SART-Principle} \label{S:SARTPrinciple}

\subsection{Preparations: notation and analysis of the projectors} \label{SS:ProjProp}

In order to analyze the Kaczmarz-iterations \eqref{eq:genKaczmarz}, we need some properties of the projection operators $\Proj \in \{\ParProj_{\btheta}, \Cone_{\src}\}$ for a single incident direction $\btheta$ or source position~$\src$, see \secref{SS:TomoModel}. First of all, we note that different values of $\btheta$ and $\src$ merely correspond to a rotation and shift of the coordinate system, respectively. We may therefore restrict our analysis  without loss of generality to the cases $\btheta = (0,0,1)^{\textup T}$ and $\src = \boldsymbol 0$, i.e.\ to the operators
 \begin{align}
 \ParProj  (f)( \bx_\perp ) &:= \int_{\OmegaX} \! f(\bx_\perp,z) \, \D z,\quad \Cone (f) ( \bvarphi )  := \int_\OmegaPhi \! f( t \bvarphi) \; \D t \MTEXT{for} \bx_\perp \in \mR^2, \bvarphi \in \mS^2 \label{eq:DefProj0} 
\end{align}
where we have defined the \emph{ray-segments} intersecting $\Omega$ in the parallel- or cone-beam setting:
 \begin{subequations} \label{eq:DefRaysDetDoms}
   \begin{align}
     \OmegaX &:= \{ z \in \mR: (\bx_\perp, z) \in \Omega\}, \qquad  \DR := \{ \bx_\perp \in \mR^2:\OmegaX \neq \emptyset \} \label{eq:DefRays} \\
    \OmegaPhi &:= \{ t \in [0;\infty): t\bvarphi \in \Omega\}, \qquad \,  \DC  := \{ \bvarphi \in \mS^2:\OmegaPhi \neq \emptyset \} \label{eq:DefDetDoms} 
    \end{align}
   \end{subequations}
By definition, $\ParProj  (f)$ and $\Cone(f)$ vanish outside the \emph{projection-domains} $\DR  \subset \mR^2, \DC \subset \mS^2$ for any $f \in L^2(\Omega)$.
Moreover, as $\Omega$ is open, so are the sets $\OmegaX, \OmegaPhi \subset \mR $ and $\DR \subset \mR^2, \DC \subset \mS^2$.
To enable a unified treatment of parallel- and cone-beam settings, we set $\mD := \mR^2$ if $\Proj = \ParProj$ and $\mD := \mS^2$ if $\Proj = \Cone$.

We define the \emph{ray-density functions} $w_{\Proj}$ and \emph{(weighted) unit-projections} $\uproj, \tuproj $: 
\begin{align}
w_{\Proj}: \Omega \to \mR_{>0}; \; \bx \mapsto \begin{cases} 
                    1&\textup{for } \Proj = \ParProj \\
                    |\bx|^{-2} &\textup{for } \Proj = \Cone
                   \end{cases},
                   \qquad \uproj := \Proj( \boldsymbol 1_\Omega ), \quad \tuproj := \Proj( w_{\Proj} )
 \label{eq:DefRayDensityUnitProj}
\end{align}
where $\boldsymbol 1_\Omega \in L^2(\Omega)$ is the constant 1-function in $\Omega$. Note that $\tuproj = \uproj$ for $\Proj = \ParProj$, whereas $\tuproj \neq \uproj$ for $\Proj = \Cone$, and that $\DP = \{ x \in \mD : \uproj(x) > 0\} = \{ x \in \mD : \tuproj(x) > 0\}$.

\paragraph{Back-projections:} We introduce the \emph{back-projections} corresponding to $\ParProj$ and $\Cone$:
 \begin{align}
  \backParProj(p)(\bx_\perp, z) &:= \begin{cases}
                                  p(\bx_\perp) &\textup{if } (\bx_\perp, z) \in \Omega \\
                                  0 &\textup{else}
                                 \end{cases}, \qquad \backCone(p)(\bx) :=  \begin{cases}
                                  p(\bx/|\bx|) &\textup{if } \bx \in \Omega \\
                                  0 &\textup{else}
                                 \end{cases} \label{eq:DefBackProj}
\end{align}
 $\backParProj(p)$ and $\backCone(p)$ are constant along all rays $\OmegaX$ or $\OmegaPhi$, which are exactly the integration-domains in the definition \eqref{eq:DefProj0} of $\ParProj$ and $\Cone$, respectively. 
 In particular, this implies that
\begin{align}
 \Proj\left( w \cdot \backProj (p) \right)  = \Proj(w) \cdot p \MTEXT{for} \Proj \in \{\ParProj, \Cone \}  \label{eq:ProjBackProjRel}
\end{align}
and all $w, p$ for which this expression makes sense. 

\paragraph{Spaces:} We study the maps  $\ParProj, \Cone$ as operators between spaces of square-integrable functions $L^2(\Omega) := \{ f \in L^2(\mR^3): \supp(f) \subset \Omega \}$ and $\DetDom := \{ p \in L^2(\mD): \supp(p) \subset \DP \}$ ($\DP \subset \mD$ as defined in \eqref{eq:DefRaysDetDoms}), with the usual inner product $\ip{f}{g}_{L^2} = \int_\Omega f(x) \cc{g(x)}\, \D x$ and corresponding norm \smash{$\norm{f}_{L^2}:= \ip{f}{f}_{L^2}^{1/2}$}. For simplicity, we restrict to {real-valued} functions throughout this work, although all results naturally carry over to the complex-valued case. To avoid confusion, it is crucial to note that our notation will follow \cref{conv:ProjOps}:

\vspace{.5em}
\begin{convention}  \label{conv:ProjOps}
 All local operations ($+$,$-$,$\cdot$, $/$, $(\cdot)^\gamma$, $\nabla$, \ldots) on functions $f \in L^2(\Omega) \subset L^2(\mR^3)$ or $p \in \DetDom\subset L^2(\mD)$ are implicitly understood to be performed only within the open domains $\Omega$ or $\DP$. For example, a quotient of $p_1, p_2 \in \DetDom$ is to be read as
 \begin{align}
  (p_1/p_2) (x) = \begin{cases}
                                  p_1(x)/p_2(x)  &\textup{for } x \in \DP \\
                                  0 &\textup{else}
                \end{cases} \MTEXT{for all} x \in \mD.
 \end{align}
\end{convention}
\optspace

\paragraph{Adjoints:} For $\Proj \in \{ \ParProj, \Cone \}$, it can be shown that $\Proj: L^2(\Omega) \to \DetDom$ is a bounded linear operator. The adjoints
are given by (weighted) back-projections  (see e.g.\ \cite{Natterer,Louis2016_ExactConeReconFormula}):
\begin{align}
 \adj\Proj: \DetDom \to L^2(\Omega); \; \Proj^\ast(p) = w_{\Proj} \cdot \backProj(p) \MTEXT{for} \Proj \in \{ \ParProj, \Cone \}. \label{eq:ProjAdj}
\end{align}
\noindent Notably, the relation \eqref{eq:ProjBackProjRel} applied to \eqref{eq:ProjAdj} yields
 \begin{align}
   \Proj \adj\Proj(p) = \Proj\left( \wproj \cdot \backProj (p) \right)  = \tuproj \cdot p \MTEXT{for all} p \in \DetDom. \label{eq:PPastMulti}
 \end{align}
\begin{remark} \label{rem:PPastEfficient}
 The formula \eqref{eq:PPastMulti} is of high computational value: it states that, while evaluations of the operators $\Proj$ and $\Proj^\ast$ \emph{alone} may be complicated and costly to compute, the composition $ \Proj \Proj^\ast$ can be implemented as a simple pointwise multiplication of projections. This observation is a key ingredient to an efficient computation of Kaczmarz-iterates.
\end{remark}
 \optspace

 \paragraph{Geometric characterization:} In the analysis, we will need to consider \emph{weighted} projectors:
  \begin{align}
 \isoProj: f \mapsto \tuproj^{-1/2} \cdot \Proj (f)  \MTEXT{for} \Proj \in \{\ParProj, \Cone \}.  
 \label{eq:DefIsoProj}
 \end{align}
 Note that the expression is well-defined by \cref{conv:ProjOps} since $\tuproj(x) > 0$ for all $x \in \DP$.
 At the first glance, the definition of $\isoProj$ may still seem artificial. Yet, it enables a compact geometric characterization of the projectors, which is probably common knowledge at least for $\Proj = \ParProj$, see e.g.\ \cite[$\S$ V.4.3]{Natterer}. 
 The proof of the following result is given in appendix~\ref{Appendix:ProjPropProof}:

\vspace{.5em}
\begin{theorem}[Geometry of the projectors] \label{thm:AdjProjClosedRange} 
For $\Proj\in  \{\ParProj, \Cone \}$, $\isoProj: L^2(\Omega) \to \DetDom$ is bounded with norm $1$ and the adjoint \smash{$\isoProj^\ast: \DetDom \to L^2(\Omega); \; p \mapsto \wproj \cdot \backProj( \tuproj^{- 1/2 } \cdot p)$} is isometric. In particular, $\isoProj$ and $\isoProj^\ast$ have closed range,  $\isoProj^\ast \isoProj$ is the orthogonal projection onto $ \Range(\isoProj^\ast) $ and $ \isoProj \isoProj^\ast = \identity_{\DetDom}$ is the identity on $\DetDom$.
 \end{theorem}
\vspace{.5em}


\subsection{$L^2$-SART: a promising example} \label{SS:L2SART-example}

As a motivation for the general result presented in the subsequent section, we consider an example of Kaczmarz-iterations, that turn out to be computable via a simple analytical formula. Let the iterates be defined by
\begin{align}
 f_{k+1} \in \argmin_{f \in L^2(\Omega)} \norm{ \Proj_{j_k} ( f) - g^\obs_{j_k} }_{L^2}^2  +  \alpha  \norm{ f - f_k }_{L^2}^2  \label{eq:L2Kaczmarz}
\end{align}
with regularization parameter $\alpha > 0$ and $g^\obs_{j_k} \in \DetDom$. The unique minimizer of this variational problem is given by the solution to the corresponding normal equation:
\begin{align}
 f_{k+1} - f_k &= \left( \adj{\Proj_{j_k}} \Proj_{j_k} + \alpha \right)^{-1} \adj{\Proj_{j_k}} ( g^\obs_{j_k} - \Proj_{j_k} (f_k) ) \nnl
 &= \adj{\Proj_{j_k}} \left(  \Proj_{j_k} \adj{\Proj_{j_k}} + \alpha \right)^{-1}  ( g^\obs_{j_k} - \Proj_{j_k} (f_k) ). \label{eq:L2KaczmarzNormalEq}
\end{align}
The second equality in \eqref{eq:L2KaczmarzNormalEq} uses an identity from functional calculus, see e.g.\ \cite{Hanke1996Regularization}. According to \eqref{eq:PPastMulti}, $\big(  \Proj_{j_k} \adj{\Proj_{j_k}} + \alpha \big)^{-1}$ is a division by $\tuprojj{j_k} + \alpha$ (with $\tuprojj{j_k}:= \tuprojj{P_{j_k}}$) so that \eqref{eq:L2KaczmarzNormalEq} yields
\begin{align}
 f_{k+1}  = f_k + \adj{\Proj_{j_k}} \left(  \frac{ g^\obs_{j_k} - \Proj_{j_k} (f_k) }{ \tuprojj{j_k} + \alpha  } \right), \label{eq:L2SARTContFormula}
\end{align}
 where all arithmetic operations are understood to be pointwise. 

 \paragraph{Relation to SART:} By exchanging the continuous object density $f_k$, projection data $g^\obs_{j_k}$, projector $\Proj_{j_k}$ and unit-projection $\tuprojj{j_k}$ in \eqref{eq:L2SARTContFormula} with suitable discretizations $\BP_{j_k} \in \mR^{m\times n}, \Bf_k \in \mR^n, \btu_{j_k}, \Bg_{j_k}^\obs  \in \mR^m$, a numerically implementable update-formula is obtained:
 \begin{align}
   \Bf_{k+1} = \Bf_k + \adj{\BP_{j_k}} \left( \left( \Bg_{j_k}^\obs - \BP_{j_k} \Bf_k \right) \oslash \left( \btu_{j_k} + \alpha \right) \right) \label{eq:L2SARTDiscrFormula}
 \end{align}
 where $\oslash$ denotes element-wise division of vectors. Notably, the iterate $\Bf_{k+1}$ from \eqref{eq:L2SARTDiscrFormula} is in general \emph{not} a solution to the discrete analogue of \eqref{eq:L2Kaczmarz}, i.e.\ typically
 \begin{align}
  \Bf_{k+1} \notin \argmin_{\Bf \in \mR^n} \norm{ \BP_{j_k} \Bf - \Bg_{j_k}^\obs }_{2}^2  +  \alpha  \norm{ \Bf - \Bf_k }_{2}^2,
 \end{align}
  because $\BP_{j_k} \adj{\BP_{j_k}}$ -- unlike $\Proj_{j_k} \adj{\Proj_{j_k}}$ -- is not diagonal for standard discretizations. Interestingly however, the update \eqref{eq:L2SARTDiscrFormula} is closely related to SART, that was originally 
 derived without reference to the continuous model. The classical SART-update \cite{AndersenKak1984SART} can be written as
 \begin{align}
     \Bf_{k+1} = \Bf_k + \adj{\BP_{j_k}} \left( \left( \Bg_{j_k}^\obs - \BP_{j_k} \Bf_k \right) \oslash \BP_{j_k}(\boldsymbol 1)   \right)  \oslash \adj{\BP_{j_k}}(\boldsymbol 1), \label{eq:ClassicalSART}
 \end{align}
 where $\boldsymbol 1 = \transp{(1,1,\ldots, 1)}$ denotes one-vectors of suitable length. In the parallel-beam case $\Proj_{j_k} = \ParProj$, one has $\tuprojj{j_k} = \Proj_{j_k}(\boldsymbol 1 _\Omega)$  and $\adj{\Proj_{j_k}}(\boldsymbol 1_{\DP}) = \boldsymbol 1_\Omega$ by \eqref{eq:DefRayDensityUnitProj} and \eqref{eq:DefBackProj}. 
 The term $\BP_{j_k}(\boldsymbol 1)$ in \eqref{eq:ClassicalSART} can thus be identified with $\btu_{j_k}$ in \eqref{eq:L2SARTDiscrFormula} and the  division by $\adj{\BP_{j_k}}(\boldsymbol 1)$ is just redundant from the perspective of the continuous model. Accordingly, the SART-update \eqref{eq:ClassicalSART} essentially corresponds to the limit $\alpha \to 0$ of  \eqref{eq:L2SARTDiscrFormula}. Since \eqref{eq:L2SARTDiscrFormula} has been derived as a discretization of \eqref{eq:L2Kaczmarz}, SART \eqref{eq:ClassicalSART} may thus be interpreted as a formula to compute the classical Kaczmarz-iterations in \eqref{eq:classicalKaczmarz}. 
 Conversely, \eqref{eq:L2SARTDiscrFormula} can be seen as an \emph{$L^2$-regularized SART-variant}.
  
 \paragraph{The SART-Scheme:}  Analogously to classical SART, the update-formula \eqref{eq:L2SARTDiscrFormula} computes the Kaczmarz-iterate in \eqref{eq:L2Kaczmarz} via a highly efficient \emph{non-iterative} scheme, that requires only a single evaluation of $\BP_{j_k}$ and $\adj{\BP_{j_k}}$ each: 
 \vspace{.5em}
  \begin{scheme}[$L^2$-regularized SART] \label{scheme:SART}  \ \vspace{-1em} \\
 \begin{itemize}
  \setlength{\itemindent}{2em}
  \item[\enum 1] Forward-project the current iterate: $\Bp_k = \BP_{j_k} (\Bf_k )$ \vspace{.25em}
  \item[\enum 2] Compute increment in projection space: $\Delta \Bp_k = (\Bg_{j_k} - \Bp_k ) \oslash( \btu_{j_k} + \alpha)$ \vspace{.25em}
  \item[\enum 3] Back-project and increment: $\Bf_{k+1} = \Bf_k + \adj{\BP_{j_k}} (\Delta \Bp _k )$
 \end{itemize}
  \end{scheme}
 \vspace{.5em}
 Notably, the actual data-fitting step \enum 2 works exclusively on \emph{2D-projections} $\Bg_{j_k}, \Bp_k, \btu_{j_k} \in \mR^m$, which are typically discretized by $m = \cO(N^2)$ pixels and thus \emph{low-dimensional} compared to the discrete \emph{3D-objects} $\Bf_k \in \mR^n$ with $n = \cO(N^3)$ ($N$ sampling-points per spatial dimension). This renders step \enum 2 cheap to compute, so that the total computational costs of \cref{scheme:SART} are essentially that of the evaluations of $\BP_{j_k}$ and $\adj{\BP_{j_k}}$  in steps \enum 1 and \enum 3. Combined with the fast (semi-)convergence of Kaczmarz-methods in $\cO(1)$ cycles (see \secref{SS:KaczmarzVSTikhonov}), this implies that $L^2$-SART-schemes allow to compute reconstructions at a favorable computational complexity of $\cO(1)$ evaluations of the full projector $\BP_{\tot} = (\BP_1,\ldots, \BP_{\Nproj})$ and $\adj \BP_{\tot}$, compare \secref{SS:RecMethComplexity}.

\subsection{Generalized SART framework} \label{SS:genSART}

In \secref{SS:L2SART-example}, it has been shown that $L^2$-Kaczmarz-iterations \eqref{eq:L2Kaczmarz} can be evaluated in a simple and efficient manner although they involve a seemingly complicated optimization problem. Motivated by this result, we 
explore in how far more general Kaczmarz-iterations permit an efficient computation analogous to the SART-like \cref{scheme:SART}. For notational convenience, we drop the subscripts in \eqref{eq:genKaczmarz} at this point and abbreviate the data-fidelity as $ \cS\left(p \right) := \cS_{k} \big( g^\obs_{j_k}; \, F_{j_k} \left(p \right) \big)$, absorbing the maps $F_{j_k}$ into the functional. Thereby, the considered Kaczmarz-iterations are cast to the generic form
\begin{align}
 f_{\new} \in \argmin_{f \in L^2(\Omega)} \cS \left( \Proj ( f) \right) +  \cR(f) =  \argmin_{f \in L^2(\Omega)} \cSAbstr \bparens{ \ProjAbstr ( f) } +  \cR(f).   \label{eq:genKaczmarzFormal}
\end{align}

On the right-hand side of \eqref{eq:genKaczmarzFormal}, it has been used that $(\Proj, \cS)$ may be replaced by any pair $(\ProjAbstr, \cSAbstr)$ of modified projectors and data-fidelities, provided that $\cSAbstr \parens{ \ProjAbstr ( f) } = \cS \parens{ \Proj(f) }$ for all $f \in L^2(\Omega)$. We have to exploit this freedom since our mathematical framework will require $\ProjAbstr$ to have \emph{closed range}, which is not satisfied for $\ProjAbstr = \Proj \in \{\ParProj, \Cone\}: L^2(\Omega) \to L^2(\DP)$, but only for suitably weighted versions like $\ProjAbstr = \isoProj$, compare \cref{thm:AdjProjClosedRange}. We consider the optimization problem in \eqref{eq:genKaczmarzFormal} for a general bounded linear operator $\ProjAbstr: X\to Y$ on Hilbert spaces $X, Y$.
In this setting, we find that the key ingredient to computing \eqref{eq:genKaczmarzFormal} via a SART-like scheme is the following assumption on the geometrical compatibility of $\ProjAbstr$ and $\cR$:

\vspace{.5em}
\begin{assump} \label{A1}
 Let $\ProjAbstr: X \to Y$ be a  bounded linear operator on Hilbert spaces $X, Y$ with null-space $\Kern(\ProjAbstr) = \{f \in X : \ProjAbstr(f) = 0\}$ such that $\Range(\ProjAbstr) = \ProjAbstr (X) \subset Y $ is closed and let $ \cR: X \to \mRinf$ be a functional. Assume that there exists an $f_{\REF} \in X$ such that 
   \begin{align}
   \cR(f_{\REF}  + \ProjAbstr^\ast(p) + f_0) \geq \cR(f_{\REF}  + \ProjAbstr^\ast(p) ) \MTEXT{for all} p \in Y, \, f_0 \in \Kern(\ProjAbstr). \tag{A} \label{eq:A1}  
  \end{align}
\end{assump}
\vspace{.5em}

\Cref{A1} ensures that the penalty functional  $\cR$ uniformly penalizes deviations from a certain reference element $f_{\REF}$ by elements from the null-space of $\ProjAbstr$.
If this condition is satisfied, \eqref{eq:genKaczmarzFormal} can be evaluated according to our principal theorem:

\vspace{.5em}
  \begin{theorem}[Generalized SART-principle] \label{thm:genSART}  
 Let \cref{A1} be satisfied and let $ \cSAbstr: Y \to \mRinf$ be any functional. Assume that there exists a minimizer  
 \begin{align}
 f_{\new} \in \argmin_{f \in X} \cSAbstr\bparens{ \ProjAbstr(f) } +  \cR( f ). \label{eq:KaczmarzGeneric}
\end{align}
  Then there is a (possibly distinct) minimizer $\tilde f_{\new} \in X$ of \eqref{eq:KaczmarzGeneric} given by  
  \begin{subequations} \label{eq:genSART}
  \begin{align}
   p_{\REF} &= \ProjAbstr (f_{\REF}  ) \label{eq:genSARTa}  \\
  \Delta p  &\in \argmin_{p \in Y}   \cSAbstr \bparens{  p_{\REF} + \ProjAbstr \ProjAbstr^\ast ( p ) } + \cR \bparens{ f_{\REF}  + \ProjAbstr^\ast ( p ) }  \label{eq:genSARTb}  \\
   \tilde f_{\new} &= f_{\REF}  + \ProjAbstr^\ast (\Delta p) \label{eq:genSARTc}  
  \end{align}
  \end{subequations}
  Conversely, any $\tilde f_{\new}$ given by \eqref{eq:genSART} minimizes \eqref{eq:KaczmarzGeneric}. 
  Furthermore, if strict inequality holds in \eqref{eq:A1} whenever $f_0 \neq 0$, then all minimizers of \eqref{eq:KaczmarzGeneric} are of the form \eqref{eq:genSART}.
 \end{theorem}
 \vspace{.5em}
 \begin{proof}
 Let $\Delta f := f_{\new}- f_{\REF} $. Since $\ProjAbstr$ has closed range, the same holds true for the adjoint $\ProjAbstr^\ast$ by the closed-range-theorem. Hence, there is an orthogonal decomposition $X = \Range(\ProjAbstr^\ast) \oplus_\perp \Kern(\ProjAbstr)$ so that, in particular, there exist $f_0 \in \Kern(\ProjAbstr)$ and $\Delta p \in Y$ such that $\Delta f  = \ProjAbstr^\ast ( \Delta p ) + f_0$. Now define $\tilde f_{\new} := f_{\REF}  + \ProjAbstr^\ast (\Delta p)$. Then \eqref{eq:A1} implies that
 \begin{align}
  \cR(\tilde f_{\new}) =  \cR(f_{\REF}  + \ProjAbstr^\ast(p) ) \leq   \cR(f_{\REF}  + \ProjAbstr^\ast(p) + f_0) =  \cR( f_{\new}).
 \end{align}
 Moreover, it holds that $ \cSAbstr \bparens{ \ProjAbstr (\tilde f_{\new}) } = \cSAbstr \bparens{ \ProjAbstr (f_{\new}) }$ since $\tilde f_{\new} - f_{\new} \in \Kern(\ProjAbstr)$. Thus, 
  \begin{align}
   \cSAbstr \bparens{\ProjAbstr (\tilde f_{\new}) } +  \cR \bparens{ \tilde f_{\new}} \leq \cSAbstr \bparens{ \ProjAbstr (f_{\new}) } +   \cR( f_{\new})
 \end{align}
 where the left-hand side is strictly smaller if $\tilde f_{\new} \neq f_{\new}$ under the additional assumption that strict inequality holds in \eqref{eq:A1} for $f_0 \neq 0$. This proves that $\tilde f_{\new}$ is a minimizer and that $\tilde f_{\new} = f_{\new}$ must hold in the case of strict inequality. Moreover, $\Delta p$ satisfies
 \begin{align}
  \ProjAbstr^\ast ( \Delta p ) &\in \argmin _ {\Delta f \in  X }  \cSAbstr \bparens{\ProjAbstr (f_{\REF}  + \Delta f ) } +  \cR(f_{\REF}  + \Delta f) \nnl
  \Rightarrow \;\; \Delta p &\in \argmin _ {p \in Y }  \cSAbstr \bparens{\ProjAbstr (f_{\REF} ) + \ProjAbstr \ProjAbstr^\ast ( p ) } +  \cR(f_{\REF}  +  \ProjAbstr^\ast ( p ) ),
 \end{align}
 which proves  \eqref{eq:genSART}. 
 For the converse statement, let  $\tilde f_{\new}$ be given by \eqref{eq:genSART}. Then $\tilde f_{\new}$ minimizes the cost-functional $ \mathcal C(f):= \cSAbstr\parens{ \ProjAbstr(f) } +  \cR( f )$ over all $f \in A := f_{\REF} + \Range(\adj \ProjAbstr)$. Since \cref{A1} implies that 
 $\mathcal C  (\tilde f_{\new} + f_0) \geq \mathcal C(\tilde f_{\new})$
 for all $f_0 \in \Kern(\ProjAbstr)$, $\tilde f_{\new}$ also minimizes $\mathcal C$ over $X = A + \Kern(\ProjAbstr)$ and hence solves \eqref{eq:KaczmarzGeneric}.
 \end{proof}
 \vspace{.5em}

 We aim to apply \cref{thm:genSART} to the tomographic Kaczmarz-iterations in \eqref{eq:genKaczmarzFormal}. Although stated in an abstract manner, the result is particularly well-suited for this application:
 \optspace
 \begin{enumerate}
  \item[(a)] The optimization problem \eqref{eq:genSARTb} is on the image-space $Y$ of $\ProjAbstr$ (\emph{projection-space}), which is much lower-dimensional than the object-space $X$ in the considered setting. 
   \optspace
  \item[(b)] The operator $\ProjAbstr \adj \ProjAbstr$ appearing in \eqref{eq:genSARTb} is trivial to evaluate for tomographic projectors $\ProjAbstr  \in \{\isoParProj, \isoCone\}$, contrary to $\ProjAbstr$ or $\adj \ProjAbstr$ \emph{alone} (see \cref{rem:PPastEfficient,thm:AdjProjClosedRange}).
 \end{enumerate}
  \optspace
 To apply \cref{thm:genSART}, \emph{literally nothing} has to be assumed on the functional $\cSAbstr$, so that the result offers complete freedom in the choice of the data-fidelity $\cS$ in \eqref{eq:genKaczmarzFormal}. What remains is to verify \cref{A1} for $\ProjAbstr  \in \{\isoParProj, \isoCone\}$ (or related choices) and suitable penalties $\cR$. This is established in \secref{S:AdmPenalties}. 
 We refer to the resulting formulas of the kind \eqref{eq:genSART} as \emph{generalized SART-} or \emph{GenSART-schemes}, aluding to their structural similarity to \cref{scheme:SART}. 

\optspace
\section{Admissible penalty functionals} \label{S:AdmPenalties}

The aim of this section is to identify penalty functionals $\cR = \cR_k$ that satisfy \cref{A1}, in which case the Kaczmarz-iterations in \eqref{eq:genKaczmarz} can be computed via GenSART-schemes by virtue of \cref{thm:genSART}.

\subsection{Preliminary insights} \label{SS3.0}

In order to gain an intuition for the admissible penalties, let us first discuss the meaning of the condition \eqref{eq:A1}. It asserts that -- relative to a certain reference object $f_{\REF}$ -- any deviation by an element from the null-space $\Kern(\ProjAbstr)$ is penalized or at least does not decrease the value of $\cR$.
For (weighted) parallel-beam projectors $\ProjAbstr = \isoParProj$,  elements in $\Kern(\ProjAbstr) = \Kern(\isoParProj) = \Kern(\ParProj)$ are exactly functions that have zero mean along each tomographic ray, whereas an element $f \in L^2(\Omega)$ is in $\Range(\ProjAbstr^\ast) = \closure{\Range(\ParProj^\ast)}$ if and only if $f$ is constant in $\Omega$ along the ray-direction, see \secref{SS:ProjProp}.
The gist of \cref{thm:genSART} for parallel-beam geometries can thus be stated as follows: if $\cR$ uniformly penalizes oscillations along the rays, the increment  computed in \eqref{eq:genKaczmarz} is constant in ray-direction and thus $f_{\new} - f_{\REF} \in \closure{\Range(\ParProj^\ast)}$. 
In cone-beam settings $\ProjAbstr = \isoCone$, the characterization of $\Kern(\ProjAbstr)$ is identical as in the parallel-beam case. However, elements in $\Range(\ProjAbstr^\ast)=\closure{\Range(\Cone^\ast)}$ are only constant along the rays up to a scaling with the ray-density, see \eqref{eq:ProjAdj}.

In a nutshell, the above indicates that \emph{smoothing} penalties, that tend to damp out variations of  $f_{\new} - f_{\REF}$ wherever possible, are promising candidates to satisfy \cref{A1}. Indeed, it will be verified for (weighted) $L^2$-, quadratic gradient- and general $L^q$-penalties in the following sections. 
 Before proceeding to the analysis of these settings, let us observe that convex combinations of admissible penalty functionals still satisfy \cref{A1}. This enables an application of \cref{thm:genSART} to Kaczmarz-iterations with mixed penalties:
 \vspace{.5em}
 \begin{lemma} \label{lem:ConvCombPenalties}
  Let $\alpha_1,\alpha_2 \geq 0$ and let $\cR_1, \cR_2: X \to \mRinf$ be functionals satisfying \eqref{eq:A1} for the same $f_{\REF} \in X$ and $\ProjAbstr:X \to Y$. Then $\cR := \alpha_1 \cR_1 + \alpha_2 \cR_2$ satisfies \eqref{eq:A1}.
 \end{lemma}
 \vspace{.5em}
 \begin{proof}
  This follows by summing scaled versions of  \eqref{eq:A1} for $\cR_1$ and $\cR_2$.
 \end{proof}
  \vspace{.5em}
\subsection{(Weighted) $L^2$-penalties} \label{SS:L2wPenalties}

As a simple candidate for admissible penalty functionals $\cR$ within the 
framework of \cref{thm:genSART}, we consider quadratic penalties of the form
\begin{align}
 \cR (f):= \norm{f - f_{\REF}}_{X}^2,  \label{eq:RkNormsq}
\end{align}
where $\norm{\cdot}_X$ denotes the norm of the Hilbert space $X$. Owing to the geometric nature of the condition \eqref{eq:A1}, it is straightforward to show that such a choice satisfies \cref{A1}. Note that the result is not limited to tomographic projectors $\Proj \in \{\ParProj, \Cone\}$, but holds within the general abstract setting of \cref{A1}:


\vspace{.5em}
  \begin{lemma}[Quadratic norm penalties] \label{lem:admissibleQuadratic}
 Let $X,Y$ be Hilbert spaces, $\ProjAbstr: X\to Y$ linear-bounded with closed range and let $ \cR (f) := \norm{f - f_{\REF}}_X^2$ for $f_{\REF} \in X$. 
 Then \cref{A1} is satisfied with strict inequality in \eqref{eq:A1} for all $f_0 \neq 0$ and it holds that
 \begin{align}
  \cR(f_{\REF} +  \ProjAbstr^\ast(p)) = \ip{p}{\ProjAbstr \ProjAbstr^\ast(p)}_Y \MTEXT{for all} p \in Y. \quad \textup{($\ip{\cdot}{\cdot}_Y$: inner product in $Y$)} \label{eq:admissibleQuadratic}
 \end{align}
 \end{lemma}
 \vspace{.5em}
 \begin{proof}
 In the considered setting, the condition \eqref{eq:A1} follows simply by the orthogonality $\Kern(\ProjAbstr)\perp \Range(\ProjAbstr^\ast)$: for all $p \in \Range(\ProjAbstr^\ast), f_0 \in \Kern(\ProjAbstr)$, it holds that
\begin{align}
 \cR(f_{\REF} + \ProjAbstr^\ast(p) + f_0)  &= \norm{\ProjAbstr^\ast(p) + f_0}_X^2 \stackrel{\ProjAbstr^\ast(p) \perp f_0}= \norm{\ProjAbstr^\ast(p)}_X^2 + \norm{f_0}_X^2  \nnl 
 &\geq  \norm{\ProjAbstr^\ast(p)}_X^2 = \cR(f_{\REF} +  \ProjAbstr^\ast(p)) \nnl
 &= \ip{\ProjAbstr^\ast(p)}{\ProjAbstr^\ast(p)}_X = \ip{p}{\ProjAbstr \ProjAbstr^\ast(p)}_Y
\end{align}
Moreover, since $\norm{f_0}_X^2 > 0$ whenever $f_0 \neq 0$, strict inequality holds in this case. The equality in the third line simply follows by the defining property of the adjoint.
 \end{proof}
 \vspace{.5em}
 
 The abstract result from \cref{lem:admissibleQuadratic} can be applied to establish GenSART-schemes for general Kaczmarz-iterations with $L^2$-penalty:
 
 \vspace{.5em}
  \begin{theorem}[Generalized SART with $L^2$-penalty] \label{thm:genL2SART}
  Let $\Proj \in \{ \ParProj , \Cone \}: L^2(\Omega) \to \DetDom$, $f_{\REF} \in L^2(\Omega)$, $\alpha > 0$ and let $ \cS: \DetDom \to \mRinf$ be any functional. Then the minimizers of
 \begin{align}
  f_{\new} \in \argmin _{ f \in L^2(\Omega) }   \cS  \left( \Proj ( f) \right) +  \alpha  \norm{f - f_{\REF}}_{L^2}^2
 \end{align}
 are uniquely determined by the GenSART-scheme
    \begin{subequations} \label{eq:genL2SART}
   \begin{align}
    p_{\REF} &= \Proj (f_{\REF} ) \label{eq:genL2SART-a}  \\
   \Delta p  &\in \argmin_{p \in \DetDom}    \cS \left(  p_{\REF} + \tuproj^{1/2} \cdot p \right) + \alpha \norm{p}^2_{L^2}  \label{eq:genL2SART-b}  \\
    f_{\new} &= f_{\REF} +  \Proj^\ast ( \tuproj^{-1/2} \cdot \Delta p). \label{eq:genL2SART-c}  
   \end{align}
   \end{subequations}
 \end{theorem}
 \vspace{.5em}
 
 The result  \eqref{eq:genL2SART} constitutes a generalization of the $L^2$-SART-scheme derived in \secref{SS:L2SART-example} for \emph{arbitrary} data-fidelities $\cS$. Indeed, the choice  $\cS(p) = \norm{ p - g^\obs }_{L^2 }^2$ in \eqref{eq:genL2SART} reproduces the formula \eqref{eq:L2SARTContFormula}. Importantly, the optimization problem \eqref{eq:genL2SART-b} in projection-space no longer contains evaluations or $\Proj$ or $\adj\Proj$ and is thus easy to solve given that $\cS $ is sufficiently simple. 
  We omit the proof of \cref{thm:genL2SART} as it is just a special case of the following result: 

 \vspace{.5em}
  \begin{theorem}[Generalized SART with $L^2$-penalty and weighted projector] \label{thm:L2SARTProjW}
 Let $\Proj \in \{ \ParProj , \Cone \}: L^2(\Omega) \to \DetDom$, $f_{\REF} \in L^2(\Omega)$, $\alpha > 0$ and let $ \cS: \DetDom \to \mRinf$ be any functional.
 Moreover, let $\lambda: \Omega \to \mR$ with $\lambda_{\min} \leq |\lambda(\bx)| \leq \lambda_{\max}$ for almost all $\bx \in \Omega$ and some  constants $0 < \lambda_{\min}  \leq \lambda_{\max} < \infty$. Then the minimizers of
 \begin{align}
  f_{\new} \in \argmin _{ f \in L^2(\Omega) }   \cS \left( \Proj ( \lambda \cdot f) \right) +  \alpha \Norm{ f - f_{\REF} }_{L^2}^2. \label{eq:L2SARTProjWOptProb}
 \end{align}
  are uniquely determined by the GenSART-scheme
   \begin{subequations} \label{eq:L2SARTProjW}
  \begin{align}
   p_{\REF} &= \Proj (\lambda \cdot f_{\REF} ) \label{eq:L2SARTProjW-a}  \\
  \Delta p  &\in \argmin_{p \in \DetDom} \cS \bparens{  p_{\REF} + \lambda_{\Proj} \cdot \tuproj^{1/2} \cdot   p } +  \alpha  \bnorm{\lambda_{\Proj}^{1/2} \cdot p}_{L^2}^2   \label{eq:L2SARTProjW-b}  \\
  f_{\new} &= f_{\REF} +  \lambda \cdot \Proj^\ast (\tuproj^{-1/2} \cdot  \Delta p) \label{eq:L2SARTProjW-c} 
  \end{align}
  \end{subequations}
  with $\lambda_{\Proj} = \Proj  ( \wproj \cdot |\lambda|^2)  / \tuproj $, $w_{\Proj}$ denoting the ray-density introduced in \secref{SS:ProjProp}.
 \end{theorem}
 \vspace{.5em}
 \begin{proof}
   The bounds for $\lambda$ ensure that $M_\lambda: L^2(\Omega) \to L^2(\Omega); \; f \mapsto \lambda \cdot f$ is an isomorphism. We set $X:= L^2(\Omega)$, $Y:= \DetDom$ and \smash{$\ProjAbstr := \isoProj \circ M_\lambda: X \to Y$} with \smash{$\isoProj(p) = \tuproj^{-1/2} \cdot \Proj(p)$} as defined in \eqref{eq:DefIsoProj}. Then $\ProjAbstr$ has closed range by \cref{thm:AdjProjClosedRange} and $M_\lambda$ being an  isomorphism. By setting \smash{$\cSAbstr(p) := \cS(\tuproj^{1/2} \cdot p)$} and $\cR(f) := \alpha \norm{ f - f_{\REF} }_{L^2}^2$, \eqref{eq:L2SARTProjWOptProb} is cast to the form \eqref{eq:KaczmarzGeneric}.

   According to \cref{lem:admissibleQuadratic,lem:ConvCombPenalties}, \cref{A1} is satisfied in the strict-inequality-version so that the GenSART-\cref{thm:genSART} is applicable. Hence, the minimizers $f_{\new}$ of \eqref{eq:L2SARTProjWOptProb} can be computed via \eqref{eq:genSART}. Substituting the expressions for $\cSAbstr$, $\ProjAbstr$, $\cR$ and exploiting that $\cR(f_{\REF} +  \ProjAbstr^\ast(p)) = \ip{p}{\ProjAbstr \ProjAbstr^\ast(p)}_Y$ according to \eqref{eq:admissibleQuadratic} yields
   \begin{subequations}
   \begin{align}
      \tilde p_{\REF} &= \ProjAbstr (f_{\REF}  ) = \tuproj^{-1/2}\cdot \ProjAbstr (\lambda \cdot f_{\REF} )  \\
     \Delta p  &\in \argmin_{p \in \DetDom}   \cS \bparens{ \tuproj^{1/2} \cdot \tilde p_{\REF} + \tuproj^{1/2} \cdot  \ProjAbstr\adj \ProjAbstr( p ) } + \alpha \ip{p}{\ProjAbstr\adj \ProjAbstr( p )}_{L^2}     \\
      f_{\new} &=  f_{\REF} + {\lambda} \cdot \adj \ProjAbstr (\tuproj^{-1/2}\cdot p).
   \end{align}
   \end{subequations}
   From this expression, the formula \eqref{eq:L2SARTProjW} is obtained by setting $p_{\REF}:= \tuproj^{1/2} \cdot \tilde p_{\REF}$ and using that, according to \eqref{eq:ProjBackProjRel} and \eqref{eq:ProjAdj}, 
   \begin{align}
    \ProjAbstr \adj \ProjAbstr(p ) &= \tuproj^{-1/2}\cdot  \Proj \bparens{ |\lambda|^2 \cdot \adj \Proj ( \tuproj^{-1/2}\cdot p ) } = \tuproj^{-1/2}\cdot  \Proj \bparens{ (\wproj \cdot |\lambda|^2) \cdot \backProj ( \tuproj^{-1/2}\cdot p ) } \nnl
    &= \tuproj^{-1} \cdot \Proj  ( \wproj \cdot |\lambda|^2) \cdot p =  \lambda_{\Proj} \cdot p \MTEXT{for all} p \in \DetDom. \qedhere
   \end{align}
 \end{proof}
  \vspace{.5em}
  
  Note that, once again, the optimization problem \eqref{eq:L2SARTProjW-b} does not involve $\Proj$ or $\adj \Proj$ except in the form of the precomputable functions $\lambda_{\Proj}$. \Cref{thm:genL2SART} is obtained from \cref{thm:L2SARTProjW} by choosing $\lambda := 1$. Moreover, the theorem enables GenSART for \emph{weighted} $L^2$-penalties:

 \vspace{.5em}
  \begin{corollary}[Generalized SART with weighted $L^2$-penalties] \label{thm:L2wSART}
 Let $\Proj \in \{ \ParProj , \Cone \}: L^2(\Omega) \to \DetDom$, $f_{\REF} \in L^2(\Omega)$, $\alpha > 0$ and let $ \cS: \DetDom \to \mRinf$ be any functional.
 Moreover, let $w: \Omega \to \mR_{> 0}$ with $w_{\min} \leq w(\bx) \leq w_{\max}$ for almost all $\bx \in \Omega$ and some  constants $0 < w_{\min}  \leq w_{\max} < \infty$. Then the minimizers of
 \begin{align}
  f_{\new} \in \argmin _{ f \in L^2(\Omega) }   \cS \left( \Proj ( f) \right) +  \alpha \bnorm{w^{-1/2} \cdot (f - f_{\REF})}_{L^2}^2. \label{eq:L2wSARTOptProb}
 \end{align}
  are uniquely determined by the GenSART-scheme
  \begin{subequations} \label{eq:L2wSART}
  \begin{align}
   p_{\REF} &= \Proj (f_{\REF} ) \label{eq:L2wSART-a}  \\
  \Delta p  &\in \argmin_{p \in \DetDom}  \cS \left(  p_{\REF} + v_{\Proj} \cdot \tuproj^{1/2} \cdot   p \right) +  \alpha  \bnorm{v_{\Proj}^{1/2} \cdot p}_{L^2}   \label{eq:L2wSART-b}  \\
  f_{\new} &= f_{\REF} + w \cdot \Proj^\ast (\tuproj^{-1/2} \cdot \Delta p) \label{eq:L2wSART-c} 
  \end{align}
  \end{subequations}
  where $v_{\Proj} = \Proj \left( w \cdot w_{\Proj} \right) / \tuproj $, $w_{\Proj}$ denoting the ray-density introduced in \secref{SS:ProjProp}.
 \end{corollary}
  \vspace{.5em}
  \begin{proof}
   $f_{\new}$ is a minimizer of \eqref{eq:L2wSARTOptProb} if and only if $\tilde f_{\new} := w^{-1/2} \cdot f_{\new}$ solves the optimization problem \eqref{eq:L2SARTProjWOptProb} for $\lambda := w^{1/2}$ and $f_{\REF}$ replaced by $\tilde f_{\REF} := w^{-1/2} \cdot f_{\REF}$. The resulting GenSART-scheme \eqref{eq:L2SARTProjW} for $\tilde f_{\new}$ yields \eqref{eq:L2wSART}.
\end{proof}
 \vspace{.5em}


The more general, weighted GenSART-schemes from \cref{thm:L2SARTProjW,thm:L2wSART} enable interesting applications, as will be seen in \secref{SSS:PolyCTSART}.
Moreover, note that the assumption that $\lambda$ or $w$ are bounded from below can be dropped at the cost of a more technical proof. Indeed, the formulas \eqref{eq:L2SARTProjW}, \eqref{eq:L2wSART} still make sense if $\lambda$ or $w$ vanishes in parts of $\Omega$.

\subsection{Gradient-penalties} \label{SS:GradPenalties}

 In order to enforce a certain smoothness of the reconstructed object, variational methods often use penalties that involve derivatives. Similarly, we seek to extend the above results to Kaczmarz-iterations \eqref{eq:genKaczmarz} with quadratic gradient penalties:
  \begin{align}
 \cR: L^2(\Omega) \to \mRinf; \; \cR(f) := \norm{ \nabla ( f - f_{\REF} ) }_{L^2}^2 = \int_{\Omega} |\nabla( f -  f_{\REF} )|^2 \; \D \bx 
    \label{eq:PenaltyGradL2}
  \end{align}
 with the convention that $\norm{\nabla h} = \infty$ for any $h \in L^2(U) \setminus W^{1,2}(U)$, where $W^{1,2}(U)$ is the usual Sobolev-space $ W^{1,2}(U) := \{ f \in L^2(U): f\textup{ weakly differentiable in $U$ with } \abs{\nabla f} \in L^2(\Omega) \}$ on an open domain $U \subset \mR^m$, with norm $\norm{f}_{W^{1,2}}^2 := \norm{f}_{L^2}^2 + \norm{\nabla f}_{L^2}^2 $. Moreover, we recall that gradients of functions $f \in L^2(\Omega)$, $p \in \DetDom$ are understood to be computed only \emph{within} $\Omega\subset \mR^3$ and $\DP \subset \mD$, \IE in particular \emph{not} across boundaries (\cref{conv:ProjOps}).  
 
 
  By \secref{SS3.0}, one might hope that the gradient penalty in \eqref{eq:PenaltyGradL2} satisfies \cref{A1} owing to its smoothing effect.
 Indeed, this  turns out to be \emph{almost} true, up to some complications arising from thickness-variations of the domain $\Omega$. This is explained in the following.
 
 For simplicity, we study the parallel-beam case, $\ProjAbstr = \isoParProj$, and continuously differentiable functions $f_0 \in \sC^1(\Omega) \cap \Kern(\ProjAbstr)$ and $p \in \DetDom$ s.t.\ $\ProjAbstr^\ast ( p ) \in \sC^1(\Omega)$. Then we have that
\begin{align}
    \cR(f_{\REF} + \ProjAbstr^\ast(p) + f_0 ) &=  \int_{\Omega} \abs{ \nabla( \isoParProj^\ast(p) + f_0 ) }^2 \; \D \bx \nnl
    &= \norm{\nabla \isoParProj^\ast(p)}_{L^2}^2 + \norm{\nabla f_0}_{L^2}^2 + 2  \int_{\Omega}  \nabla \isoParProj^\ast(p) \cdot { \nabla f_0 } \; \D \bx. \label{eq:GradPen-1}
\end{align}
By \eqref{eq:GradPen-1}, the inequality \eqref{eq:A1} holds true if and only if  the mixed integral on the right-hand side  vanishes for all admissible $p$ and $f_0$. 

Now let us assume that $\Omega$ is convex with $\sC^1$-boundary. Then we may write it in the form $\Omega = \{ (\bx_\perp, z) \in \mR^3 : a(\bx_\perp) < z < b(\bx_\perp)\}$ for continuous functions $a,b:\mR^2 \to \mR$, $a \leq b$, that give the entrance- and exit points of the tomographic rays into $\Omega$. Hence, it holds that
\begin{align}
 \int_{\Omega}  \nabla \isoParProj^\ast(p) \cdot {\nabla f_0} \; \D \bx &= \int_{\mR^2} \bbparens{ \int_{a(\bx_\perp)}^{b(\bx_\perp)}  \nabla \isoParProj^\ast(p)(\bx_\perp, z) \cdot {\nabla f_0(\bx_\perp, z)} \; \D z } \;\D \bx_\perp \nnl
 &= \int_{\mR^2} \nabla \bparens{ p / \uparproj^{1/2} }(\bx_\perp) \cdot \bbparens{ \int_{a(\bx_\perp)}^{b(\bx_\perp)}  {\nabla_\perp f_0(\bx_\perp, z)} \; \D z } \;\D \bx_\perp. \label{eq:GradPen-2}
\end{align}
In the second step, we exploited that $\isoParProj^\ast(p)(\bx_\perp, z) = \backParProj(p / \uparproj^{1/2})(\bx_\perp, z) = p(\bx_\perp)$ is constant in $z$ so that the $z$-component of the gradient $\nabla =  (\nabla_\perp, \partial_z)$ vanishes and $\nabla \adj \isoParProj(p)(\bx_\perp,z) = \nabla(p / \uparproj^{1/2})(\bx_\perp)$ can be pulled out of the inner integral. Since $a$ and $b$ are continuously differentiable within the open set $U := \{ \bx_\perp \in \mR^2 : a(\bx_\perp) < b(\bx_\perp)\}$ by the assumptions on $\Omega$, we can apply Leibniz' rule to the inner integrals in \eqref{eq:GradPen-2}:
\begin{align}
  \int_{a(\bx_\perp)}^{b(\bx_\perp)}  &\nabla_\perp f_0(\bx_\perp, z) \; \D z = \nabla_\perp \int_{a(\bx_\perp)}^{b(\bx_\perp)}   f_0(\bx_\perp, z) \; \D z \nnl
  &- \left( f_0(\bx_\perp, b(\bx_\perp)) \nabla b(\bx_\perp) - f_0(\bx_\perp, a(\bx_\perp))\nabla a(\bx_\perp) \right) \MTEXT{for all} \bx_\perp \in U. \label{eq:GradPen-3}
\end{align}
The first integral on the right-hand side of \eqref{eq:GradPen-3} can be identified as $\nabla \ParProj(f_0) (\bx_\perp)$ and vanishes since $f_0 \in \Kern(\isoParProj) = \Kern(\ParProj)$. The second term, arising from variations of the domain-boundary-functions $a$ and $b$, does \emph{not} vanish in general. Yet, it vanishes if \emph{either}
\vspace{.5em}
\begin{itemize}
 \item $a,b$ are constant within $U$. This is only  the case for domains of the form $\Omega = \Omega_\perp \times \Omega_\parallel$ for some $\Omega_\perp \subset \mR^2$ and $\Omega_\parallel \subset \mR$, such as cuboids aligned with the ray-direction. \vspace{.25em}
 \item $f_0$ vanishes on the boundary of $\Omega$.
\end{itemize}
\vspace{.5em}


 \eqref{eq:GradPen-3}, \eqref{eq:GradPen-2} and \eqref{eq:GradPen-1} indicate that either of the above restrictions may ensure that the inequality \eqref{eq:A1} is satisfied. We focus on the second one in order to retain geometrical flexibility of the object-domain $\Omega$.
In order to verify \cref{A1} in a Sobolev-space setting, we consider the  subspace of functions with trace zero, $W^{1,2}_0(\Omega) := \closure{\CcInfO} \subset W^{1,2}(\Omega)$, defined as the closure of the smooth and compactly supported functions $\CcInfO$ in the $W^{1,2}$-topology. Firstly, we need to analyze the mutual behavior of the projectors and the gradient operator, as established by the following lemma, which is proven in appendix~\ref{Appendix:ProjGradCompatibility}: 
 \vspace{.5em}
 \begin{lemma} \label{lem:ProjGradient}
  Let $\Proj \in \{ \ParProj, \Cone \}$, $f \in W^{1,2}_0(\Omega)$ and $p \in \DetDom$. Denote by $\nabla_{\mD}$ the gradient on the detection domain $\mD \in \{ \mR^2 , \mS^2\}$ and by $\nablaProj$ the component of gradient in $\mR^3$ perpendicular to the local ray-direction of the projector $\Proj$. Then it holds that
  \begin{align}
 \Proj(w_{\Proj}^{-1/2} \cdot \nablaProj f)  &= \nabla_{\mD} \Proj(f) \label{eq:lem-ProjGradient-1} \\
 \norm{\tuproj^{-1/2} \cdot \nabla_{\mD} \Proj(f)}_{L^2} &\leq  \norm{w_{\Proj}^{-1/2} \cdot \nablaProj f}_{L^2}, \label{eq:lem-ProjGradient-2}
   \end{align}
   where $w_{\Proj}$ is the {ray-density} defined in \secref{SS:ProjProp}.
 Moreover, $\backProj(p) \in W^{1,2}(\Omega)$ holds true if and only if \smash{$\uproj^{1/2} \cdot \nabla_{\mD}  \, p  \in \DetDom$} and in this case 
     \begin{align}
 \norm{\nabla  \backProj(p) }_{L^2} &= \bnorm{\uproj^{1/2}  \cdot \nabla_{\mD} \, p }_{L^2} \label{eq:lem-ProjGradient-3} \\
 \ip{\nabla \backProj(p) }{ \nabla f }_{L^2} &=  \bip{\uproj^{1/2} \cdot \nabla_{\mD} p }{ \uproj^{-1/2} \cdot \nabla_{\mD} \Proj(f) }_{L^2}.  \label{eq:lem-ProjGradient-4}
  \end{align}
 \end{lemma}
  \vspace{.5em}  
  
It should be emphasized that the occurrence of both $\tuproj$ and $\uproj$ in \cref{lem:ProjGradient} is \emph{not} a typo, but is indeed necessary to correctly cover the cone-beam case $\Proj = \Cone$. In a parallel-beam setting $\Proj = \ParProj$, on the other hand, one has $\tuproj = \uproj$ and $\wproj =1$, so that the expressions can be simplified.
\Cref{lem:ProjGradient} permits to prove the admissibility of gradient-penalties, yet in a slightly restricted setting due to the ``geometrical complications'' explained above:

\vspace{.5em}
 \begin{lemma}[Admissibility of $L^2$-gradient-penalties] \label{lem:gradientL2Admissible}
 Let $\Proj \in \{ \ParProj, \Cone  \}$,  $f_{\REF} \in L^2(\Omega)$ and let $\cR$ be defined by \eqref{eq:PenaltyGradL2}.
Let $X := ( L^2(\Omega), \ip{\cdot}{\cdot}_\Proj)$ be equipped with the inner product $ \ip{f_1}{f_2}_\Proj := \ip{w_\Proj \cdot f_1}{f_2}_{L^2}$.
and define $\ProjAbstr : X \to \DetDom; \, f \mapsto \uproj^{- 1 /2  } \cdot \Proj\left(f\right)$. Then \cref{A1} is satisfied if we restrict to elements $f_0 \in \Kern(\ProjAbstr) \cap W^{1,2}_0(\Omega) $ in \eqref{eq:A1} and strict inequality  holds for all $f_0 \neq 0$. Moreover, 
  \begin{align}
    \cR(f_{\REF} + \ProjAbstr^\ast(p) ) = \bnorm{ \uproj^{1/2} \cdot \nabla_{\mD} \bparens{ \uproj^{-1/2} \cdot p  } }_{L^2}^2  \label{eq:thm3-1}
  \end{align}
  for all $p \in \DetDom$ where possibly both sides may attain the value $\infty$.
 \end{lemma}

 \vspace{.5em}
 \begin{proof}
 By definition, it holds that $\ProjAbstr = M \circ \isoProj \circ \iota_X$, where $\isoProj$ denotes the weighted projector from \cref{thm:AdjProjClosedRange}, $M: \DetDom \to \DetDom; \; p \mapsto (\tuproj /\uproj)^{ 1/2 } \cdot p $ and  $\iota_{X}: X \hookrightarrow L^2(\Omega); \; f \mapsto f$ is the canonical embedding. 
 The norms of $X$ and $L^2(\Omega)$ are equivalent since $c \leq w_{\Proj} \leq C$ is bounded from below and above by constants $0 < c \leq C$. The bounds also imply that $c^{1/2} \leq (\tuproj /\uproj)^{ 1/2 } \leq C^{1/2}$ by \eqref{eq:DefRayDensityUnitProj}. Hence, $\iota_X$ and $M$ are isomorphisms so that $\ProjAbstr = M \circ \isoProj \circ \iota_X$ has closed range by \cref{thm:AdjProjClosedRange}. Moreover, it holds that $\iota_X^\ast(f) = w_{\Proj}^{-1} \cdot f$ for all $f \in L^2(\Omega)$. Hence, we have for all $p \in \DetDom$
 \begin{align}
      \adj \ProjAbstr( p ) = \iota_X^\ast \bparens{ \adj \Proj \bparens{ \uproj^{-1/2} \cdot p }  } = w_{\Proj}^{-1} \cdot \adj \Proj \bparens{  \uproj^{-1/2} \cdot p } \stackrel{\eqref{eq:ProjAdj}}= \backProj \bparens{  \uproj^{-1/2} \cdot p }. \label{eq:gradientL2Admissible-pf1}
 \end{align} 
  Now let $p\in \DetDom$ and $ f_0 \in \Kern(\ProjAbstr) \cap W^{1,2}_0(\Omega) $ be arbitrary. If $\ProjAbstr^\ast( p ) \not\in W^{1,2}(\Omega)$, then also $\ProjAbstr^\ast( p ) + f_0 \notin W^{1,2}(\Omega)$ and so $\cR(f_{\REF} + \ProjAbstr^\ast(p) + f_0 ) = \infty$ by \eqref{eq:PenaltyGradL2}, in which case  \eqref{eq:A1} is trivially satisfied. Hence, we assume that $\ProjAbstr^\ast( p ) \in W^{1,2}(\Omega)$ in the following. Then it holds that
    \begin{align} 
  \cR(f_{\REF} + \ProjAbstr^\ast(p) + f_0 )  &=  \norm{\nabla\parens{\ProjAbstr^\ast( p ) + f_0 } }_{L^2}^2  \nnl
  &=   \norm{\nabla \ProjAbstr^\ast( p ) }_{L^2}^2  +  2\ip{\nabla \ProjAbstr^\ast( p )}{\nabla f_{0}}_{L^2} +  \norm{\nabla f_0}_{L^2}^2  \nnl 
  &\geq  \norm{\nabla \ProjAbstr^\ast( p ) }_{L^2}^2  +  2\ip{\nabla \ProjAbstr^\ast( p )}{\nabla f_{0}}_{L^2}\nnl 
  &= \cR(f_{\REF} + \ProjAbstr^\ast(p))  +  2\ip{\nabla \ProjAbstr^\ast( p )}{\nabla f_{0}}_{L^2}.  \label{eq:gradL2-proof-1}	
   \end{align}
  We need to show that the mixed term on the right-hand side of \eqref{eq:gradL2-proof-1} vanishes. However, this is a simple consequence of \cref{lem:ProjGradient}: since $\adj \ProjAbstr(p) \in W^{1,2}(\Omega)$ and
  \begin{align}
    \adj \ProjAbstr(p)  = \adj \iota_X \circ \adj \isoProj \circ \adj M (p) = w_{\Proj}^{-1} \cdot w_{\Proj} \cdot \backProj \bparens{  \tuproj^{-1/2} \cdot (\tuproj /\uproj)^{ 1/2 } \cdot p }    = \backProj ( \uproj^{-1/2}  \cdot p ),
  \end{align}
   it follows that $  \uproj^{ 1/2} \cdot \nabla_\mD \big( \uproj^{-1/2}  \cdot p \big) \in  \DetDom $ and we obtain by application of \eqref{eq:lem-ProjGradient-4}
  \begin{align}
   \ip{\nabla \ProjAbstr^\ast( p )}{\nabla f_{0}}_{L^2} &= \ip{\nabla \backProj ( \uproj^{-1/2}  \cdot p )}{\nabla f_{0}}_{L^2}\nnl
   &= \bip{\uproj^{1/2}  \cdot \nabla_{\mD} \parens{\uproj^{-1/2} \cdot p} }{ \uproj^{-1/2} \cdot  \nabla_{\mD} \Proj(f_0) }_{L^2} \stackrel{f_0 \in \Kern(\ProjAbstr) = \Kern(\Proj)}= 0. 
  \end{align}
  Inserting this result into \eqref{eq:gradL2-proof-1} shows that the inequality \eqref{eq:A1} is satisfied for all $ f_0 \in \Kern(\ProjAbstr) \cap W^{1,2}_0(\Omega)$. Moreover, if $f_0 \neq 0$, then $f_0$ is necessarily \emph{non-constrant} in $\Omega$ so that $\norm{\nabla f_0}_{L^2} > 0$. By the second line in \eqref{eq:gradL2-proof-1}, this implies that strict inequality holds in \eqref{eq:A1} if $f_0 \neq 0$.
  Finally, combining \eqref{eq:gradientL2Admissible-pf1} and \eqref{eq:lem-ProjGradient-3} yields
  \begin{align}
   \cR(f_{\REF} + \Proj^\ast(p)) = \bnorm{ \nabla \bparens{ \backProj (\uproj^{-1/2} \cdot p) } }_{L^2}^2 = \bnorm{ \uproj^{1/2} \cdot \nabla_{\mD} \bparens{ \uproj^{-1/2} \cdot p } }_{L^2}^2,
  \end{align}
  which, by \cref{lem:ProjGradient}, remains valid if the expressions attain the value $\infty$.
 \end{proof}
  \vspace{.5em}

\Cref{lem:gradientL2Admissible} enables GenSART-schemes for Kaczmarz-iterations of the form \eqref{eq:KaczmarzGeneric} with $\cR(f) = \norm{\nabla(f-f_{\REF})}_{L^2}^2$ -- if the optimization is restricted to a slightly smaller choice set: 
\begin{align}
 f - f_{\REF} \in W^{1,2}_{\!\Proj} (\Omega) := \bparens{ \closure{\Range(\backProj)} \cap W^{1,2}(\Omega)} \oplus \bparens{ \Kern(\Proj) \cap W^{1,2}_0(\Omega)} \subsetneq W^{1,2} (\Omega)
\end{align}
By \cref{lem:ConvCombPenalties,lem:admissibleQuadratic}, this remains true for more general Sobolev-$W^{1,2}$-like penalties,
\begin{align}
 \cR(f) := \alpha \norm{f - f_{\REF}}_{W^{1,2}_{\gamma, \Proj}}^2 \MTEXT{with}
 \norm{h}_{W^{1,2}_{\gamma, \Proj}}^2 := (1-\gamma) \bnorm{w_\Proj^{1/2} \cdot h}_{L^2}^2 + \gamma \Norm{\nabla h}_{L^2}^2 \label{eq:DefSobolevW12Penalty}
\end{align}
for some $0 \leq \gamma \leq 1$. Note that $h\mapsto \norm{w_\Proj^{1/2} \cdot h}_{L^2}$ is the Hilbert-space-norm of $X$ in \cref{lem:gradientL2Admissible}, which is identical to the $L^2$-norm in the parallel-beam case $\Proj = \ParProj$ but \emph{not} in a cone-beam setting $\Proj = \Cone$, since $w_\ParProj = \boldsymbol 1 _\Omega$ but $w_\Cone \neq \boldsymbol 1 _\Omega$.
The general result reads as follows:

 \vspace{.5em}
  \begin{theorem}[Generalized SART with $W^{1,2}$-penalties] \label{thm:W12SART}
   Let $\Proj \in \{ \ParProj , \Cone \}: L^2(\Omega) \to \DetDom$, $f_{\REF} \in L^2(\Omega)$, $\alpha > 0$ and let $ \cS: \DetDom \to \mRinf$ be any functional. Then the minimizers of
 \begin{align}
  f_{\new} \in \argmin _{ f \in f_{\REF} + W^{1,2}_{\!\Proj} (\Omega)} \cS \left( \Proj ( f) \right) +  \alpha \norm{ f - f_{\REF} }_{W^{1,2}_{\gamma, \Proj}}^2. \label{eq:W12Kaczmarz}
 \end{align}
  are uniquely determined by the GenSART-scheme
   \begin{subequations} \label{eq:W12SART}
  \begin{align}
   p_{\REF} &= \Proj (f_{\REF} ) \label{eq:W12SARTa}  \\
  \Delta p  &\in  \argmin_{p \in \DetDom}   \cS \bparens{  p_{\REF} +   \uproj^{1/2}  \cdot p } +  \alpha (1-\gamma) \Norm{  p}_{L^2}^2 + \alpha \gamma \bnorm{ \uproj^{1/2} \cdot \nabla_{\mD} ( \uproj^{-1/2} \cdot p)}_{L^2}^{2}   \label{eq:W12SARTb}  \\
  f_{\new} &= f_{\REF} + \backProj (\uproj^{-1/2} \cdot \Delta p). \label{eq:W12SARTc} 
  \end{align}
  \end{subequations}
 \end{theorem}
 \vspace{.5em}
 \begin{proof}
  Let $X := ( L^2(\Omega), \ip{\cdot}{\cdot}_\Proj)$  and $P: X \to \DetDom$ be defined as in \cref{lem:gradientL2Admissible} and \smash{$\tilde \cS (p) := \cS( \uproj^{1/2} \cdot p)$}. Then the optimization problem \eqref{eq:W12Kaczmarz} can be written in the form
  \begin{align}
   f_{\new} \in \argmin _{ f \in X } \tilde \cS \left( P(f)  \right) +  \alpha ( 1- \gamma ) \cR_1 (f) + \alpha \gamma \cR_2 (f) . \label{eq:W12SART-pf1}
  \end{align}
  where the functionals $\cR_1, \cR_2: X \to \mRinf$ are given by
  \begin{align}
    \cR_1 (f) := \Norm{ f- f_{\REF}}_{X}^2, \quad \cR_2 (f) := \Norm{\nabla  (f- f_{\REF})}_{L^2}^2 + \chi_{\tilde X}( f- f_{\REF}).
  \end{align}
 $\chi_{\tilde X}: X \to \mRinf$ is the indicator functional of $\tilde X:= W^{1,2}_{\!\Proj} (\Omega)$, defined by $\chi(h) = 0$ if $h \in \tilde X$ and $\chi(h) = \infty$ otherwise. 
 As \smash{$\closure{\Range(\backProj)} = \Range(\adj \ProjAbstr)$} and $\Kern(\Proj) = \Kern(\ProjAbstr)$, it holds that $\tilde X = (\Range(\ProjAbstr) \cap W^{1,2}(\Omega)) \oplus ( \Kern(\ProjAbstr) \cap W_{0}^{1,2}(\Omega) )$.
 By \cref{lem:admissibleQuadratic,lem:gradientL2Admissible}, \cref{A1} is thus satisfied for $\ProjAbstr: X \to \DetDom$ and $\cR = \cR_j$ for $j = 1,2$. According \cref{lem:ConvCombPenalties}, the same is true for $\cR := \alpha (1-\gamma) \cR_1 + \alpha \gamma \cR_2$. Hence, the GenSART-\cref{thm:genSART} is applicable to \eqref{eq:W12SART-pf1}
  so that a minimizers $f_{\new} $ can be found via the scheme \eqref{eq:genSART}:
   \begin{subequations} \label{eq:W12SARTpf3}
     \begin{align}
    \tilde p_{\REF} &= \ProjAbstr (f_{\REF} ) = \uproj^{-1/2}  \cdot \Proj ( f_{\REF} ) \label{eq:W12SARTpf3a}  \\
   \Delta p  &\in \argmin_{p \in \DetDom} \tilde \cS \bparens{  \tilde p_{\REF} + \ProjAbstr \ProjAbstr ^\ast ( p ) } + \cR( f_{\REF}  + \ProjAbstr^\ast ( p ))   \label{eq:W12SARTpf3b}  \\
   f_{\new} &= f_{\REF} + \ProjAbstr^\ast (\Delta p) \label{eq:W12SARTpf3c} 
   \end{align}
   \end{subequations}
  By \cref{lem:admissibleQuadratic,lem:gradientL2Admissible}, the penalty term can be rewritten to
  \begin{align}
   \cR( f_{\REF}  + \ProjAbstr^\ast ( p )) &= \alpha (1-\gamma) \cR_1 ( f_{\REF}  + \ProjAbstr^\ast ( p )) + \alpha \gamma \cR_2 ( f_{\REF}  + \ProjAbstr^\ast ( p )) \nnl
   &= \alpha (1-\gamma) \ip{p}{\ProjAbstr \ProjAbstr^\ast(p)}_{L^2} + \alpha \gamma \bnorm{ \uproj^{1/2} \cdot \nabla_{\mD} ( \uproj ^{-1/2} \cdot p)}_{L^2}^{2}.
  \end{align}
  Moreover, as seen in the proof of \cref{lem:gradientL2Admissible}, it holds that $\adj \ProjAbstr (p ) = \backProj( \uproj^{-1/2} \cdot p)$ and thus \smash{$\ProjAbstr \adj \ProjAbstr (p ) = \uproj^{-1/2} \cdot \Proj  \backProj( \uproj^{-1/2} \cdot p) =   p $} for all $p \in \DetDom$ by \eqref{eq:ProjBackProjRel}.
  Substituting these expressions into \eqref{eq:W12SARTpf3} along with \smash{$p_{\REF} = \uproj^{1/2 } \cdot \tilde p_{\REF} $} and \smash{$\cS(p ) = \tilde \cS( \uproj^{-1/2 } \cdot p)$} yields \eqref{eq:W12SART}.
   \end{proof}
 \vspace{.5em}

 
  \noindent To conclude this section, we make a few remarks on peculiarities of \cref{thm:W12SART}:
 \vspace{.5em}
 \begin{itemize}
  \item The derived SART formula \eqref{eq:W12SART} only involves the \emph{unweighted} unit-projection $\uproj$ and back-projector $\backProj$ and \emph{not} the $L^2$-adjoint $\adj\Proj: p \mapsto w_\Proj \cdot \backProj(p)$. Accordingly, the ray-density-weighting of the back-projection in the cone-beam case is omitted. This is quite intuitive since back-projecting uniformly along the rays results in smaller values of the gradient-penalty functional \eqref{eq:PenaltyGradL2} and thus a non-weighted back-projection can be regarded as the natural one in the considered setting. \vspace{.25em}
  \item In \eqref{eq:W12Kaczmarz}, the increment $\Delta f := f_{\new} - f_{\REF}$ is not optimized over the whole feasible set $W^{1,2}(\Omega)$, but only within the closed subspace $ W^{1,2}_{\!\Proj} (\Omega) \subsetneq W^{1,2}(\Omega)$. This may seem like a fundamental flaw of the result. Notably, however,  $ W^{1,2}_{\!\Proj} (\Omega)$ is \emph{much larger} than the mere space of back-projections $\closure{\Range(\backProj)} \cap W^{1,2}(\Omega) \subsetneq W^{1,2}_{\!\Proj} (\Omega) $, over which the scheme \eqref{eq:W12SART} trivially computes the optimal increment $\Delta f = \backProj(\uproj^{-1/2}\cdot \Delta p)$. In this sense, \cref{thm:W12SART} still provides a non-trivial simplification of the optimization problem in \eqref{eq:W12Kaczmarz}.
  Indeed, $W^{1,2}_{\!\Proj} (\Omega)$ is arguably \emph{almost} as large as $W^{1,2}(\Omega)$: \smash{$\closure{\Range(\backProj)} \cap W^{1,2}(\Omega)$} contains all functions that are constant along the rays and  $\Kern(\Proj) \cap W_{0}^{1,2}(\Omega)$ all those, which have zero mean along these and vanish on the boundary of $\Omega$. Accordingly, the ``missing subspace'', \IE the orthogonal complement of $W^{1,2}_{\!\Proj} (\Omega)$ within $W^{1,2}(\Omega)$, must be composed solely of functions in $\Kern(\Proj)$ that are \emph{linear} along all ray-segments. In practice, the GenSART-iterates defined by \eqref{eq:W12SART} are thus expected to provide \emph{almost} the optimum of the objective in \eqref{eq:W12Kaczmarz} over all $f\in W^{1,2}(\Omega)$ (or $f \in L^2(\Omega)$).
 \end{itemize}
 \vspace{.5em}

 \subsection{$L^q$-penalties} \label{SS:LqPenalties}

The aim of this section is to demonstrate that \cref{A1} does not restrict the choice of penalty functionals $\cR$ to quadratic ones. We consider $L^q$-penalties
\begin{align}
 \cR(f) &:= \norm{f - f_{\REF}}_{L^q}^q \MTEXT{with} \norm{h}_{L^q}^q := \int_{\mR^m} |h(\bx)|^q \, \D \bx \in \mRinf, \quad 1 \leq q < \infty \label{eq:RkLq}
\end{align}
 and define $L^q(\Omega) := \{ f : \norm{f}_{L^q} < \infty \}$ as usual.
 We prove the admissibility of such penalties only for a parallel-beam setting $\Proj = \ParProj$. An extension to the cone-beam case may be possible.

\vspace{.5em}
  \begin{lemma}[Admissibility of $L^q$-penalties] \label{thm:admissibleLq}
 Let $\ProjAbstr = \isoParProj : L^2(\Omega) \to \DetDom$ and let $\cR: L^2(\Omega) \to \mRinf$ be defined by \eqref{eq:RkLq}. Then \cref{A1} is satisfied and
 \begin{align}
   \cR(f_{\REF} + \ProjAbstr^\ast(p)) = \begin{cases} 
                                            \bnorm{ \uparproj ^{1/q-1/2} \cdot p }_{L^q}^q &\text{if }  \ProjAbstr^\ast(p) \in L^q(\Omega) \\
                                            \infty &\text{else}
                                           \end{cases} \MTEXT{for all} p\in \DetDom.
 \end{align}
 \end{lemma}
\vspace{.5em}

 The proof of \cref{thm:admissibleLq} is given in appendix~\ref{Appendix:LqPenaltyProof}.
 For completeness, we state the GenSART-scheme that is obtained by applying \cref{thm:genSART} to the setting in \cref{thm:admissibleLq}:
 
 \vspace{.5em}
  \begin{theorem}[Generalized SART with $L^q$-penalties] \label{cor:LqSART}
 Let  $f_{\REF} \in L^2(\Omega)$, $\alpha > 0$ and let $ \cS: \DetDom \to \mRinf$ be any functional. Assume that there exists a minimizer  
 \begin{align}
  f_{\new} \in \argmin _{ f \in L^2(\Omega) }  \cS \left( \ParProj ( f) \right) +  \alpha \norm{f - f_{\REF}}_{L^q}^q. \label{eq:LqKaczmarz}
 \end{align}
 Then any $ \tilde f_{\new}$ determined by the GenSART-scheme
   \begin{subequations} \label{eq:LqSART}
  \begin{align}
   p_{\REF} &= \ParProj (f_{\REF} ) \label{eq:LqSART-a}  \\
  \Delta p  &\in \argmin_{p \in \DetDom} \cS \bparens{ p_{\REF} + \uparproj^{1/2} \cdot p } +  \alpha  \bnorm{ \uparproj ^{1/q - 1/2} \cdot p }_{L^q}^q   \label{eq:LqSART-b}  \\
  \tilde f_{\new} &= f_{\REF} +  \ParProj^\ast \bparens{ \uparproj^{-1/2} \cdot  \Delta p } \label{eq:LqSART-c} 
  \end{align}
  \end{subequations}
  also minimizes \eqref{eq:LqKaczmarz}.
  If $q > 1$, then any minimizer $f_{\new}$ of \eqref{eq:LqKaczmarz} is of the form \eqref{eq:LqSART}.
 \end{theorem}
 \vspace{.5em}

\optspace
\section{Applications} \label{S:Applications}

In the preceding sections, it has been analyzed in which abstract situations Kaczmarz-iterations of the form \eqref{eq:genKaczmarz} can be computed via a generalized SART-scheme. In the following, the principal theory is applied to design tailored methods for various settings of tomographic imaging. Specifically, the aim is to exploit the extraordinary freedom that the GenSART-\cref{thm:genSART} offers in choosing the data-fidelity functionals $\cS_k $ and image-formation operators $F_j$.
Despite differences between the specific applications, it should be emphasized that all of the proposed methods are applicable for both parallel- and cone-beam acquisition-geometries, without any requirements on the incident directions or source positions.  


\subsection{Noise-model-adapted GenSART} \label{SS:NoiseAdaptedSART}

As outlined in \secref{SS:RecMethComplexity}, variational- and Kaczmarz-type reconstruction methods may account for the expected statistics of the data errors $\bepsilon$ in \ipref{ip1} by suitably choosing the data-fidelity functionals $\cS_k$ in \eqref{eq:genKaczmarz}. 
We illustrate this for Kaczmarz-iterations with a simple $L^2$-penalty and fixed $\cS_k = \cS$:
\begin{align}
 f_{k+1} &\in \argmin_{f \in L^2(\Omega)} \cS \left( g_{j_k}^\obs; \,  F_{j_k} \left( \Proj_{j_k} ( f) \right)   \right) +  \alpha \norm{ f - f_k }_{L^2}^2  \label{eq:NoiseAdaptKaczmarz}
 \end{align}
By \cref{thm:genL2SART}, the minimizer can be computed via the SART-like scheme
\begin{subequations} \label{eq:NoiseAdaptSART}
  \begin{align}
  p_k &= \Proj_{j_k} ( f_k ) \label{eq:NoiseAdaptSART-a} \\
  \Delta p_k &\in  \argmin_{p \in \DetDomj{j_k}}   \cS \left( g_{j_k}^\obs; \, F_{j_k} \left(  p_k  + \tuprojj{j_k}^{1/2} \cdot p  \right) \right) + \alpha \norm{p}_{L^2}^2 \label{eq:NoiseAdaptSART-b} \\
 f_{k+1}  &= f_k + \Proj_{j_k}^\ast \left( \tuprojj{j_k}^{-1/2} \cdot \Delta p_k  \right)  \label{eq:NoiseAdaptSART-c}
\end{align}
\end{subequations}
where $\DPj{j}$ and $\tuprojj{j}$ denote the projection-domain and weighted unit-projection of the projector $\Proj_{j}$, respectively, see \secref{SS:ProjProp}.
According to Bayesian theory, the Kaczmarz-iterations in \eqref{eq:NoiseAdaptKaczmarz} can be tailored for a specific (probabilistic) model of the data errors $\bepsilon$ by choosing $\cS$ as the negative log-likelihood of the fitted data $g_{j_k} =  F_{j_k} \left( \Proj_{j_k} ( f) \right) $ given the observations $g_{j_k}^\obs$:
\begin{align}
  \cS \bparens{ g_{j}^\obs; \, g_j } := - \ln \mathbb P ( g_{j}^\obs | g_j ), \label{eq:LogLikelihood}
\end{align}
where $\mathbb P ( g_{j}^\obs | g_j )$ denotes the probability of measuring $g_{j}^\obs$ given that the true data is $g_j$. By inserting \eqref{eq:LogLikelihood} into \eqref{eq:NoiseAdaptSART}, generic noise-model-adapted GenSART-schemes are obtained.

\paragraph{Efficient closed-form optimization in projection-space:}
For general noise-models and image-formation operators $F_j$, the optimization problem in  \eqref{eq:NoiseAdaptSART-b} could still be hard to solve, in spite of being cast to the low-dimensional projection-space via the GenSART-approach. In the following, we therefore outline practically relevant settings where the optimization-step in the GenSART-scheme \eqref{eq:NoiseAdaptSART} may be performed at negligible computational costs.

Often, the data-errors $\bepsilon$ at different spatial positions can be assumed to be stochastically independent. 
If the $F_{j}$ are pointwise operators, i.e.\ ``$F_{j_k}(p)(x) = F_{j_k}(p(x))$'' for all $p$ and $x$ (in particular if  $F_j = \identity$), the data-fidelity obtained via \eqref{eq:LogLikelihood} is then of \emph{integral-form}:
\begin{align}
  \cS \bparens{ g_{j}^\obs; \,F_j(p) } = \int_{\DPj j } s_j(x, p(x)) \, \D x + c  \MTEXT{for some} s_j: \DPj j \times \mR \to \mRinf \label{eq:DataFidIntegral} 
\end{align}
with some additive constant $c\in \mR$ that does not affect the minimizer in \eqref{eq:NoiseAdaptKaczmarz}. 
By substituting \eqref{eq:DataFidIntegral} into the the objective-functional to be minimized in \eqref{eq:NoiseAdaptSART}, 
we obtain
\begin{align}
 \cS &\Parens{ g_{j_k}^\obs; \,  F_{j_k} \bparens{ p_k  + \tuprojj{j_k}^{1/2}\cdot p } } + \alpha \norm{p}_{L^2}^2 \nnl
 &=  \int_{\DPj{j_k}} \Big( s_{j_k}\big(x, p_k(x) +  \tuprojj{j_k}(x)^{1/2}  p(x)  \big)+ \alpha  |p(x)|^2 \Big) \, \D x    \label{eq:SARTCostFunctionalPointwise}
\end{align}
The integrand in \eqref{eq:SARTCostFunctionalPointwise} depends only on point-evaluations of $p$, i.e.\ only on the local function value. As a consequence, the optimization  in \eqref{eq:NoiseAdaptSART-b} is equivalent to a family of \emph{scalar} problems:
\begin{align}
 \Delta p _k &\in  \argmin_{p \in \DetDomj{j_k}}  \cS \bparens{ g_{j_k}^\obs; \,  F_{j_k} \bparens{ p_k  + \tuprojj{j_k}^{1/2}\cdot p } } + \alpha \norm{p}_{L^2}^2  \nnl
 \Leftrightarrow \; \Delta p _k(x) &\in  \argmin_{y \in \mR} s_{j_k}\big(x, p_k (x) + \tuprojj{j_k}(x)^{1/2} y  \big) + \alpha  y^2 \nnl
                                       &= \tuprojj{j_k}(x)^{-1/2} \cdot \Big(  \prox\left(s_{j_k}(x, \cdot)\right)\left(p_k (x), 2\tuprojj{j_k}(x)^{1/2} /\alpha\right) - p_k (x) \Big)  \label{eq:SARTScalarOptProblems}
\end{align}
for almost all $x \in \DPj {j_k}$. Here, the usual \emph{proximal map} of a functional has been introduced:
\begin{align}
 \prox(s)(y, \tau) := \argmin_{x \in \mR} s(x) + \frac{(x-y)^2}{2 \tau}   \label{eq:DefProx}
\end{align}

The prox in \eqref{eq:SARTScalarOptProblems} can be typically evaluated numerically in $\cO(1)$ floating-point operations. Hence, a discretized form of the optimization problem in \eqref{eq:NoiseAdaptSART} can be solved in $\cO(\Mpx)$, where $\Mpx$ is the number of degrees-of-freedom of a discretized projection. This enables evaluations of \eqref{eq:NoiseAdaptSART} at literally the same computational costs as classical SART-iterations, compare \secref{SS:L2SART-example} -- even for highly non-trivial choices of $\cS$, as demonstrated by the subsequent examples.

\subsubsection{(Weighted) $L^2$-fidelities}  \label{SSS:L2DataFid}

For completeness, we mention the case of (weighted) $L^2$-data fidelities $\cS \bparens{ g_{j}^\obs; \, F_j(p) } = \norm{ ( F_j(p) - g_{j}^\obs)/ \sigma_j }_{L^2}^2$, which are adapted to data errors  caused by Gaussian white noise of possibly spatially varying variance $\sigma^2_j$. For $F_j = \identity$, this choice of $\cS$ is of the integral form \eqref{eq:DataFidIntegral} with $s_j(x, y) := (y - g_{j}^\obs(x))^2 / \sigma_j(x)^2$ and  a simple proximal map,   $\prox\left(s_{j}(x, \cdot)\right)\left(y, \tau \right) = ( 2 \sigma_j(x)^2 y  + \tau g_{j}^\obs(x) ) / ( 2 \sigma_j(x)^2 + \tau )$. 

\subsubsection{Robust GenSART}  \label{SSS:RobustDataFid}

As a first non-standard application, we consider the problem of robust tomographic reconstruction: systematic errors in the acquisition geometry or modeling-inaccuracies due to nonlinear effects, as arising from metal-inclusions in soft tissue for example \cite{BarrettEtAl2004ArtifactsInCT}, tend to produce large outliers in the data, i.e.\ errors with highly non-Gaussian statistics.
In such a setting, an $L^2$-data-fidelity  is far too greedy in fitting the data. 

The problem has been addressed by employing a more robust (Huber-)$L^1$-term, see e.g.\  \cite{AndersenHansen2014_genARTProx}, or even non-convex data-fidelity functionals such as the negative log-likelihood of the Student's t-distribution, as proposed in \cite{Bleichrodt2015_STTomo}. For trivial image-formation maps $F_j = \identity$, these choices correspond to a data-term $\cS ( g_j^\obs; \, F_j(p) )$ of the integral-form \eqref{eq:DataFidIntegral} with
\begin{subequations} \label{eq:DataFidRobust}
 \begin{align}
   s_j( x, y) &= s_-(y - g_{j}^\obs), \qquad  s_- \in  \{  s_{L^1_{\textup H}, \nu }, s_{\textup{s-t}, \nu}\} \label{eq:DataFidRobust-a} \\
   s_{L^1_{\textup H}, \nu }(y) &:= \begin{cases}
                              |y|^2 &\textup{if } |y| \leq \nu \\
                              2\nu|y| - \nu^ 2 &\textup{else}
                            \end{cases}, \qquad  s_{\textup{s-t}, \nu}(y) := \nu^2 \ln( 1 + |y|^2 / \nu^ 2 )   \label{eq:DataFidRobust-b}
\end{align}
\end{subequations} 

Both functions in \eqref{eq:DataFidRobust-b} show the same quadratic growth behavior as $y \mapsto |y|^2$ for $|y| \ll \nu$ whereas,  for $|y| \geq \nu$, the growth is only linear in the $L^1$-Huber- and  logarithmic in the Student's-t-case. Accordingly, the resulting data fidelities $\cS$ defined by \eqref{eq:DataFidRobust-a} and \eqref{eq:DataFidIntegral} behave like an $L^2$-term for small deviations between the measured and fitted data, but penalize much less strongly when the deviations are larger, thereby yielding an increased robustness towards outliers.
%
The proximal maps of $ s_{L^1_{\textup H}, \nu }$, $s_{\textup{s-t}, \nu}$ are given by
\begin{subequations} \label{eq:ProxL1HST}
 \begin{align}
 \prox( s_{L^1_{\textup H}, \nu })(y, \tau)&= y - \frac{2\nu\tau y}{\max\{|y|, 2\nu\tau + 1\}}  \label{eq:ProxL1H}\\
 \prox(s_{\textup{s-t}, \nu})(y, \tau)&= \argmin\!\left\{ s_{\textup{s-t}, \nu}(x_0): x_0\in \mR, \,  (x_0^2 + \nu^2) (x_0-y) + 2\tau \nu^ 2 x_0 = 0\right\} \label{eq:ProxST}
\end{align}
\end{subequations}
The cubic equation in \eqref{eq:ProxST} has one or three real roots $x_0$, which can be computed analytically in $\cO(1)$. In the case of three solutions, $\argmin\left\{ s_{\textup{s-t}, \nu}(x_0) \right\}$ may be determined among these by simple trial-and-error. Consequently, the GenSART-scheme \eqref{eq:NoiseAdaptSART} can be evaluated efficiently via \eqref{eq:SARTScalarOptProblems} for the considered robust data-fidelity terms. Notably, the simplicity of the $L^1$-Huber-prox has already been exploited in \cite{AndersenHansen2014_genARTProx} to construct efficient robust \emph{ART}-iterations.

\subsubsection{Poisson-noise-adapted GenSART} \label{SSS:PoiNoise}

In many practical applications of X-ray- or electron tomography, the data errors are primarily due to the \emph{Poisson-statistics} of the detection process: detector pixels actually count a \emph{discrete} number of incident photons or electrons over some exposure time $t>0$, where the counts follow a Poisson-distribution. Disregarding this effect may lead to severe anisotropic noise in the reconstruction, for example when a sample in CT is so strongly absorbing along certain incident directions that only very few counts are detected (\emph{photon starvation}, see e.g.\ \cite{Mori2013PhotonStarvation}).

If the detector is composed of $\Mpx \in \mN$ pixels with spatially varying detection sensitivities $\omega_i: \mD \to \mR_{\geq 0}$, the measured data is given by a vector 
\begin{align}
 g^{\textup{obs}}_j = ( g^{\textup{obs}}_{ji} )_{i = 1}^{  \Mpx} \in \mR^{ \Mpx}, \quad g^{\textup{obs}}_{ji} \sim \textup{Poi}\left( t g_{ji} \right), \quad g_{ji} = \cM_i(g_j) := \int_{\mD} \omega_i g_j \, \D x \label{eq:PoissonData} 
\end{align}
where $g_j$ denotes the exact continuous data and $X \sim \textup{Poi}(\lambda)$ means that $X$ is a Poisson-distributed random variable of intensity $\lambda \geq 0$. In this setting, the log-likelihood in \eqref{eq:LogLikelihood} leads to the discrete Kullback-Leibler-divergence, see \cite{HohageWerner2016PoissonReview} for details:
\begin{align}
  \cS^{\textup{Poi}} \left( g^{\textup{obs}}_j;  g_j \right) &:=  \sum_{ i = 1}^{ \Mpx} \kl (g^{\obs}_{ji}; t \cM_i \left( g_j )  \right) ), \quad \kl( b; a) := \begin{cases}
                                                     a - b - b \ln \left( \frac{a}{b} \right) & a,b \geq 0 \\
                                                     \infty &\text{else}
                                                    \end{cases} \label{eq:KLDataFid}
\end{align}
with the conventions that $\ln(0) = - \infty$ and $0\cdot \ln(a/0) = 0$ for all $a \geq 0$.

Under the assumption that variations of the true data $g_j$ are negligible within the support of $\omega_i$ (i.e.\ within a single pixel!), $\cS^{\textup{Poi}}$ can be approximated in the form
\begin{align}
  \cS^{\textup{Poi}}  \left( g^{\textup{obs}}_j;  g_j \right) \approx \int_{\mD} \kl \left(g^{\obs}_{j, \textup{cont}}(x) ; t g_j(x)  \right) \cdot \omega(x)  \, \D x + c  =: \cS( g^{\obs}_{j }  ; g_j ) \label{eq:KLDataFidCont}
\end{align}
where  $\omega := \sum_{i = 1}^\Mpx \omega_i$, $g^{\obs}_{j, \textup{cont}}(x):= \big( \sum_{i = 1}^\Mpx g^{\obs}_{ji} \omega_i(x) / \int_{\mD}  \omega_i \, \D x \big) /  \omega(x)$ (with the convention that $0/0 = 0$) and  $c$ is independent of $g_j$. The derivation is given in appendix~\ref{Appendix:PoiDataFid}.

So far, we have not specified the image-formation operators $F_j$, relating the data $g_j$ to the tomographic projections $\Proj_j(f)$. We consider two different acquisition modes:
\vspace{.5em}
\begin{itemize}
 \item \emph{Dark-field imaging:} The exact data is directly proportional to the projections, i.e.\ $F_j(p_j) = I_j \cdot p_j$ where $I_j$ is the illumination intensity. This applies for example to  HAADF-STEM, a state-of-the-art electron tomography technique, see e.g.\  \cite{MidgleyWeyland2002_ETomoContrastMechanisms,Oektem2015MathETomo}.  \vspace{.25em}
 \item \emph{Bright-field imaging:} The tomographic data gives the relative attenuation experienced by the illuminating beam, i.e.\ $F_j(p_j) = I_j \cdot \exp( - p_j )$.  This is the model for classical (monochromatic) X-ray computed tomography.
\end{itemize}
\vspace{.5em}
Inserting these models into \eqref{eq:KLDataFidCont}, it can be seen that the resulting data-term $\cS( g^{\obs}_{j }; \, F_j(p) )$ is of the integral-form \eqref{eq:DataFidIntegral} with
\begin{align}
   s_j(x,y) = 
\begin{cases} 
 \kl \left(g^{\obs}_{j, \textup{cont}}(x) ;\, t I_j(x)\cdot  y  \right) &\textup{(dark-field)} \\
 t I_j(x) \exp(-y)  + g^{\obs}_{j, \textup{cont}}(x) \cdot y  &\textup{(bright-field)}
\end{cases}.
\end{align}
In the dark-field case, $\prox ( s_j(x, \cdot) )$ has a well-known closed form, see e.g.\ \cite{HohageWerner2016PoissonReview}. In the bright-field case, $\prox ( s_j(x, \cdot) )$  may be evaluated numerically by a few iterations of Newton's method. This admits an efficient implementation of Poisson-noise-adapted GenSART-schemes.


\subsection{Regularized Newton-Kaczmarz-GenSART} \label{SS:NewtonKaczmarz}

Regularized Newton-Kaczmarz  methods have been proposed in \cite{Burger2006NewtonKaczmarz} for the solution of general block-structured inverse problems $G(f) = (G_1(f), \ldots, G_N(f)) = (g_1^\obs, \ldots, g_N^\obs)$ with \emph{nonlinear} forward operators $G_j: X \to Y_j$ between Hilbert spaces $X, Y_1, \ldots, Y_N$. In its simplest form, the approach boils down to performing Levenberg-Marquardt iterations on the different sub-problems $G_j(f) = g_j^\obs$:
\begin{align}
 f_{k+1} = \argmin_{f \in X} \Norm{G_{j_k}(f_k) + G_{j_k}'[f_k] ( f - f_k) - g _{j_k}^ \obs }_{Y_j}^2 + \alpha \Norm{ f - f_k}_{X}^2  \label{eq:NewtonKaczmarz}
\end{align}
where $G_j'[f_k] : X \to Y_j$ denotes the Fr\'echet-derivative of the operator at $f_k$.
In the following, two nonlinear tomographic reconstruction problems are presented for which the iterations \eqref{eq:NewtonKaczmarz} can be computed efficiently via generalized SART-schemes.

\subsubsection{Propagation-based X-ray phase contrast tomography}  \label{SSS:PCTSART}

We consider the setting of \mbox{(propagation-based)}  {X-ray phase contrast tomography} (XPCT), see e.g.\ \cite{Cloetens1999holotomography,Bartels2015,Krenkel2014BCAandCTF,Ruhlandt2016_RadonPCICommute,MaretzkeEtAl2016OptExpr}.
In this experimental setup, the recorded data is given by near-field diffraction patterns, that relate to tomographic projections of the object density via a highly non-trivial image-formation operator $F_j = F $: under the standard assumptions of an ideal, fully coherent X-ray beam and negligible absorption (often a very good approximation at high X-ray energies), the measured parallel-beam(!)  data under the $j$th tomographic incident direction is modeled by
\begin{align}
 g_{j} := F  \left( \ParProj_j (f) \right) \MTEXT{with} F(p) := \left| \cD \left( \exp\left( - \I  p  \right) \right) \right|^2 -1.  \label{eq:PCTNLModel}
\end{align}
Here, $| \cdot |^2$ denotes the pointwise squared modulus of a complex-valued field and $\cD$ is the \emph{Fresnel propagator}, which is given by a unitary Fourier-multiplier ($\cF$: Fourier transform):
\begin{eqnarray}
 \cD ( \psi ) :=   \cF^{-1}\left( \mF \cdot \cF(\psi) \right), \quad \mF(\bxi )  := \exp \Parens{ -{ \I \bxi^2 }/\parens{ 4 \pi \FAlt } } \MTEXT{for} \bxi \in \mR^m, \label{eq:FresnelPropDef}
\end{eqnarray}
where $\F$ is the Fresnel-number of the imaging setup.
The example matches the Newton-Kaczmarz setting if we define $G_j(f)  := F  \left( \ParProj_j (f) \right)$. The Fr\'echet-derivative is given by
\begin{align}
 G_j[f]h &= F'[\ParProj_j (f)]  \ParProj_j(h), \quad F'[p]h_p = 2\imag \left( \closure{\cD \left(  \exp\left( - \I p  \right) \right) } \cdot \cD \left( \exp\left( - \I p  \right) \cdot h_p \right) \right) \label{eq:PCTDerivative}
\end{align}
where the overbar denotes complex conjugation and $\imag$ the pointwise imaginary part.  Newton-Kaczmarz iterations for this problem with $L^2$-data-fidelity and Sobolev-$W^{1,2}$-penalty, as first proposed in \cite{MaretzkeEtAl2016OptExpr}, are  of the form 
\begin{align}
 f_{k+1} =   \argmin_{f \in L^2(\Omega) } \bnorm{F  \left( \ParProj_{j_k} (f_k) \right) + F'[\ParProj_{j_k} (f_k) ] \left( \ParProj_{j_k}(f-f_k) \right) - g_{j_k}^{\textup{obs}}}_{L^2}^2  \nnl
 +  \alpha \left( (1-\gamma)  \norm{f-f_k}_{L^2}^2  + \gamma  \norm{\nabla( f-f_k )}_{L^2}^2\right). \label{eq:NewtonKaczmarzPCT}
\end{align}
If we take $\cS(p):=  \norm{F  \left( \ParProj_{j_k} (f_k) \right) + F'[\ParProj_{j_k} (f_k) ] \left( p - \ParProj_{j_k} (f_k) \right) - g_{j_k}^{\textup{obs}}}_{L^2}^2$,   
\eqref{eq:NewtonKaczmarzPCT} matches the setting of \cref{thm:W12SART}. Hence, within the minor approximation discussed in the end of \secref{SS:GradPenalties}, the minimizer may be computed via the GenSART-scheme
\begin{align}
 f_{k+1} \approx f_k + \ParProj_{j_k}^\ast \bigg( \uprojj{j_k}^{-1/2} \cdot  \argmin_{p \in \DetDomj{j_k}} \bnorm{ F'[\ParProj_{j_k} (f_k) ]\bparens{ \uprojj{j_k}^{1/2} \cdot p } - r_{j_k}  }_{L^2}^2 \nnl
 + \alpha_k  (1-\gamma) \norm{  p }_{L^2}^2 + \alpha_k  \gamma \bnorm{ \uprojj{j_k}^{1/2} \cdot \nabla_{\mD} \bparens{ \uprojj{j_k}^{-1/2 } \cdot p }}_{L^2}^2  \bigg)   \label{eq:PCTSART} 
\end{align}
with residual $r_{j_k} := g_{j_k}^{\textup{obs}} - F  \left( \ParProj_{j_k} (f_k) \right)$. 
The quadratic optimization problem in \eqref{eq:PCTSART} 
can be solved for example by a conjugate-gradient method applied to the normal equation.

\subsubsection{Polychromatic CT}  \label{SSS:PolyCTSART}

If the polychromatic nature of the X-rays in conventional CT-scanners is neglected, so called \emph{beam-hardening} artifacts may arise \cite{BarrettEtAl2004ArtifactsInCT}.
In \cite{DeMan2001_PolyCTModel,Humphries2014_PolyCTModel}, a simplified model for polychromatic CT has been proposed, which partially accounts for the arising nonlinear effects. Within this model, the detected intensity data $g_j$ for the $j$th tomographic projection is predicted as
\begin{align}
 g_j = \int \underbrace{ I_0(\varepsilon)  \exp \Big( -  \Phi(\varepsilon) \Proj_j \big(  \phi (f) \big) -  \Theta(\varepsilon)    \Proj_j \big( \theta (f) \big)  \Big) }_{=:G_{j, \varepsilon}(f)} \, \D \varepsilon =: G_{j} (f).   \label{eq:PolyCTModel}
\end{align}
Here, $I_0(\varepsilon)$ is the emitted intensity of the X-ray source at photon-energy $\varepsilon$
and $f$ is the spatially varying attenuation at some reference-energy $\varepsilon_0$. $ \Phi(\varepsilon)   \phi(f)$ and $ \Theta(\varepsilon)  \theta(f)$ model the attenuation's photo-electric- and Compton-scattering-components, respectively. The main approximation is that these functions are assumed to be representable as a function of $f$ multiplied by  energy-dependent scaling factors: 
\begin{align}
 \Phi(\varepsilon) = \frac{ \varepsilon _0^3 }{\varepsilon^3} \MTEXT{and} \Theta(\varepsilon) = \frac{ f_{\textup{KN}} ( \varepsilon ) }{ f_{\textup{KN}} ( \varepsilon_0 ) }, \qquad \text{$f_{\textup{KN}}$: Klein-Nishina function} \label{eq:PolyCTEnergyDep}
\end{align}
Note that the expressions $\phi (f)$ and $\theta (f)$ are to be understood pointwise, \IE $\phi (f)(x) = \phi(f(x))$ and $\theta (f)(x) = \theta(f(x))$ (with a slight abuse of notation). The scalar functions $\phi, \theta: \mR_{\geq 0} \to \mR_{\geq 0}$ interpolate known value-pairs $\{(f_m, \phi_m)\}$ and $\{(f_m, \theta_m)\}$ for different materials $m = 1,2,\ldots$ in the imaged object, such as water, fat and bone, see \cite{DeMan2001_PolyCTModel,Humphries2014_PolyCTModel} for details. We assume that $\phi$ and $\theta$ are continuously differentiable.

Notably, the nonlinearity in \eqref{eq:PolyCTModel} \emph{cannot} be described by a nonlinear map acting on the projections, i.e.\ $G_j(f) \neq F_j(\Proj_j(f))$, so that the setting does not seem to match our tomographic model \eqref{eq:TomoModel}.
However, we may still compute the Fr\'echet-derivative:
\begin{subequations}
 \begin{align}\label{eq:PolyCTFrechetDeriv}
 G_j'[f]h_f &= -  G_j^\Phi(f)  \cdot \Proj_j\left(  \phi' (f) \cdot h_f \right) -  G_j^\Theta(f)  \cdot  \Proj_j\left(  \theta' (f) \cdot h_f \right) \\
G_j^\Phi(f) &:= \int \Phi(\varepsilon) G_{j, \varepsilon }(f) \, \D \varepsilon , \quad  G_j^\Theta(f) :=   \int \Theta(\varepsilon) G_{j, \varepsilon}(f) \, \D \varepsilon 
 \end{align}
\end{subequations}
By \eqref{eq:ProjBackProjRel}, it holds that $ G_j^\Phi(f)  \cdot \Proj_j \left(  \phi' (f) \cdot h_f \right) =   \Proj_j\bparens{ \phi' (f) \cdot h_f \cdot \backProj_j(  G_j^\Phi(f)) }$, which yields
\begin{subequations} \label{eq:PolyCTFrechetDeriv2}
  \begin{align}
  G_j'[f]h_f 
  &= \Proj_j( \lambda_j(f)  \cdot h_{f} ), \quad \lambda_j(f) = \phi' (f) \cdot \backProj_j \left( G_j^\Phi(f) \right)  +  \theta' (f) \cdot  \backProj_j \left(  G_j^\Theta(f) \right) 
  \label{eq:PolyCTFrechetDeriv2-a}
 \end{align}
\end{subequations}
 Using \eqref{eq:PolyCTFrechetDeriv2}, Newton-Kaczmarz iterations \eqref{eq:NewtonKaczmarz} for the considered problem with an $L^2$-penalty can be written in the form
\begin{align}
 f_{k+1} &=   \argmin_{f \in L^2(\Omega) } \Norm{  \Proj_{j_k} \left( \lambda_{j_k}(f_k)  \cdot (f-f_k) \right)   - r_{j_k}}_{L^2}^2  + \alpha_k  \Norm{f-f_k}_{L^2}^2 \label{eq:NewtonKaczmarzPolyCT} 
\end{align}
with  residual $r_{j_k}= g_{j_k}^{\textup{obs}} - G_{j_k} (f_k)$. \eqref{eq:NewtonKaczmarzPolyCT} matches the setting of \cref{thm:L2SARTProjW}. By rearranging the resulting GenSART-scheme \eqref{eq:L2SARTProjW} for \eqref{eq:NewtonKaczmarzPolyCT}, we obtain the update-formula
\begin{align}
   f_{k+1} 
 &= f_k + \lambda_{j_k}(f_k)  \cdot \adj{\Proj_{j_k}} \bigg( \frac{ r_{j_k} }{ \Proj_{j_k} \big( \wprojj{j_k} \cdot |\lambda_{j_k}(f_k)|^2  \big)  + \alpha_k  }  \bigg) \label{eq:NewtonKaczmarzPolyCTSART} 
\end{align}
where $\wprojj{j_k}$ denotes the ray-density to the projector $\Proj_{j_k}$.
Including the necessary computations of $G_{j_k} (f_k)$ and $\lambda_{j_k}(f_k)$, evaluating  \eqref{eq:NewtonKaczmarzPolyCTSART} requires three evaluations of the forward- and back-projectors $\Proj_{j_k}$ and $\adj{\Proj_{j_k}}$, plus computationally inexpensive pointwise operations. 

Similarly efficient formulas may be obtained if the $L^2$-data-fidelity in \eqref{eq:NewtonKaczmarzPolyCT} is replaced by the Poisson-noise-adapted Kullback-Leibler term from \secref{SSS:PoiNoise}.


\subsection{Extensions} \label{SS:Extensions}

The following section outlines different ideas on how to extend the generalized SART-approach to  an even more versatile tool for devising efficient Kaczmarz-type reconstruction methods.

\subsubsection{Box constraints} \label{SSS:BoxConstraints}

Analogously as in other Kaczmarz-methods, box constraints $f_{\min} \leq f \leq f_{\max}$ on the admissible values of the object $f$ may be incorporated in GenSART-schemes simply by setting
\begin{align}
 f_{k+1}  \leftarrow \max \big\{ \min \{ f_{k+1}, f_{\max} \}, f_{\min} \big\}.
\end{align}
after each iteration. This approach can be interpreted as interlacing \emph{projections} onto the convex set $\{f \in L^2(\Omega): f_{\min} \leq f \leq f_{\max}\}$ 
and is standard in Kaczmarz-type tomographic reconstructions ever since the introduction of ART \cite{Gordon1970ART}.

\subsubsection{Additional quadratic regularizer} \label{SSS:AddQuadReg}

So far, the penalization in the considered Kaczmarz-iterations was always with respect to the preceding iterate $f_k$. In addition, it might be desirable to impose a static regularizer such that the total penalty  is given by
\begin{align}
 \cR_k(f) =  \alpha_1 \norm{T(f - f_k)}_Z^2 + \alpha_2 \norm{T(f - f_{\REF})}_Z^2,
\end{align}
for $\alpha_1,\alpha_2 > 0$, Hilbert spaces $X,Z$ and a bounded linear operator $T: X \to Z$. By writing the squared norms as inner products, the expression can be cast to the form
\begin{align}
 \cR_k(f) =  (\alpha_1 + \alpha_2) \norm{T(f - f_{k,\REF})}_Z^2 + c
\end{align}
where $f_{k,\REF} := ( \alpha_1 f_k + \alpha_2 f_{\REF}) / (\alpha_2 + \alpha_1)$ and the constant $c \in \mR$ is independent of $f$ and thus irrelevant for the computation of the minimizer. Hence, the resulting Kaczmarz-iterations are of the same form as before, up to a modified reference solution $f_k \to f_{k,\REF}$ and regularization parameter $\alpha \to \alpha_1 + \alpha_2$. This means that also the derived GenSART-formulas only have to be modified by exchanging these parameters.

\subsubsection{Kaczmarz-type splitting and primal-dual methods} \label{SSS:Interlacing} 

Variational reconstruction methods seek to minimize terms of the form $ \cS_\tot\left( \Proj_\tot (f)  \right) + \cR(f)$ with $\cS_\tot(p_1, \ldots, p_\Nproj) = \sum_{j=1}^\Nproj \cS(p_j)$, compare \secref{SS:RecMethComplexity}. Often, this is achieved by some \emph{splitting} method, alternating (sub-)gradient-descent- (``forward-'') or proximal (``backward-'') iterations with respect to the data-fidelity $\cS_\tot$ and the penalty functional $\cR$ (or their duals), see e.g.\ \cite{CombettesEtAl2011ProximalSplitting}. It is straightforward to combine such an approach with a Kaczmarz-type strategy that exploits the block-structure of $\cS_\tot$ to reduce computational costs of the individual iterations. Examples of (primal) Kaczmarz-type splitting methods are given by iterations ($\partial \cR$: subdifferential)
\begin{subequations} \label{eq:FWBWSplit}
 \begin{align}
 f_{k+\frac 1 2} &\in \argmin_{f \in L^2(\Omega) } \cS_{j_k}\left( \Proj_{j_k} (f) \right)  + \frac 1 {2\tau_k}  \norm{f - f_k}_{L^2}^2 \label{eq:FWBWSplit-a} \\
 f_{k+1} &\in  \begin{cases} 
                  \argmin_{f \in L^2(\Omega) } \cR(f)  + \frac 1 {2\sigma_k}  \norm{f - f_{k+\frac 1 2}}_{L^2}^2 &\text{(backward-backward)} \\
                  f_{k+\frac 1 2} - \sigma_k \partial \cR(f_{k+\frac 1 2}) &\text{(backward-forward)}
                 \end{cases}  \label{eq:FWBWSplit-b}  %
\end{align}
\end{subequations}
with stepsize-parameters $\tau_k, \sigma_k$. Algorithms of this kind have been proposed and analyzed in \cite{Bertsekas2011_IncrProxMethods,Bertsekas2011_IncrProxMethods2,AndersenHansen2014_genARTProx}.
Importantly, the proximal iteration for the data-fidelity, \eqref{eq:FWBWSplit-a}, can be computed efficiently via the GenSART-scheme from \cref{thm:genL2SART}. Similarly, GenSART-formulas may be used to compute proximal steps in block-primal-dual methods as considered in \cite{ChambolleEhrhardtEtAl2017StochasticPDHG}. 

%
%
\optspace
\section{Numerical examples} \label{S:NumExamples}

All of the GenSART-schemes from \secref{S:Applications} have been successfully implemented as numerical algorithms.
In the following, exemplary results are presented.

\subsection{Implementation} \label{SS:Implementation}


In previous studies, Kaczmarz-type  reconstruction methods have usually been derived for a \emph{discretized} tomographic model. On the contrary, the theory in this work relies on properties of the parallel- or cone-beam projectors $\Proj \in \{\ParProj, \Cone\}$ that are valid only in \emph{continuous space}. In particular, while the generalized SART \cref{thm:genSART} is equally valid in finite dimensions, the tentative simplicity of $\Proj \Proj^\ast$ does typically \emph{not} carry over to discretizations of these operators. Moreover, discrete analogues of the penalty functionals considered in \secref{S:AdmPenalties} will in general no longer satisfy \cref{A1}.

For this reason, we pursue a \emph{post-discretization} strategy: we propose Kaczmarz-iterations within a \emph{continuous} tomographic model, then devise GenSART-schemes for their computation using the results from the preceding sections and \emph{finally} discretize these schemes to obtain numerically implementable SART-like iterations. If $\BP \in \mR^{m\times n}, \Bf \in \mR^{n}, S, \bu\in \mR^{m}$ are suitable discretizations of $\Proj, f, \cS, \uproj$, then this would look as follows for $L^2$-penalized Kaczmarz:
 \begin{align}
  f_{\new} &\in \argmin _{f \in L^2(\Omega)}  \cS( \Proj( f ) ) + \alpha \norm{ f - f_{\REF} }_{L^2}^2 \label{eq:Post-discretization-1} \\
  \stackrel{\textup{(GenSART)}}\rightsquigarrow f_{\new} &\in f_{\REF} + \adj \Proj \Bparens{ \tuproj^{-1/2} \cdot  \argmin_{ p \in \DetDom  } \cS\bparens{ \Proj( f_{\REF} ) + \tuproj^{1/2} \cdot p } + \alpha \norm{ p }_{L^2}^2  } \label{eq:Post-discretization-2} \\
  \stackrel{\textup{(discretize)}}\rightsquigarrow \Bf_{\new} &\in \Bf_{\REF} + \adj \BP \Bparens{ \btu^{-1/2} \odot  \argmin_{ \Bp \in \mR^m  } S\bparens{  \BP( \Bf_{\REF} ) + \btu^{1/2} \odot \Bp } + \alpha \norm{ \Bp }_{2}^2  } \label{eq:Post-discretization-3}
 \end{align}
 A drawback of the approach is that the discrete GenSART-update \eqref{eq:Post-discretization-3} in general does not \emph{exactly} solve a discretized version of \eqref{eq:Post-discretization-1}, \IE $\Bf_{\new} \notin \argmin_{ \Bf \in \mR^n }  S( \BP( \Bf ) ) + \alpha \norm{ \Bf - \Bf_{\REF} }_{2}^2$, but only up to \emph{discretization errors}. While this inexactness might lead to numerical instabilities in principle, no such effects are observed for the  examples in \secref{SS:RobustExample} and \secref{SS:PCTExample}.


As usual, we assume discrete objects $\Bf \in \mR^{\Nvx}$ and data $\Bg_j \in \mR^{\Mpx}$ to provide samples of the continuous quantities $f \in L^2(\Omega)$ and $g_j$ on equidistant Cartesian grids. 
Pointwise operations on functions are then represented by element-wise operations of vectors.
Integrals in continuous space can be approximated by summation over the entries of the corresponding vectors in the discretized model and derivatives can be implemented via finite differences. For example, $L^q$-norms are then identified with $q$-norms in $\mR^n$, i.e.\
\begin{align}
\int_{\Omega} |f(\bx)|^q \, \D \bx = \norm{f}_{L^q}^q \sim \norm{\Bf}_q^q = \sum_{i = 1}^n |\Bf_i|^q \MTEXT{if} \Bf = (\Bf_i)_{i =1}^n \text{ discretizes } f.
\end{align}
Discrete and continuous quantities can be related via sampling operators,
\begin{align}
 S_{\textup O}: L^2(\Omega) \to \mR^\Nvx; \; f \mapsto \big( {\textstyle \int_{\mV_i} f \, \D x } \big)_{i=1}^\Nvx
\end{align}
for voxels $\mathbb V_i \subset \Omega$ disjointly covering the object-domain $\Omega$, and an analogous map $S_{\textup D}: L^2(\mD) \to \mR^\Mpx$ in the projection-domain.
Projectors $\BP_{j} \in \{\ParProj_{\btheta}, \Cone_{\src}\}$  may then be naturally discretized via $\BP_{j} := S_{\textup D} \Proj_j \adj{S_{\textup O}}$. See also \cite{XuMueller2006_TomoInterpMethods} for alternative discretizations.

\paragraph{Unit-projections and precomputations:} Recall that all of the generalized SART formulas in \secref{S:AdmPenalties} and \secref{S:Applications} involve (weighted) unit-projections $\uprojj{j}$ or $\tuprojj{j}$. Clearly, discrete approximations $\bu_{j}, \tilde \bu_{j}$ of these are needed for numerical computations. For general object-domains $\Omega$, these may be precomputed via $\bu_{j} = \BP_j(\boldsymbol 1_{\mR^{\Nvx}})$ and $\tilde \bu_{j} = \BP_j(\bw_j)$, where $\bw_j$ is a suitable discretization of the ray-density. In this case, GenSART-methods thus require one additional evaluation of the full projector $(\BP_1, \ldots, \BP_{\Nproj})$ prior to the actual iterations. However, for geometrically simple domains $\Omega$ such as boxes, cylinders and balls,  $\uprojj{j} = \Proj_j( \boldsymbol 1 _\Omega)$, $\tuprojj{j} = \Proj_j( w_j)$ can also be computed analytically. Thereby, costly precomputations may be avoided.

Furthermore, discretized GenSART-schemes typically involve element-wise \emph{divisions} by powers of $\bu_{j}$ or $ \tilde \bu_{j}$. In accordance with \cref{conv:ProjOps} for the underlying continuous-space model, such operations should only be performed for \emph{non-vanishing} entries of $\bu_{j}, \tilde \bu_{j}$. For all other entries, the result of expressions of the form $\Bp \oslash \bu_{j}^\mu$ may simply be taken to be zero.

\subsection{Robust tomography test case} \label{SS:RobustExample}

\begin{figure}[hbt!]
 \centering
 \includegraphics[width=.9\textwidth]{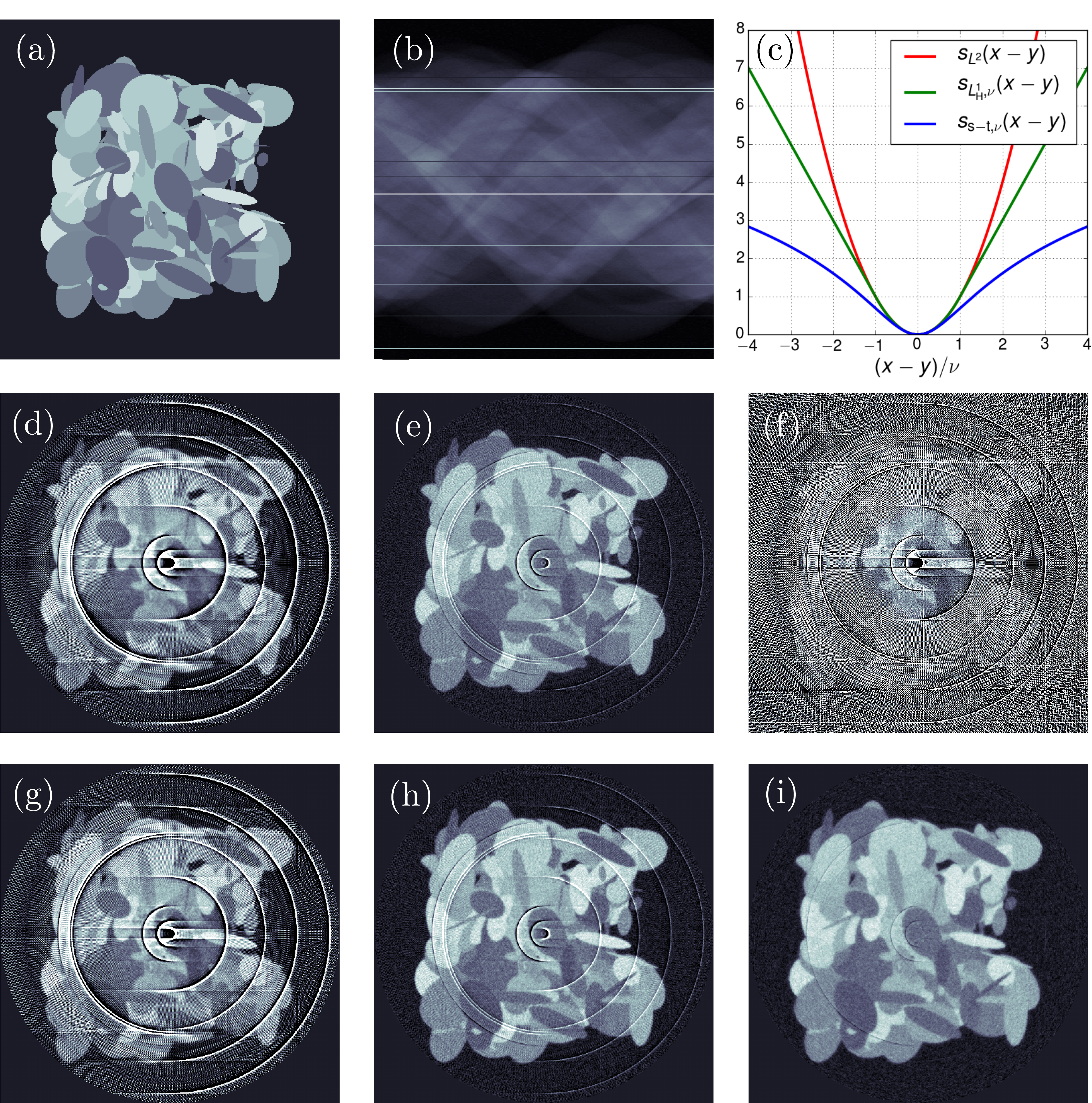}
 \caption{Numerical robust tomography test case: (a)~$512\times 512$-sized object phantom. (b)~Simulated parallel-beam data. 
 (c)~Plot of the different data-fidelity functions. (f) FBP-reconstruction.  (d),(e)~Tikhonov-reconstruction (\cref{alg:RobustTikh}) with $L^2$- and $L^1$-Huber-data-fidelity.  (g)--(i) SART-reconstruction (\cref{alg:RobustSART}) with $L^2$- and $L^1$-Huber- and Student's-t data-fidelity. The linear color scale in all object-plots (subfigures (a) and (d)--(i)) is identical. See text for details.}
 \label{fig:RobustSART}
\end{figure}

As a first numerical example, we consider the application of robust reconstrution from tomographic projection data,   as introduced in \ref{SSS:RobustDataFid}. To this end, we compare $L^2$-regularized Tikhonov-regularization and Kaczmarz-iterations
\begin{subequations}
\begin{align}
 f^{\textup{Tik}} &= \argmin_{f \in L^2(\Omega) } \cS_{\tot}\left( g_\tot^\obs; \, \Proj_{\tot}(f)\right) + \alpha_{\textup{Tik}} \norm{f}_{L^2}^2 \label{eq:AppRobustTikh} \\
  f_{k+1} &= \argmin_{f \in L^2(\Omega) } \cS \left( g_{j_k}^\obs ; \Proj_{j_k} (f)\right) + \alpha \norm{f-f_k}_{L^2}^2 \nnl
  &= f_k + \adj \Proj_{j_k} \bigg( \tuprojj{j_k}^{-1/2} \cdot \argmin_{p \in \DetDomj{j_k}} \cS \left( g_{j_k}^\obs ; \Proj_{j_k} (f_k) + \tuprojj{j_k}^{1/2} \cdot p \right) + \alpha \norm{p}_{L^2}^2 \bigg). \label{eq:AppRobustKaczarz}
\end{align}
\end{subequations}
with data-fidelities given by $ \cS_{\tot}\left( g^\obs_{\tot}; \, \Proj_{\tot}(f)\right) = \sum_{j = 1}^{\Nproj} \cS  ( g_j^\obs; \, \Proj_{j}(f) )$ where
\begin{align}
 \cS \left( g_j^\obs; \, p \right) = \int_{\DPj {j}} s(p(\bx) - g_j^\obs (\bx)) \, \D \bx, \qquad s  \in  \{s_{L^2}, s_{L^1_{\textup H}, \nu }, s_{\textup{s-t}, \nu}\}.
\end{align}
The $L^1$-Huber- and Student's t-functions $s_{L^1_{\textup H}, \nu }$ and $s_{\textup{s-t}, \nu}$ are defined in \eqref{eq:DataFidRobust} and $s_{L^2}(y) := |y|^2$  simply corresponds to an $L^2$-data-fidelity. $s_{L^2}$, $s_{L^1_{\textup H}, \nu }$ and $s_{\textup{s-t}, \nu}$ are plotted in \cref{fig:RobustSART}(c).

As numerical reconstruction methods, we consider for once a discretized version of the Tikhonov regularization in \eqref{eq:AppRobustTikh}, stated as \cref{alg:RobustTikh}. On the other hand, we design a generalized SART-analogue, \cref{alg:RobustSART}, by discretizing the update-formula \eqref{eq:AppRobustKaczarz}. Note that the discrete analogue of the optimization problem in \eqref{eq:AppRobustKaczarz} factorizes into a family of scalar problems just like in the continuous setting, see \secref{SS:NoiseAdaptedSART}. This enables a highly efficient implementation of this step regardless of the choice $ s  \in  \{s_{L^2}, s_{L^1_{\textup H}, \nu }, s_{\textup{s-t}, \nu}\}$.
 
\begin{algorithm}
\caption{Robust Tikhonov reconstruction}
\label{alg:RobustTikh}
\begin{algorithmic}
\vspace{.1em}
\REQUIRE{Data $ \Bg^\obs_\tot \in \mR^{\Mpx  \Nproj}$, projector $\BP_{\tot} \in  \mR^{( \Mpx  \Nproj)  \times  \Nvx}$, regularization parameter $\alpha_{\textup{Tik}} > 0$, data-fidelity $s \in \{s_{L^2}, s_{L^1_{\textup H}, \nu} , s_{\textup{s-t}, \nu}\}$, $\cS_\tot(\Bg ;\,\Bp ):= \sum_{i=1}^{\Mpx  \Nproj} s(p_i - g_i)$}
\vspace{.25em}
\ENSURE{$\Bf^{\textup{Tik}} \in  \argmin_{\Bf \in \mR^{\Nvx }} \cS_\tot\left( \Bg^\obs_\tot;\, \BP_{\tot} \Bf\right) + \alpha_{\textup{Tik}} \norm{\Bf}_{2}^2$}
\end{algorithmic}
\end{algorithm}

\begin{algorithm}
\caption{Robust SART reconstruction}
\label{alg:RobustSART}
\begin{algorithmic}
\vspace{.1em}
\REQUIRE{Data $ \Bg^\obs_j \in \mR^{\Mpx}$, projectors $\BP_{j} \in  \mR^{\Mpx \times \Nvx}$, regularization parameter $\alpha > 0$, initial guess $\Bf_0 \in \mR^\Nvx$, data-fidelity $s \in \{s_{L^2}, s_{L^1_{\textup H}, \nu} , s_{\textup{s-t}, \nu}\}$, $\cS (\Bg;\,\Bp):= \sum_{i=1}^\Mpx s(p_i - g_i)$, weighted unit-projections $\tilde \bu_{j}$.}
\vspace{.25em}
\FOR{$k = 0,\ldots, k_{\textup{stop}}-1$}
\STATE{$\Bp_k = \BP_{j_k}(\Bf_k)$}
\STATE{$\Delta \Bp_k  \in \argmin_{\Bp \in \mR^{\Mpx}} \cS \Bparens{ \Bg_{j_k}^\obs; \, \Bp_k  + \tilde \bu_{j_k}^{1/2} \odot \Bp } +  \alpha  \norm{ \Bp }_{2}^2$}
\STATE{$\Bf_{k+1} = \Bf_k + \adj{\BP_{j_k}} \Bparens{ \Delta \Bp_k \oslash \tilde \bu_{j_k}^{1/2} }$}
\ENDFOR
\vspace{.25em}
\ENSURE{final object-iterate $\Bf_{k_{\STOP}}$}
\end{algorithmic}
\end{algorithm}

We test the different tomographic reconstruction algorithms for a 2D-phantom $\Bf^{\textup{exact}}$ of image-size $512\times 512$ pixels composed of random ellipses, shown in \cref{fig:RobustSART}(a).
Note that such a 2D-setting can be treated as a limiting case of the 3D-geometries studied in this work. Parallel-beam data is simulated from the phantom under tomographic incident angles $\theta = 0^\circ, 1^\circ, \ldots, 179^\circ$. To avoid ``inverse crime'', different discretizations of the tomographic projectors are used in the data-simulation and the reconstruction. The simulated exact data is corrupted by Gaussian white noise of relative magnitude $2\,\%$ in 2-norm and additionally  by systematic outliers: following \cite{Bleichrodt2015_STTomo}, we simulate $2\,\%$ randomly chosen dead detector pixels that measure a constant incorrect value in all projections. Such measurement errors manifest as  horizontal stripes in the sinogram-plot of the tomographic data shown in \cref{fig:RobustSART}(b). 

The  data is reconstructed with \cref{alg:RobustTikh,alg:RobustSART} using the data-fidelities $ s_{L^2}, s_{L^1_{\textup H}, \nu }$, $s_{\textup{s-t}, \nu} $ with regularization parameters $\alpha_{\textup{Tik}} = 300$ and $\alpha = 2\alpha_{\textup{Tik}}$, the latter choice being motivated by \secref{SS:KaczmarzVSTikhonov}. The robustness-parameter $\nu$ (compare \secref{SSS:RobustDataFid}) is taken as $20\,\%$ of the standard-deviation of the data. \Cref{alg:RobustSART} is run for one symmetric Kaczmarz-cycle, i.e.\ $k_{\STOP} = 2 \Nproj$ and $j_k = j_{k_{\STOP}-1-k}$ for all $k$. The processing order $j_0, \ldots, j_{\Nproj-1}$ of the tomographic projections is chosen according to a multi-level-scheme from \cite{GuanGordon1994_MLS}.  

The Tikhonov results for the data-fidelities $s \in \{ s_{L^2}, s_{L^1_{\textup H}, \nu} \}$ (computed by a conjugate-gradient method and a linearly convergent primal-dual algorithm from  \cite{ChambollePock2011}, respectively) are plotted in \cref{fig:RobustSART}(d) and (e).  It is not clear how to  reliably minimize the Tikhonov-functional for the non-convex Student's-t functional $s = s_{\textup{s-t}, \nu}$ so that this case is omitted here. The GenSART-results for all(!) data-fidelities are plotted in \cref{fig:RobustSART}(f)-(h). For comparison, \cref{fig:RobustSART}(c) also shows a reconstruction by filtered back-projection (FBP), computed using an implementation from the ASTRA-toolbox \cite{ASTRAToolbox2015,ASTRAToolbox2016} with default-parameters.

As expected, the results for $L^2$-data fidelities (\cref{fig:RobustSART}(d),(f)) show pronounced ring- and pattern-artifacts arising from the dead pixels -- not to speak of the FBP-solution in \cref{fig:RobustSART}(c). The artifacts are significantly reduced for the $L^1$-Huber-reconstructions (\cref{fig:RobustSART}(e),(g)).
Finally, the reconstruction with the non-convex and even more robust Student's-t data-fidelity in \cref{fig:RobustSART}(h) is almost completely free of artifacts.

Notably, for both the $L^2$- and the non-quadratic(!) $L^1$-Huber-term, the Tikhonov- and GenSART-reconstructions turn out to be qualitatively indistinguishable. This suggests that variational methods can indeed be emulated by Kaczmarz-iterations -- even in  settings where this is not predicted by the convergence theory from \secref{SS:KaczmarzVSTikhonov}.
At the same time, it should be emphasized that the computational costs are practically identical for all three GenSART-methods, essentially amounting to \emph{two} evaluations of the full projector $\BP_{\tot}$ and its adjoint $\adj \BP_{\tot}$, which is significantly less than for the Tikhonov-reconstructions: for the non-quadratic $L^1$-Huber-case the primal-dual algorithm requires 73 iterations, and thus 73(!) evaluations of $\BP_{\tot}, \adj \BP_{\tot}$, to converge to a prescribed accuracy of $1\,\%$ according to a criterion from \cite{ChambollePock2011}.

\subsection{Newton-Kaczmarz-GenSART for experimental XPCT-data} \label{SS:PCTExample}

For a second and somewhat more involved numerical test case, we implement the Newton-Kaczmarz-iterations for X-ray phase contrast tomography (XPCT) from \secref{SSS:PCTSART}. To this end, the obtained GenSART-formula \eqref{eq:PCTSART} is discretized, where the gradients $\nabla$ are replaced by finite-difference operators $\bnabla \in \mR^{M_{\grad}\times \Mpx}$, $\big(\bnabla ((p_i)_{i=1}^\Mpx)\big)_m = p_{k_m} - p_{i_m}$ for index pairs $(i_m, k_m)$ of neighboring pixels. Importantly, as implied by the treatment of the continuous-space setting in \secref{SS:GradPenalties}, the set $\cN_{\bnabla} = \{ ( i_m, k_m) \}_{m = 1}^{M_\nabla}$ must be restricted to such neighbor-index-pairs, for which  the unit-projection $\boldsymbol u _j = (u_{ji})_{i=1}^\Mpx$ does not vanish for either of the indices, i.e.\ $u_{ji_m}\neq 0$ and $ u_{jk_m} \neq 0$. In a nutshell, this means that the discrete gradient $\bnabla$ must  only be computed \emph{within} the support of the unit-projection in accordance with \cref{conv:ProjOps}. Upon discretizing \eqref{eq:PCTSART}, the optimization problem in projection-space assumes the form
\begin{align}
 \Delta \Bp_k = \argmin_{\Bp \in \mR^\Mpx} \Norm{\BT_k (\bu_{j_k}^{1/2} \odot \Bp ) - \boldsymbol r_k }_{2}^2  + \alpha(1-\gamma) \norm{\Bp}_2^2 + \alpha \gamma \Norm{\bU_{\bnabla, j_k}^{1/2} \bnabla \bparens{ \Bp \oslash \bu_{j_k}^{1/2} }  }_2^2 \label{eq:PCTDiscrProjOptProblem} 
\end{align}
where $\boldsymbol r_k$ is the residual and $\BT_k = F'[\BP_{j_k} ( \Bf_k )] \in \mR^{\Mdata \times \Mpx}$ is the Fr\'echet-derivative of the discretized image-formation map $  F:  \mR^{\Mpx}  \to \mR^{\Mdata}$. The $\bU_{\bnabla, j} \in \mR^{M_{\grad} \times M_{\grad}}$ are diagonal, positive-semidefinite matrices that implement a discrete analogue $\bnabla( \Bp) \mapsto \bU_{\bnabla, j}  \bnabla( \Bp)$ of the multiplication $\nabla(p) \mapsto u_j \cdot \nabla(p)$ for gradients of continuous projections $p$.

With $\bU_{j} := \textup{diag}(\bu_{j})  \in \mR^{\Mpx \times \Mpx}$ 
and $\boldsymbol I \in \mR^{\Mpx \times \Mpx}$ denoting the identity, 
\eqref{eq:PCTDiscrProjOptProblem} then amounts to solving the associated normal equation:
\begin{align}
 \Delta \Bp_k = \Parens{  \bU_{j_k}^{1/2} \adj{\BT_k}  \BT_k \bU_{j_k}^{1/2}  + \alpha(1-\gamma) \boldsymbol I + \alpha \gamma \bU_{j_k}^{-1/2} \adj \bnabla \bU_{\bnabla, j_k} \bnabla \bU_{j_k}^{-1/2}  }^{-1} \adj{\BT_k} \bparens{ \bu_{j}^{1/2} \odot  \boldsymbol r_k }
\end{align}
The map $F$ as well as its derivative $ F'[\Bp]$ can be implemented using only fast Fourier transforms and pointwise operations, leading to a computational complexity of $\cO(\Mdata \log(\Mdata))$. Accordingly, the above symmetric positive-definite problem can be solved efficiently by a matrix-free conjugate-gradient (CG) method. The obtained Newton-Kaczmarz-GenSART method is summarized in \cref{alg:PCTSART}, which is notably not limited to XPCT but may be adapted for a wide range of other image-formation operators $F$.

\begin{algorithm}
\caption{Newton-Kaczmarz-GenSART (for X-ray phase contrast tomography)}
\label{alg:PCTSART}
\begin{algorithmic}
\REQUIRE{Data $ \Bg^\obs_j \in \mR^{\Mdata}$, parallel-beam projectors $\BP_{j} \in  \mR^{\Nvx \times \Mpx}$, initial guess $\Bf_0 \in \mR^\Nvx$, regularization $\alpha \! > \!0$, $\gamma\!\in \![0;1]$, unit-projections $ \bu_{j}$, $\bU_{j} := \textup{diag}(\bu_{j})$ and  $\bU_{\bnabla, j}  \in \mR^{m_\grad \times m_\grad}$, 
image-formation operator $F: \mR^{\Mpx} \to \mR^{\Mdata}$, Fr\'echet-derivative $F'$.} \vspace{.25em}
\FOR{$k = 0,\ldots, k_{\STOP}-1$}
\STATE{$\Bp_k = \BP_{j_k}  \Bf_k$, $\boldsymbol r_{k} = \Bg^{\textup{obs}}_{j_k} - F \left(  \Bp_k \right)$, $\BT_k = F'[\Bp_k]$}
\STATE{$ \Delta \Bp_k \stackrel{\textup{CG}}=  \Parens{  \bU_{j_k}^{1/2} \adj{\BT_k}  \BT_k \bU_{j_k}^{1/2}  + \alpha(1-\gamma) \boldsymbol I + \alpha \gamma \bU_{j_k}^{-1/2} \adj \bnabla \bU_{\bnabla, j_k} \bnabla \bU_{j_k}^{-1/2}  }^{-1} \adj{\BT_k} \bparens{ \bu_{j}^{1/2} \odot  \boldsymbol r_k }$}
\STATE{$ \Bf_{k+1} =   \Bf_k +  \adj{\BP_{j}} \big( \bu_{j}^{-1/2}\odot  \Delta  \Bp_k \big)$}
\STATE{$\big[ \Bf_{k+1} =  \max \{ 0, \Bf_{k+1} \} \big] \qquad  \qquad \qquad  \text{\% optional non-negativity constraint}$}
\ENDFOR \vspace{.25em}
\ENSURE{final object-iterate $\Bf_{k_{\STOP}}$}
\end{algorithmic}
\end{algorithm}

\begin{figure}[hbt!]
 \centering
 \includegraphics[width=\textwidth]{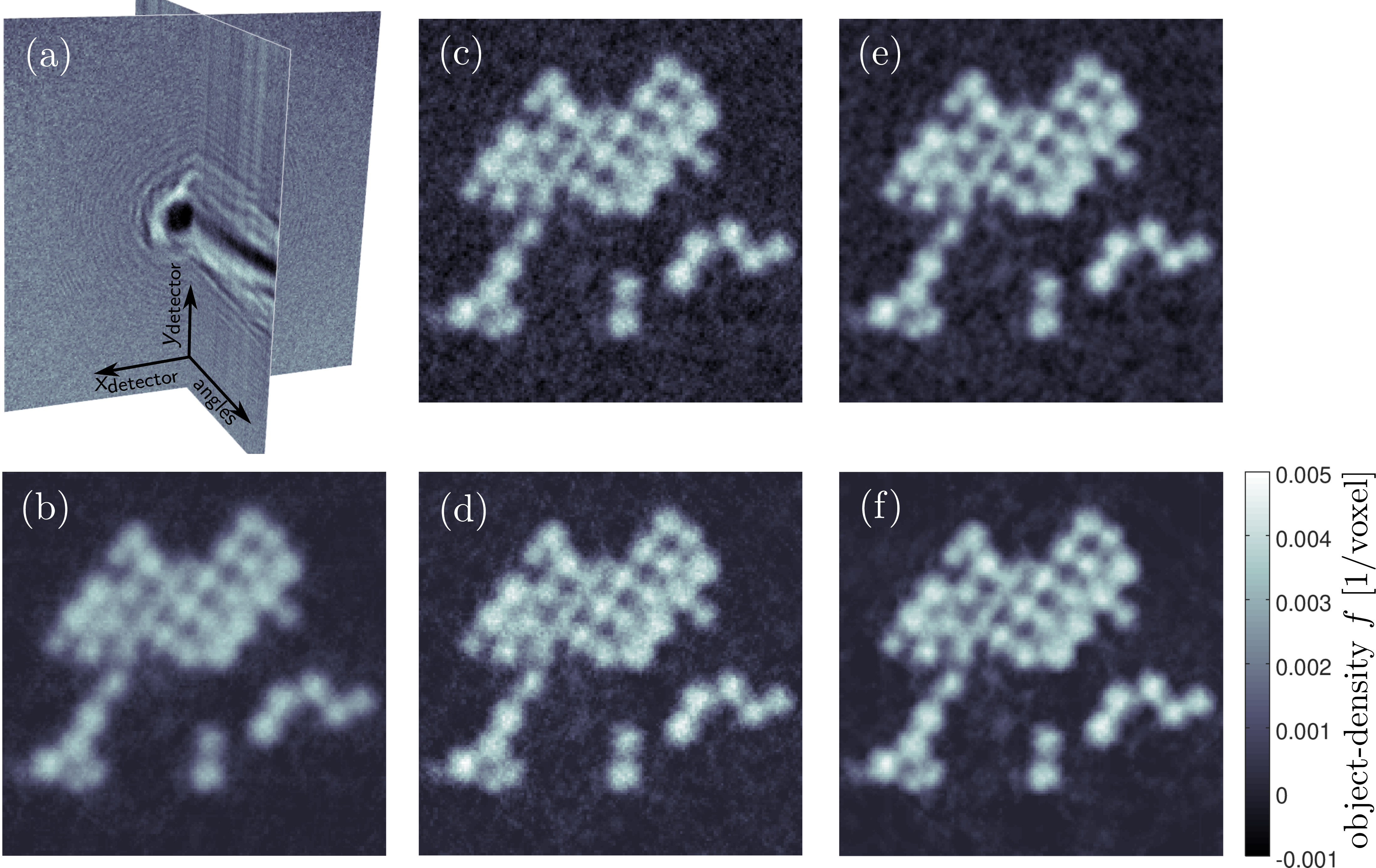}
 \caption{X-ray phase contrast tomography test case: (a)~Orthoslice-plot of the 3D tomographic data $\Bg^{\obs} = (\Bg^{\obs}_1, \ldots, \Bg^{\obs}_{\Nproj})$, composed of $1024\times1024$-sized diffraction patterns for 249 incident angles. (b)~Reconstruction from \cite{MaretzkeEtAl2016OptExpr}: central slice of the $256^3$-voxel volume perpendicular to the tomographic axis (zoomed to the object-containing region). (c)--(f)~Same slice of the reconstructed object after one cycle of \cref{alg:PCTSART}: (c)~$\gamma =0$, no non-negativity constraint. (d)~$\gamma=0$, with non-negativity. (e)~$\gamma=0.8$, no non-negativity. (f)~$\gamma=0.8$, with non-negativity. The color scale in (b)--(f) is identical.}
 \label{fig:PCTSART}
\end{figure}

We apply \cref{alg:PCTSART} to experimental XPCT-data $\Bg^\obs_{\tot} \in \mR^{\Nproj \Mdata}$ from \cite{MaretzkeEtAl2016OptExpr}, composed of $\Nproj = 249$ near-field diffraction patterns of size $\Mdata = 1024^2$ pixels each, acquired in an (effective) parallel-beam setting under tomographic incident angles between $0^\circ$ and $173^\circ$. The 3D tomographic data, measured at the synchrotron light-source  PETRAIII (see \cite{Salditt2015GINIX} for experimental details), is visualized by orthoslices in \cref{fig:PCTSART}(a).
 
It is known that the  object $\Bf$ -- a colloidal nano-crystal -- is contained in a centered cube of size $256^3$ voxels, which we take as the reconstruction volume $\Omega$ (\emph{support constraint}). Moreover, the true values of $\Bf$ are known to be non-negative for physical reasons. 
We compute reconstructions using \cref{alg:PCTSART} with $\Bf_0 = 0$ and regularization parameter $\alpha = 500$ in four different setups: with $\gamma = 0$ (i.e.\ pure $L^2$-penalty) and $\gamma = 0.8$ (primarily gradient-penalization), with and without a non-negativity constraint, respectively.
For each of these settings, \emph{one} (non-symmetric) Kaczmarz-cycle is run, again following the MLS-scheme from \cite{GuanGordon1994_MLS} for the processing order of the data. \Cref{fig:PCTSART}(c)-(f) plots the central slices of the reconstructed $256^3$-voxel volumes perpendicular to the tomographic axis. For comparison,  \cref{fig:PCTSART}(b) shows the same view for the reconstruction from \cite{MaretzkeEtAl2016OptExpr}, obtained by a similar $L^2$-regularized Newton-Kaczmarz method, but without using a GenSART-scheme.

By nature of the imaged sample, the slices should show circular cross-sections of the colloids of constant density within a zero background, i.e.\ a binary image.
It can be seen that all methods perform widely similarly in terms of contrast and resolution where the four reconstructions with \cref{alg:PCTSART} (\cref{fig:PCTSART}(c)-(f)) seem to yield somewhat less blurry colloids with a higher contrast to the background than the reference-result in \cref{fig:PCTSART}(b). Moreover, comparing \cref{fig:PCTSART}(c),(d) to (e),(f) respectively, reveals the anticipated noise-suppressing effect of the gradient-penalty compared to pure $L^2$-regularization, at essentially the same sharpness and contrast. Indeed, the reconstruction in \cref{fig:PCTSART}(f) using both a gradient-penalty and non-negativity constraint seems closest to the expected ideal binary object and is arguably  best suited for further image-segmentation and -analysis -- despite the slight approximation involved in the computation via GenSART-schemes (compare \secref{SS:GradPenalties}).

\optspace
\section{Conclusions} \label{S:Conclusions}

In this work, efficient solution formulas have been proposed for the computation of regularized Kaczmarz-iterations (also known as ``Tikhonov-Kaczmarz'' or ``incremental proximal iterations'') for tomographic reconstruction. By their structural analogy and similar computational efficiency to classical SART-iterations \cite{AndersenKak1984SART}, the derived schemes are termed generalized SART (GenSART). Notably, the approach strongly exploits mathematical structures specific to tomographic inverse problems, namely the well-known fact that the acquired data for a single tomographic view only contains information 
on the \emph{projection} of the three-dimensional object to a 2D-manifold. As a consequence, Kaczmarz-updates may often be computed in the \emph{lower-dimensional} projection-space, reducing the computational costs.

As demonstrated in the manuscript, this enables efficient Kaczmarz-methods for various non-standard tomographic settings such as robust reconstruction in the presence of large outliers (\secref{SSS:RobustDataFid}, \secref{SS:RobustExample}), photon-starvation- (\secref{SSS:PoiNoise}) and beam-hardening-resistant tomography (\secref{SSS:PolyCTSART}) or phase contrast tomography (\secref{SSS:PCTSART}, \secref{SS:PCTExample}).
The key feature of the GenSART-approach here lies in its flexibility in the choice of the data-fidelity and image-formation model (including even \emph{non-convex} choices) -- essentially without effect on the computational costs! Note that this constitutes a crucial difference to (bulk-) variational reconstruction methods, where the number of required iterations typically grows dramatically as soon as non-quadratic functionals are to be minimized. 
The proposed methods are furthermore applicable to both parallel-beam- and cone-beam-tomography, absolutely regardless of the sampling in the tomographic incident directions and naturally incorporate support constraints on the unknown object, thus retaining the geometrical flexibility of SART.

In addition to (weighted) $L^2$-penalties and Sobolev-$W^{1,2}$-regularizers that can be directly incorporated into GenSART-schemes (yet, the latter with the flaw of a slight approximation, see \secref{SS:GradPenalties}), the approach may be combined with interlacing proximal- or gradient-descent iterations to enable advanced regularizers, as outlined in \secref{SS:Extensions}. Likewise, GenSART-formulas might serve to efficiently evaluate proximal maps in Kaczmarz-type primal-dual methods of the kind proposed in \cite{ChambolleEhrhardtEtAl2017StochasticPDHG}. Such applications constitute promising prospects for future work.

Importantly, this work is \emph{not} concerned with convergence of Kaczmarz-methods beyond known results (e.g.\ \cite{KindermannLeitao2014_genARTConvergence}) but focusses on an efficient computation of the iterates. Yet, it clearly motivates a further analysis of Kaczmarz-convergence with non-quadratic data-fidelities.

Another important goal is \emph{parallelization}: while GenSART-iterations are efficient in terms of computational complexity, Kaczmarz-methods are more \emph{sequential} in nature than bulk-variational methods (though, by far not as sequential as ART). This complicates accelerating the approach by parallel computing, which for bulk methods can be readily achieved by using massively parallelized implementations of tomographic (back-)projections from the ASTRA-toolbox \cite{ASTRAToolbox2015,ASTRAToolbox2016} for example. If the issue of parallelization is resolved, the results of this work show that Kaczmarz-methods computed via GenSART-schemes may provide a highly efficient alternative to variational reconstruction methods -- at almost the same flexibility.

\appendix
\optspace
\section{Geometry of the Projectors} \label{Appendix:ProjPropProof}
  \begin{proof}[Proof of \Cref{thm:AdjProjClosedRange}:]
  We show that $\isoBProjAbstr: \DetDom \to L^2(\Omega); \; p \mapsto \wproj \cdot \backProjAbstr(\tuproj^{-1/2} \cdot p)$ is well-defined and isometric.
   Let $p \in \DetDom$ be arbitrary. For $\Proj = \ParProj$, we have
   \begin{align}
    \norm{ \isoBProjAbstr (p)  }_{L^2}^2 &= \int_{\Omega} \left|  \backProj(\tilde u_{\Proj}^{-1/2} \cdot  p)(\bx) \right|^2 \, \D \bx \stackrel{\eqref{eq:DefBackProj}}= \int_{\DR} \left|   (\tilde u_{\Proj}^{-1/2} \cdot  p)(\bx_\perp) \right|^2 \cdot \!\!\!\! \underbrace{ \Bparens{ {\textstyle \int}_{\OmegaX}  1  \, \D z }}_{= \Proj(\boldsymbol 1_\Omega) (\bx_\perp) = u_{\Proj} (\bx_\perp) } \!\!\!\! \D \bx_\perp \nnl
    &= \int_{\Omega_{\Proj}} \left(  u_{\Proj}(\bx_\perp)/ \tilde u_{\Proj}(\bx_\perp) \right) \cdot   |p  (\bx_\perp)|^2 \, \D \bx_\perp \stackrel{\tilde u_{\Proj} = u_{\Proj}}= \norm{p}_{\DetDomR}^2 .
   \end{align}
   By introducing polar coordinates, we obtain the same for the cone-beam case $\Proj = \Cone$:
    \begin{align}
    \norm{ \isoBProjAbstr(p)  }_{L^2}^2 &= \int_{\Omega} \left| w_{\Cone} \cdot \backCone(\tilde u_{\Cone}^{-1/2} \cdot p)(\bx) \right|^2 \, \D \bx = \int_{\DC} |(\tilde u_\Cone^{-1} \cdot p ) (\bvarphi)|^2 \int_{\OmegaPhi}  |w_\Cone(t \bvarphi)|^2  t^2 \, \D t \, \D \bvarphi \nnl
    &\stackrel{w_\Cone(t \bvarphi) = t^{-2}}= \int_{\Omega_{\Cone}} \tilde u_\Cone (\bvarphi)^{-1}  \cdot   |p  (\bvarphi)|^2 \cdot \underbrace{ \Bparens{ {\textstyle \int}_{ \OmegaPhi } t^{-2}  \, \D t } }_{\Cone(w_{\Cone})(\bvarphi) = \tilde u_\Cone  (\bvarphi) } \, \D \bvarphi  = \norm{p}_{\DetDom}^2.  
   \end{align}
   
   Accordingly, $\isoBProjAbstr: \DetDom \to L^2(\Omega)$ is well-defined and isometric for $\Proj \in  \{\ParProj, \Cone \}$. Since \smash{$\isoBProjAbstr(p) = \adj\Proj(\tuproj^{-1/2} \cdot p)$}, its adjoint is given by \smash{$\adj \isoBProjAbstr(f) = \tuproj^{-1/2} \cdot \Proj(f) = \isoProj(f)$}. Hence, $\isoProj = \adj \isoBProjAbstr:L^2(\Omega) \to \DetDom$ is bounded with norm $\norm{\isoProj} = \norm{\isoBProjAbstr} = 1$. By the isometry-property, $\adj \isoProj = \isoBProjAbstr$ has closed range. According to the closed range theorem, the same thus holds true for $\isoProj$. As $\adj\isoProj$ is isometric, $\adj\isoProj: \DetDom \to \Range(\adj\isoProj)$ defines a unitary operator, which implies the characterizations of $\adj \isoProj \isoProj$ and $\isoProj \adj \isoProj = \identity_{\DetDom}$.
%
 \end{proof}

\optspace
\section{Projectors and Gradients} \label{Appendix:ProjGradCompatibility}
\begin{proof}[Proof of \Cref{lem:ProjGradient}:]
  Let $f \in \CcInfO$ be smooth and compactly supported. Then $f$ can be identified with a function in $\CcInf{\mR^3}$ by simply extending it with $0$ outside $\Omega$ (notably, this would \emph{not} be true if we only assumed $f \in \CInfO$).
  
  We decompose the gradient into mutually orthogonal components along- and perpendicular to the rays: in the parallel-beam case $\Proj = \ParProj$, the ray-direction is constant and points along the $z$-axis (without loss of generality). Accordingly, $\nablaProj$ is simply the gradient with respect to the remaining Cartesian coordinates $\bx_\perp$, i.e.\ $\nablaProj = \nabla_\perp$. Since $\Omega$ is bounded, there exist constants $a_{\Proj} < b_{\Proj}$ such that $\Omega$ is contained in the stripe $\mR^2 \times (a_{\Proj};\, b_{\Proj}) \supset \Omega$. Accordingly, we obtain by Leibniz' integral rule
  \begin{align}
   \Proj(\wproj^{-1/2} \cdot  \nablaProj f) (\bx_\perp)  &= \Proj(  \nablaProj f) (\bx_\perp) = \int_{a_{\Proj}}^{b_{\Proj}} \nabla_\perp f(\bx_\perp, z) \, \D z \nnl
   &= \nabla_{\mR^2}  \int_{a_{\Proj}}^{b_{\Proj}}  f(\bx_\perp, z) \, \D z = \nabla_{\mR^2}  \Proj(f)(\bx_\perp)  = \nabla_{\mD}  \Proj(f)(\bx_\perp)  \label{eq:lem-ProjGradient-pf1} 
  \end{align}
  for all $\bx_\perp \in \mR^2$. Here, it was used that the ray-density is constant in the considered parallel-beam case, $\wproj = 1$ for $\Proj = \ParProj$. For the cone-beam case $\Proj = \Cone$, the gradient can be expressed in polar coordinates: if $f^{(\textup p)}$ is defined by $f^{(\textup p)}( \bvarphi , t ) := f( t \bvarphi )$, then 
  \begin{align} 
   \nabla f (t \bvarphi ) &= t^{-1} \nabla\!_{\bvarphi} f^{(\textup p)} (\bvarphi, t) + \bvarphi \partial_t f^{(\textup p)} (\bvarphi, t) \nnl
   &= w_{\Cone}(t \bvarphi)^{1/2}  \nabla\!_{\bvarphi} f^{(\textup p)} (\bvarphi, t) + \bvarphi \partial_t f^{(\textup p)} (\bvarphi, t) \MTEXT{for all} t > 0, \bvarphi \in \mS^2.
  \end{align}
  where the ray-density $w_{\Cone}(\bx) =  |\bx|^{-2} $ has been inserted.  The gradient-component $\nablaProj f$ perpendicular to the ray-direction $\bvarphi$ is given by the first term on the right-hand side. Since $\Omega$ is bounded and $0 \notin \closure \Omega$, the domain is contained in some annulus, $\Omega \subset \{ t\bvarphi: \bvarphi \in \mS^2$, $a_{\Cone} < t < b_{\Cone} \}$ for $0 < a_{\Cone} < b_{\Cone}$. By another application of Leibniz' rule, this implies
  \begin{align}
   \Proj(\wproj^{-1/2} \cdot \nablaProj f) (\bvarphi) &= \int_{a_{\Cone}}^{b_\Cone} \left( \wproj^{-1/2} \cdot \nablaProj f \right) (t \bvarphi) \, \D t = \int_{a_{\Cone}}^{b_\Cone} \nabla\!_{\bvarphi} f^{(\textup p)} (\bvarphi, t)   \, \D t \nnl 
   &= \nabla_{\mS^2} \int_{a_{\Cone}}^{b_\Cone} f^{(\textup p)} (\bvarphi, t)  \, \D t = \nabla_{\mS^2} \Cone(f)(\bvarphi) = \nabla_{\mD} \Proj(f)(\bvarphi) \label{eq:lem-ProjGradient-pf2}  
  \end{align}
  {for all} $\bvarphi \in \mS^2$.
  Equations \eqref{eq:lem-ProjGradient-pf1} and \eqref{eq:lem-ProjGradient-pf2} show that \eqref{eq:lem-ProjGradient-1} holds for all $f \in \CcInfO$ and $\Proj \in \{\ParProj, \Cone\}$. Since \smash{$\isoProj: L^2(\Omega) \to \DetDom;\; p \mapsto \tuproj^{-1/2} \cdot \Proj(p)$} is bounded with norm 1 according to \Cref{thm:AdjProjClosedRange}, this furthermore implies that
  \begin{align}
    \bnorm{\tuproj^{-1/2}\cdot \nabla_{\mD} \Proj(f)}_{L^2}  &= \bnorm{\tuproj^{-1/2}\cdot \Proj(\wproj^{-1/2} \cdot \nablaProj f)}_{L^2}\leq  \bnorm{ \wproj^{-1/2} \cdot \nablaProj f }_{L^2}  \nnl 
    &\leq c^{-1/2} \bnorm{ \nablaProj f }_{L^2} \leq c^{-1/2} \bnorm{ \nabla f }_{L^2} \leq c^{-1/2} \bnorm{  f }_{W^{1,2}} \label{eq:lem-ProjGradient-pf3}  
  \end{align}
  for all $ f \in \CcInfO$ and $c := \min_{\bx \in \Omega} \wproj(\bx) > 0$.
  As $\CcInfO$ is dense in $W^{1,2}_0(\Omega)$ with respect to the $W^{1,2}$-norm, \eqref{eq:lem-ProjGradient-pf3} shows that \smash{$p \mapsto \tuproj^{- 1 / 2 } \cdot \nabla_{\mD} \Proj (p)$} defines a bounded linear map $W^{1,2}_0(\Omega) \to \DetDom$. Hence, \eqref{eq:lem-ProjGradient-1} and \eqref{eq:lem-ProjGradient-2} remain valid for all $f \in W^{1,2}_0(\Omega)$.

  Now let $ \backProj(p) \in W^{1,2}(\Omega)$.
  According to \eqref{eq:DefBackProj}, $\backProj(p)$ is constant within $\Omega$ along the ray-direction. Hence, the corresponding gradient-components vanish in an $L^2$-sense so that $\nabla \big( \backProj(p) \big) = \nablaProj \big( \backProj(p) \big)$. By inserting the expressions for $\nablaProj$ in the parallel- and cone-beam geometry, this yields for almost all $(\bx_\perp, z), t\bvarphi \in \Omega$ 
  \begin{align}
   \nabla \big( \backParProj(p) \big) (\bx_\perp, z) &= \nablaProj  \big( \backParProj(p) \big) (\bx_\perp, z) \stackrel{\eqref{eq:DefBackProj}}=   \nabla_{\mR^2} ( p ) (\bx_\perp) \\
   \nabla \big( \backCone(p) \big) (t\bvarphi) &= \nablaProj \big( \backCone(p) \big) (t\bvarphi)  \stackrel{\eqref{eq:DefBackProj}}=   t^{-1}\nabla_{\mS^2} ( p ) (\bvarphi).
  \end{align}
  By decomposing the norm-integral into integrals perpendicular and along the ray-direction, respectively, we may now compute the $L^2$-norm of $\nabla \left( \backProj(p) \right)$:
  \begin{align}
    \Norm{\nabla \left( \backProj(p) \right)}_{L^2}^2  &= 
      \begin{cases} 
        \int_{\DR} \int_{\OmegaX} \left| \nabla_{\mR^2}  \left(  p \right) \left( \bx_\perp \right)  \right|^2  \, \D z \,  \D \bx_\perp &\text{for } \Proj  = \ParProj \\
        \int_{\DC}  \int_{\OmegaPhi}    \left| \frac 1 t \nabla_{\mS^2} \left(  p \right) \left( \bvarphi \right)  \right|^2  t^2 \, \D t  \, \D \bvarphi &\text{for } \Proj  = \Cone
       \end{cases} \nnl
       &= 
      \begin{cases} 
        \int_{\DR} \left| \nabla_{\mR^2} \left(  p \right) \left( \bx_\perp \right)  \right|^2 \left( \int_{\OmegaX}  \D z \right) \, \D \bx_\perp &\text{for } \Proj  = \ParProj \\
        \int_{\DC}  \left|  \nabla_{\mS^2} \left(  p \right) \left( \bvarphi \right)  \right|^2 \left( \int_{\OmegaPhi} \D t \right) \, \D \bvarphi &\text{for } \Proj  = \Cone
       \end{cases} \nnl
    &= \bnorm{ \uproj ^{1/2} \cdot \nabla_{\mD} (p)  }_{L^2}^2, \label{eq:lem-ProjGradient-pf5}  
  \end{align}
  which proves \eqref{eq:lem-ProjGradient-3}. The equality \eqref{eq:lem-ProjGradient-pf5} furthermore shows the equivalences
      \begin{align}
   \backProj(p) \in W^{1,2}(\Omega) \;  \Leftrightarrow \;  \nabla \left( \backProj(p) \right) \in L^2(\Omega) \; \Leftrightarrow \; \uproj ^{1/2} \cdot \nabla_{\mD}  (p) \in \DetDom.   
    \end{align}

  By continuity of $\nabla: W^{1,2}(\Omega) \to L^2(\Omega)$ and  $\nabla_{\mD} \Proj$, it is sufficient to show \eqref{eq:lem-ProjGradient-4} for $f \in \CcInfO$. Using the expressions for $\nabla \big( \backProj(p) \big)$ derived above, we obtain
  \begin{align}
    &\bip{\nabla \left(  \backProj(p)\right) }{ \nabla f }_{L^2} =  \bip{\nablaProj \left(  \backProj(p)\right) }{ \nabla f }_{L^2} =  \bip{\nablaProj \left(  \backProj(p)\right) }{ \nablaProj f }_{L^2} \nnl
    &\qquad\;=  \begin{cases} 
        \int_{\DR} \nabla_{\mR^2}  \left(  p \right) \left( \bx_\perp \right) \cdot \int_{\OmegaX} { \nabla_\perp f ( \bx_\perp, z) }    \, \D z \,  \D \bx_\perp &\text{for } \Proj  = \ParProj \\
        \int_{\DC}  \nabla_{\mS^2} \left(  p \right) \left( \bvarphi \right) \cdot \int_{\OmegaPhi} \frac{1}{t^2} {\nabla_{\bvarphi} f^{(p)} (\bvarphi, t) } t ^2 \, \D t  \, \D \bvarphi &\text{for } \Proj  = \Cone
       \end{cases} \nnl
   &\;\;\,\stackrel{\eqref{eq:lem-ProjGradient-pf1},\,\eqref{eq:lem-ProjGradient-pf2}}= \bip{ \nabla_{\mD} \left( p \right)}{\Proj \bparens{ \wproj^{-1/2} \cdot \nablaProj (f) } }_{L^2} \nnl
   &\quad\;\;\stackrel{\eqref{eq:lem-ProjGradient-1}}=  \bip{  \uproj ^{1/2} \cdot \nabla_{\mD} \left( p \right)}{\uproj ^{-1/2} \cdot \nabla_{\mD} \Proj\left( f  \right)}_{L^2} \MTEXT{for all} f \in \CcInfO. 
  \end{align}
   \end{proof}

   \optspace
\section{Admissibility of $L^q$-Penalties} \label{Appendix:LqPenaltyProof}
\begin{proof}[Proof of \Cref{thm:admissibleLq}:]
 Let $p \in \DetDom$ and $f_0 \in \Kern(\isoParProj)$ be arbitrary. If $ \cR (f_{\REF} + \isoParProj^\ast(p) + f_0) = \norm{\isoParProj^\ast(p) + f_0}_{L^q}^q = \infty$, then \eqref{eq:A1} trivially holds true. Hence, we may assume that $f := \isoParProj^\ast(p) + f_0 \in L^q(\Omega)$. With the rays $\OmegaX$ defined as in \eqref{eq:DefRaysDetDoms}, we then have
 \begin{align}
    \cR(f_{\REF} + \isoParProj^\ast(p) + f_0) &= \int_{\Omega} |f|^q \, \D x = \int_{\mR^2} \left( \int_{\OmegaX} \left| f(\bx_\perp, z) \right|^q \, \D z \right) \, \D \bx_\perp \label{eq:thm:admissibleLq-1} 
 \end{align}
 where the inner integrals are finite for almost all $\bx_\perp \in \mR^2$.
 An application of Jensen's inequality to the convex function $\mR \to \mR; \, x \mapsto |x|^q$ in these shows that, for almost all $\bx_\perp\in\mR^2$,
 \begin{align}
  \int_{\OmegaX} \left| f(\bx_\perp, z) \right|^q \, \D z  \geq \underbrace{ \bigg(\int_\OmegaX 1 \, \D z  \bigg)^{1-q} }_{= \uparproj(\bx_\perp)^{1-q}}\bbabs{ \int_{\OmegaX}  f(\bx_\perp,z)  \, \D z }^q. \label{eq:LqProof-0}
 \end{align}
 Comparing to \eqref{eq:DefRadon}, we find that the integral within the modulus exactly corresponds to an evaluation of \smash{$\ParProj = \uparproj^{1/2} \cdot \isoParProj$}. Since $f = \adj \isoParProj(p) + f_0$ with $f_0 \in \Kern(\isoParProj)$, this implies
  \begin{align}
  \int_{\OmegaX} \left| f(\bx_\perp, z) \right|^q \, \D z  \geq \uparproj(\bx_\perp)^{1-q} \left| \ParProj(f) (\bx_\perp) \right|^q =   \uparproj(\bx_\perp)^{1-q} \big| \ParProj \isoParProj^\ast (p) (\bx_\perp) \big|^q \nnl
  = \bigg(\int_\OmegaX 1 \, \D z  \bigg)^{1-q}\bigg| \int_\OmegaX \underbrace{\isoParProj^\ast(p)(\bx_\perp, z)}_{\textup{constant in $z$}} \, \D z \bigg|^q =  \int_\OmegaX |\isoParProj^\ast(p)(\bx_\perp, z)|^q \, \D z  \label{eq:LqProof-1} 
 \end{align}
 {for almost all} $\bx_\perp\in\mR^2$.
 Substituting the estimate \eqref{eq:LqProof-1} into \eqref{eq:thm:admissibleLq-1} proves \eqref{eq:A1}:
 \begin{align}
  \cR(f_{\REF} + \isoParProj^\ast(p) + f_0) \geq \int_{\mR^2} \bbparens{ \int_\OmegaX |\isoParProj^\ast(p)(\bx_\perp, z)|^q \, \D z }  \, \D \bx_\perp = \cR(f_{\REF} + \isoParProj^\ast(p)). \label{eq:LqProof-2}
 \end{align}
 In particular, \eqref{eq:LqProof-2} shows that $\norm{\isoParProj^\ast(p)}_{L^q}^q = \cR(f_{\REF} + \isoParProj^\ast(p)) < \infty$, i.e.\ $\isoParProj^\ast(p)\in L^q(\Omega)$.
 Since $\ParProj \isoParProj^\ast(p) = \uparproj^{1/2} \cdot \isoParProj \isoParProj^\ast(p) = \uparproj^{1/2} \cdot p$ by \Cref{thm:AdjProjClosedRange}, \eqref{eq:LqProof-1} furthermore yields
 \begin{align}
   \cR(f_{\REF} + \isoParProj^\ast(p)) &=  \int_{\mR^2}  \int_\OmegaX |\isoParProj^\ast(p)(\bx_\perp, z)|^q \, \D z  \nnl
   &=  \int_{\mR^2} \uparproj(\bx_\perp)^{1-q/2} \left|p(\bx_\perp)\right|^q \, \D \bx_\perp = \bnorm{ \uparproj ^{1/q-1/2} \cdot p }_{L^q}^q.
 \end{align}
 
 Now let $q > 1$ and $f_0\neq 0$. Then $x \mapsto |x|^q$ is strictly convex so that strict inequality holds in \eqref{eq:LqProof-0} whenever the integrand $z \mapsto |f(\bx_\perp, z) |$ with $f = \adj \isoParProj(p) + f_0$ is non-constant. The latter is clearly the case for some $\bx_\perp\in \mR^2$, as elements in $\Kern(\isoParProj)$ are,  on each ray $\OmegaX$, either non-constant or identically zero. Hence, strict inquality holds in \eqref{eq:LqProof-2} if $f_0 \neq 0$.
 \end{proof}

 \optspace
 \section{Poisson-noise-adapted data fidelity} \label{Appendix:PoiDataFid}
 In the following, it is shown that the log-likelihood for Poisson-noisy data given in \eqref{eq:KLDataFid} can be approximated by the functional in \eqref{eq:KLDataFidCont} if variations of the true data $g_j$ within the supports of the $\omega_i$ are negligible. Specifically, we assume that $g_j$ is ``constant enough'' in $\supp(\omega_i)$ such that it may be pulled in and out of the integrals, i.e.\ for some $x_{ji} \in \mD$
\begin{align}
 \ln \left( t\cM_i(g_j) \right)  &= \ln \left(\int_{\mD} t\omega_i g_j \, \D x \right) \approx \ln \left( g_j (x_{ji}) \int_{\mD} t\omega_i  \, \D x \right)  = \ln \left( g_j (x_{ji}) \right) + \tilde c_{ij} \nnl
 &\approx  \frac{\int_{\mD} \ln \left( g_j \right) t\omega_i \, \D x}{ \int_{\mD} t\omega_i \, \D x } + \tilde c_{ij} 
\end{align}
where the $\tilde c_{ij}$ are independent of $g_j$.
Inserting this approximation into \eqref{eq:KLDataFid} yields
\begin{align}
  \cS^{\textup{Poi}}  \left( g^{\textup{obs}}_j;  g_j \right) &\approx \sum_{i = 1}^\Mpx \int_{\mD}  t g_j \omega_i \, \D x  - g^{\textup{obs}}_{ji} \cdot \left(  \frac{\int_{\mD} \ln \left( g_j  \right) t \omega_i \, \D x}{ \int_{\mD} t\omega_i \, \D x } + \tilde c_{ij}  \right) - g^{\textup{obs}}_{ji} \ln( g^{\textup{obs}}_{ji})  \nnl
  &=   \int_{\mD} \left( \left( \sum_{i = 1}^\Mpx \omega_i \right) t g_j - \left( \sum_{i = 1}^\Mpx  \frac{  g^{\textup{obs}}_{ji} \omega_i }{ \int_{\mD} \omega_i \, \D x }  \right)  \cdot  \ln \left( t g_j  \right) \right) \, \D x + \tilde c \nnl 
  &= \int_{\mD} \kl \left(g^{\obs}_{j, \textup{cont}}(x) ; t g_j(x)  \right) \cdot \omega(x)  \, \D x + c  \label{eq:KLDataFidCont-Appendix} 
\end{align}
with $\omega(x) := \sum_{i = 1}^\Mpx \omega_i(x)$, $g^{\obs}_{j, \textup{cont}}(x):= \big( \sum_{i = 1}^\Mpx g^{\obs}_{ji} \omega_i(x) / \int_{\mD}  \omega_i \, \D x \big) /  \omega(x)$ and some $c$ that is independent of $g_j$.

\bibliographystyle{siamplain}
\bibliography{/home/simon/SharedObjects/literature.bib}

\begin{thebibliography}{10}

\bibitem{AndersenKak1984SART}
{\sc A.~H. Andersen and A.~C. Kak}, {\em Simultaneous algebraic reconstruction
  technique ({SART}): a superior implementation of the {ART} algorithm},
  Ultrasonic imaging, 6 (1984), pp.~81--94.

\bibitem{AndersenHansen2014_genARTProx}
{\sc M.~S. Andersen and P.~C. Hansen}, {\em Generalized row-action methods for
  tomographic imaging}, Numerical Algorithms, 67 (2014), pp.~121--144.

\bibitem{BarrettEtAl2004ArtifactsInCT}
{\sc J.~F. Barrett and N.~Keat}, {\em Artifacts in {CT}: recognition and
  avoidance}, Radiographics, 24 (2004), pp.~1679--1691.

\bibitem{Bartels2015}
{\sc M.~Bartels, M.~Krenkel, J.~Haber, R.~Wilke, and T.~Salditt}, {\em {X}-ray
  holographic imaging of hydrated biological cells in solution}, Physical
  review letters, 114 (2015), p.~048103.

\bibitem{Beck2009Teboulle_FISTA}
{\sc A.~Beck and M.~Teboulle}, {\em A fast iterative shrinkage-thresholding
  algorithm for linear inverse problems}, SIAM Journal on Imaging Sciences, 2
  (2009), pp.~183--202.

\bibitem{Bertsekas2011_IncrProxMethods}
{\sc D.~P. Bertsekas}, {\em Incremental gradient, subgradient, and proximal
  methods for convex optimization: A survey}, Optimization for Machine
  Learning, 2010 (2011), pp.~1--38.

\bibitem{Bertsekas2011_IncrProxMethods2}
{\sc D.~P. Bertsekas}, {\em Incremental proximal methods for large scale convex
  optimization}, Mathematical programming, 129 (2011), p.~163.

\bibitem{Bleichrodt2015_STTomo}
{\sc F.~Bleichrodt, T.~van Leeuwen, and K.~J. Batenburg}, {\em Robust artefact
  reduction in tomography using {S}tudent’s--t data fitting}, in The 13th
  International Meeting on Fully Three-Dimensional Image Reconstruction in
  Radiology and Nuclear Medicine, 2015, pp.~395--398.

\bibitem{Burger2006NewtonKaczmarz}
{\sc M.~Burger and B.~Kaltenbacher}, {\em Regularizing {N}ewton--{K}aczmarz
  methods for nonlinear ill-posed problems}, SIAM Journal on Numerical
  Analysis, 44 (2006), pp.~153--182.

\bibitem{Burger2016VariationalImageProcessing}
{\sc M.~Burger, A.~Sawatzky, and G.~Steidl}, {\em First order algorithms in
  variational image processing}, in Splitting Methods in Communication,
  Imaging, Science, and Engineering, Springer, 2016, pp.~345--407.

\bibitem{ChambolleEhrhardtEtAl2017StochasticPDHG}
{\sc A.~Chambolle, M.~J. Ehrhardt, P.~Richt{\'a}rik, and C.-B. Sch{\"o}nlieb},
  {\em Stochastic primal-dual hybrid gradient algorithm with arbitrary sampling
  and imaging applications}, SIAM Journal on Optimization, 28 (2018),
  pp.~2783--2808.

\bibitem{ChambollePock2011}
{\sc A.~Chambolle and T.~Pock}, {\em A first-order primal-dual algorithm for
  convex problems with applications to imaging}, Journal of Mathematical
  Imaging and Vision, 40 (2011), pp.~120--145.

\bibitem{Cloetens1999holotomography}
{\sc P.~Cloetens, W.~Ludwig, J.~Baruchel, D.~Van~Dyck, J.~Van~Landuyt,
  J.~Guigay, and M.~Schlenker}, {\em Holotomography: Quantitative phase
  tomography with micrometer resolution using hard synchrotron radiation
  {X}-rays}, Applied Physics Letters, 75 (1999), pp.~2912--2914.

\bibitem{CombettesEtAl2011ProximalSplitting}
{\sc P.~L. Combettes and J.-C. Pesquet}, {\em Proximal splitting methods in
  signal processing}, in Fixed-point algorithms for inverse problems in science
  and engineering, Springer, 2011, pp.~185--212.

\bibitem{Cormack1963CT}
{\sc A.~M. Cormack}, {\em Representation of a function by its line integrals,
  with some radiological applications}, Journal of Applied Physics, 34 (1963),
  pp.~2722--2727.

\bibitem{Cormack1963CTII}
{\sc A.~M. Cormack}, {\em Representation of a function by its line integrals,
  with some radiological applications. ii}, Journal of Applied Physics, 35
  (1964), pp.~2908--2913.

\bibitem{CezaroEtAl2011_TikhonovKaczmarz}
{\sc A.~De~Cezaro, J.~Baumeister, and A.~Leitao}, {\em Modified iterated
  {T}ikhonov methods for solving systems of nonlinear ill-posed equations},
  Inverse Problems and Imaging, 5 (2011), pp.~1--17.

\bibitem{DeMan2001_PolyCTModel}
{\sc B.~De~Man, J.~Nuyts, P.~Dupont, G.~Marchal, and P.~Suetens}, {\em An
  iterative maximum-likelihood polychromatic algorithm for {CT}}, IEEE
  Transactions on Medical Imaging, 20 (2001), pp.~999--1008.

\bibitem{DefriseEtAl2011_TVSurrogateSART}
{\sc M.~Defrise, C.~Vanhove, and X.~Liu}, {\em An algorithm for total variation
  regularization in high-dimensional linear problems}, Inverse Problems, 27
  (2011), p.~065002.

\bibitem{EkebergEtAl2015_Virus3DCDI}
{\sc T.~Ekeberg, M.~Svenda, C.~Abergel, F.~R. Maia, V.~Seltzer, J.-M. Claverie,
  M.~Hantke, O.~J{\"o}nsson, C.~Nettelblad, G.~Van Der~Schot, et~al.}, {\em
  Three-dimensional reconstruction of the giant mimivirus particle with an
  {X}-ray free-electron laser}, Physical review letters, 114 (2015), p.~098102.

\bibitem{ElfvingHansen2014SemiconvKaczmarz}
{\sc T.~Elfving, P.~C. Hansen, and T.~Nikazad}, {\em Semi-convergence
  properties of {K}aczmarz’s method}, Inverse Problems, 30 (2014), p.~055007.

\bibitem{ElfvingNikazad2009KaczmarzAsPrecondLandweber}
{\sc T.~Elfving and T.~Nikazad}, {\em Properties of a class of block-iterative
  methods}, Inverse Problems, 25 (2009), p.~115011.

\bibitem{Hanke1996Regularization}
{\sc H.~W. Engl, M.~Hanke, and A.~Neubauer}, {\em Regularization of inverse
  problems}, Springer, 1996.

\bibitem{FeldkampDavisKress1984_FDKAlgorithm}
{\sc L.~Feldkamp, L.~Davis, and J.~Kress}, {\em Practical cone-beam algorithm},
  JOSA A, 1 (1984), pp.~612--619.

\bibitem{Gabay1983ADMM}
{\sc D.~Gabay}, {\em Chapter {IX}: Applications of the method of multipliers to
  variational inequalities}, in Augmented Lagrangian Methods: Applications to
  the Numerical Solution of Boundary-Value Problems, Elsevier, 1983, pp.~299 --
  331.

\bibitem{Gordon1970ART}
{\sc R.~Gordon, R.~Bender, and G.~T. Herman}, {\em Algebraic reconstruction
  techniques {(ART)} for three-dimensional electron microscopy and {X}-ray
  photography}, Journal of theoretical Biology, 29 (1970), pp.~471--481.

\bibitem{GuanGordon1994_MLS}
{\sc H.~Guan and R.~Gordon}, {\em A projection access order for speedy
  convergence of {ART} (algebraic reconstruction technique): a multilevel
  scheme for computed tomography}, Physics in Medicine and Biology, 39 (1994),
  p.~2005.

\bibitem{Guizar2015_QuantitativeROITomo}
{\sc M.~Guizar-Sicairos, J.~J. Boon, K.~Mader, A.~Diaz, A.~Menzel, and
  O.~Bunk}, {\em Quantitative interior {X}-ray nanotomography by a hybrid
  imaging technique}, Optica, 2 (2015), pp.~259--266.

\bibitem{HohageWerner2016PoissonReview}
{\sc T.~Hohage and F.~Werner}, {\em Inverse problems with poisson data:
  statistical regularization theory, applications and algorithms}, Inverse
  Problems, 32 (2016), p.~093001.

\bibitem{Hounsfield1973CT}
{\sc G.~N. Hounsfield}, {\em Computerized transverse axial scanning
  (tomography): Part~1. description of system}, The British journal of
  radiology, 46 (1973), pp.~1016--1022.

\bibitem{Humphries2014_PolyCTModel}
{\sc T.~Humphries and A.~Faidani}, {\em Segmentation-free quasi-{N}ewton method
  for polyenergetic {CT} reconstruction}, in 2014 IEEE Nuclear Science
  Symposium and Medical Imaging Conference (NSS/MIC), IEEE, 2014, pp.~1--5.

\bibitem{Johansen2015GammaTomo}
{\sc G.~Johansen}, {\em Gamma-ray tomography}, in Industrial Tomography,
  Elsevier, 2015, pp.~197--222.

\bibitem{Kaczmarz1937ART}
{\sc S.~Kaczmarz}, {\em Angen{\"a}herte {A}ufl{\"o}sung von {S}ystemen linearer
  {G}leichungen}, Bulletin International de {l'}Academie Polonaise des Sciences
  et des Lettres, 35 (1937), pp.~355--357.

\bibitem{KindermannLeitao2014_genARTConvergence}
{\sc S.~Kindermann and A.~Leitao}, {\em Convergence rates for {K}aczmarz-type
  regularization methods}, Inverse Problems and Imaging, 8 (2014),
  pp.~149--172.

\bibitem{Kostenko2013_AllAtOncePCTWithTV}
{\sc A.~Kostenko, K.~J. Batenburg, A.~King, S.~E. Offerman, and L.~J. van
  Vliet}, {\em Total variation minimization approach in in-line {X}-ray
  phase-contrast tomography}, Optics express, 21 (2013), pp.~12185--12196.

\bibitem{Krenkel2014BCAandCTF}
{\sc M.~Krenkel, M.~T{\"o}pperwien, M.~Bartels, P.~Lingor, D.~Schild, and
  T.~Salditt}, {\em {X}-ray phase contrast tomography from whole organ down to
  single cells}, in Developments in {X}-ray Tomography IX, vol.~9212,
  International Society for Optics and Photonics, 2014, p.~92120R.

\bibitem{Louis2016_ExactConeReconFormula}
{\sc A.~K. Louis}, {\em Exact cone beam reconstruction formulae for functions
  and their gradients for spherical and flat detectors}, Inverse Problems, 32
  (2016), p.~115005.

\bibitem{MaretzkeEtAl2016OptExpr}
{\sc S.~Maretzke, M.~Bartels, M.~Krenkel, T.~Salditt, and T.~Hohage}, {\em
  Regularized {N}ewton methods for {X}-ray phase contrast and general imaging
  problems}, Optics express, 24 (2016), pp.~6490--6506.

\bibitem{MidgleyWeyland2002_ETomoContrastMechanisms}
{\sc P.~Midgley and M.~Weyland}, {\em 3d electron microscopy in the physical
  sciences: the development of {Z}-contrast and {EFTEM} tomography},
  Ultramicroscopy, 96 (2003), pp.~413--431.

\bibitem{Mori2013PhotonStarvation}
{\sc I.~Mori, Y.~Machida, M.~Osanai, and K.~Iinuma}, {\em Photon starvation
  artifacts of {X}-ray {CT}: their true cause and a solution}, Radiological
  physics and technology, 6 (2013), pp.~130--141.

\bibitem{MorishimaEtAl2017PyramidMyonCT}
{\sc K.~Morishima, M.~Kuno, A.~Nishio, N.~Kitagawa, Y.~Manabe, M.~Moto,
  F.~Takasaki, H.~Fujii, K.~Satoh, H.~Kodama, et~al.}, {\em Discovery of a big
  void in khufu’s pyramid by observation of cosmic-ray muons}, Nature, 552
  (2017), p.~386.

\bibitem{Natterer}
{\sc F.~Natterer}, {\em {The mathematics of computerized tomography}}, vol.~32
  of {Classics of Applied Mathematics}, Society for Industrial and Applied
  Mathematics, 2001.

\bibitem{Oektem2015MathETomo}
{\sc O.~{\"O}ktem}, {\em Mathematics of electron tomography}, in Handbook of
  Mathematical Methods in Imaging, Springer, 2015, pp.~937--1031.

\bibitem{Ruhlandt2016_RadonPCICommute}
{\sc A.~Ruhlandt and T.~Salditt}, {\em Three-dimensional propagation in
  near-field tomographic {X}-ray phase retrieval}, Acta Crystallographica
  Section A: Foundations and Advances, 72 (2016), pp.~215--221.

\bibitem{Salditt2015GINIX}
{\sc T.~Salditt, M.~Osterhoff, M.~Krenkel, R.~N. Wilke, M.~Priebe, M.~Bartels,
  S.~Kalbfleisch, and M.~Sprung}, {\em Compound focusing mirror and {X}-ray
  waveguide optics for coherent imaging and nano-diffraction}, Journal of
  synchrotron radiation, 22 (2015), pp.~867--878.

\bibitem{Setzer2011ADMMForImageProcessing}
{\sc S.~Setzer}, {\em Operator splittings, {B}regman methods and frame
  shrinkage in image processing}, International Journal of Computer Vision, 92
  (2011), pp.~265--280.

\bibitem{SidkyEtAl2012TomoWithChambollePock}
{\sc E.~Y. Sidky, J.~H. J{\o}rgensen, and X.~Pan}, {\em Convex optimization
  problem prototyping for image reconstruction in computed tomography with the
  {Chambolle}--{Pock} algorithm}, Physics in Medicine and Biology, 57 (2012),
  p.~3065.

\bibitem{SidkyEtAl2008CBTomoWithTV}
{\sc E.~Y. Sidky and X.~Pan}, {\em Image reconstruction in circular cone-beam
  computed tomography by constrained, total-variation minimization}, Physics in
  Medicine and Biology, 53 (2008), p.~4777.

\bibitem{ASTRAToolbox2016}
{\sc W.~van Aarle, W.~J. Palenstijn, J.~Cant, E.~Janssens, F.~Bleichrodt,
  A.~Dabravolski, J.~De~Beenhouwer, K.~J. Batenburg, and J.~Sijbers}, {\em Fast
  and flexible {X}-ray tomography using the {ASTRA} toolbox}, Optics express,
  24 (2016), pp.~25129--25147.

\bibitem{ASTRAToolbox2015}
{\sc W.~van Aarle, W.~J. Palenstijn, J.~De~Beenhouwer, T.~Altantzis, S.~Bals,
  K.~J. Batenburg, and J.~Sijbers}, {\em The {ASTRA} toolbox: A platform for
  advanced algorithm development in electron tomography}, Ultramicroscopy, 157
  (2015), pp.~35--47.

\bibitem{XuMueller2006_TomoInterpMethods}
{\sc F.~Xu and K.~Mueller}, {\em A comparative study of popular interpolation
  and integration methods for use in computed tomography}, in 3rd IEEE
  International Symposium on Biomedical Imaging: Nano to Macro, 2006,
  pp.~1252--1255.

\bibitem{YangEtAl2010_ROIByGenTV}
{\sc J.~Yang, H.~Yu, M.~Jiang, and G.~Wang}, {\em High-order total variation
  minimization for interior tomography}, Inverse Problems, 26 (2010),
  p.~035013.

\end{thebibliography}

\end{document}